\input amstex

\documentstyle{amsppt}

\magnification1200

\tolerance=10000

\def\n#1{\Bbb #1}

\def\Tr{\hbox{Tr }}

\def\im{\hbox{im }}

\def\re{\hbox{re }}

\def\Gal{\hbox{Gal }}

\def\Ker{\hbox{Ker }}

\def\opr{\hbox{def}}

\def\Res{\hbox{Res }}

\def\ord{\hbox{ord }}

\def\rank{\hbox{ rank }}

\def\diag{\hbox{ diag }}

\def\rd{\sqrt{-\Delta}}

\def\d{\Delta}

\def\s{\sigma}

\def\e11{E_{11}}

\def\ve{\varepsilon}

\def\lw{\Longrightarrow }

\topmatter
\title
Kolyvagin's trace relations for Siegel sixfolds
\endtitle
\author
D. Logachev
\endauthor
\NoRunningHeads
\date July, 2005 \enddate
\address
Departamento de Matem\'aticas
Universidad Sim\'on Bol\'\i var
Apartado Postal 89000
Caracas, Venezuela
\endaddress
\thanks The author is grateful to A. Zelevinskij who indicated him the proof of the Lemma
4.1.12, and to Th. Berry for linguistic correstions.  \endthanks
\abstract In his earlier preprints the author offered a program of generalization of
Kolyvagin's result of finiteness of SH to the case of some motives which are quotients of
cohomology motives of Shimura and Drinfeld varieties. The present paper is devoted to the
first step of this program --- finding of an analog of Kolyvagin's trace relations. We
solve it for Siegel sixfolds and for the Hecke correspondences related to the matrices
diag($1,1,1,p,p,p$) and diag($p,1,1,p,p^2,p^2$). This is the first non-trivial case for
Shimura varieties. Some results for other types of Siegel varieties and Hecke
correspondences are obtained.

Ideas and methods of the present paper open a large new area of research: results given
here constitute a tiny part of what can be done. Particularly, maybe it is possible to
realise for Drinfeld varieties of any even rank the program of generalization of
Kolyvagin's result.
\endabstract
\endtopmatter
\document
{\bf Contents.}

\medskip

\settabs 20 \columns

\+ 1. Introduction and detailed description of results. &&&&&&&&&&&&&&&&&&&  \cr

\medskip

\+ 1.1. Statement of result, and justification of the subject. &&&&&&&&&&&&&&&&&&& 2 \cr

\+ 1.2. Kolyvagin's trace relations.&&&&&&&&&&&&&&&&&&& 4 \cr

\+ 1.3. Definitions related to Hecke correspondences.&&&&&&&&&&&&&&&&&&&  5\cr

\+ 1.4. Definitions related to partitions.&&&&&&&&&&&&&&&&&&&  6\cr

\+ 1.5. Contents of Section 3: Case $\goth T_p=T_p$.&&&&&&&&&&&&&&&&&&&  7\cr

\+ 1.6. Contents of Section 4: Case $\goth T_p=T_{p,\goth i}$.&&&&&&&&&&&&&&&&&&&  9\cr

\+ 1.7. Contents of Section 5: Miscellaneous.&&&&&&&&&&&&&&&&&&&  11\cr

\+ 1.8. Conjectures and possibilities of further investigation.&&&&&&&&&&&&&&&&&&&  11\cr

\medskip

\+ Section 2. D-equivalence. &&&&&&&&&&&&&&&&&&&  12\cr

\medskip

\+ Section 3. Case of Hecke correspondence $T_p$. &&&&&&&&&&&&&&&&&&&  \cr

\medskip

\+ 3.1. Definition of $V$. &&&&&&&&&&&&&&&&&&&  16\cr

\+ 3.2. Preliminary notations and lemmas. &&&&&&&&&&&&&&&&&&&  16\cr

\+ 3.3. Finding of D-equivalence and of the field of definition. &&&&&&&&&&&&&&&&&&& 
18\cr

\+ 3.4. Finding of the degree of the second projection. &&&&&&&&&&&&&&&&&&&  22\cr

\medskip

\+ Section 4. Case of Hecke correspondence $T_{p,\goth i}$. &&&&&&&&&&&&&&&&&&&  \cr

\medskip

\+ 4.1. Some matrix equalities. &&&&&&&&&&&&&&&&&&&  27\cr

\+ 4.2. Finding of D-equivalence and of the field of definition.&&&&&&&&&&&&&&&&&&&  31\cr

\medskip

\+ Section 5. Miscellaneous.&&&&&&&&&&&&&&&&&&&  \cr

\medskip

\+ 5.1. Structure of $T_{p,2}(V)$.&&&&&&&&&&&&&&&&&&&  37\cr

\+ 5.2. Action of $T_p$ on components of $T_{p}(V)$.&&&&&&&&&&&&&&&&&&&  40\cr

\+ 5.3. A theorem on a weak ``equivalence'' relation.&&&&&&&&&&&&&&&&&&&  42\cr

\+ 5.4. Non-coincidence of components of $T_{p,\goth i}(V)$.&&&&&&&&&&&&&&&&&&&  47\cr

\medskip

\medskip

{\bf 1. Introduction.} \nopagebreak

\medskip

{\bf 1.1. Statement of result, and justification of the subject.} \nopagebreak

\medskip

 Kolyvagin ([K89] and subsequent papers) proved finiteness of Tate-Shafarevich group and
of group of rational points of modular\footnotemark \footnotetext{After 1999 we can omit
the word ``modular'' but in 1989 the Taniyama--Shimura conjecture was not proved yet.
Modularity of the elliptic curve is used essentially in Kolyvagin's proof.} elliptic
curves over $\n Q$ of analytic rank 0. There is a natural problem

\medskip

{\bf (*)} To generalize this result to the case of some submotives of high-dimensional
Shimura varieties and/or Drinfeld modular varieties.

\medskip

The present paper is devoted to generalization of the first step of Kolyvagin's proof ---
trace relations for Heegner poins (see (1.2.1)) --- to the case of Siegel varieties of
genus 3 (some results are valid for any $g$).

\medskip

{\bf Remark.} The author offers in preprint [L04.2] a program of proof of (*) for Shimura
varieties, he solves in this paper some further steps of the proof and indicates the
obtacles. The results of the present paper and of [L04.2] give some evidence that maybe
we shall be able to extend Kolyvagin's result to the functional case.

\medskip

The subject of the present paper is the following. Let $X$ be an irreducible component of
a Shimura variety of a fixed level, $V\subset X$ its Shimura subvariety and $\goth T_p$ a
$p$-Hecke correspondence on $X$ ($p$ is a prime fixed throughout the paper). We get some
information on the structure of $\goth T_p(V)$ (which a priori is a cycle on $X$) for the
case when $X$ is a Siegel variety of genus $g$, $V$ its subvariety corresponding to a
reductive group $GU(\goth r, \goth s)$, $\goth r+\goth s=g$ (points of $V$ parametrize
abelian $g$-folds with multiplication by an imaginary quadratic field $K$) and $\goth
T_p=T_p$ or $T_{p,\goth i}$ (see (1.3.1) for the definition). Namely, we find (roughly
speaking) the set of irreducible components of $\goth T_p(V)$, their fields of definition
and Galois action on them. The results are complete for the case $\goth T_p=T_p$, $g=3$
and near to complete for any $g$. For the case $\goth T_p=T_{p,1}$ only the general
components are described completely.

\medskip

The answers are given in terms of geometry of the finite set $\goth T_p(t)$, where $t\in
V$ is a generic point. $\goth T_p(t)$ is some kind of Grassmann variety, see (1.3.4),
(1.3.6). We introduce 2 partitions D and I of ( = equivalence relations on) $\goth
T_p(t)$: roughly speaking, two points of $\goth T_p(t)$ are equivalent iff they belong to
the same irreducible component of $\goth T_p(V)$, see 1.4, 2.8 for the exact definitions.
We get a simple description of these partitions D and I in terms of geometry of $\goth
T_p(t)$. Namely, we introduce in 1.5 elementary ``geometric'' partitions $\goth D$,
$\goth H$ on $\goth T_p(t)$, and we describe D, I and the Galois actions on the set of
parts of I in terms of $\goth D$, $\goth H$. Results for the case $\goth T_p=T_p$ are
given in subsection 1.5 (theorems 3.3.1, 3.3.3, 1.5.3, conjectures 1.5.4, 1.5.6, 1.5.7),
and for the case $\goth T_p=T_{p,\goth i}$ in subsection 1.6 (conjecture 1.6.3, theorems
4.2.15, 4.2.19, conjecture 4.2.20).

\medskip

The structure of calculations shows that we have a deep and beautiful theory. It is easy
to get an algebraic criterion of equivalence of points of $\goth T_p(t)$ (Theorem 2.16;
see also (3.3.1.3) for its explicit form), but it is not clear how to solve the
corresponding equations. Paper [L01] contains a solution of this problem for Siegel
varieties of genus 2. Its method is based on explicit formulas for the image of $GU(1,1)$
in $GSp_{4}$. These formulas do not exist for $g>2$.

\medskip

Nevertheless, application of the crucial formula (3.3.1.4) permits us to solve the
problem for the simplest Hecke correspondence $\goth T_p=T_p$ uniformly for all $g$ (for
a generic equivalence class) much easier than it was made in [L01] for $g=2$.

\medskip

Situation for other Hecke correspondences $\goth T_p=T_{p,\goth i}$ is more complicated.
Formulas (4.2.3) (analogs of (3.3.1.4) for $T_{p,\goth i}$) permit us to reduce the
problem of finding of a generic equivalence class to solving of equations (4.2.15.10 --
4.2.15.12) which look rather hopeless from the first sight.

\medskip

But here we get the second wonder: the equations can be greatly simplified!
The process of their simplification reminds assembling a puzzle when
different parts unexpectedly and beautifully fit together.
For example, sometimes in the process of proof it was seen that if some polynomial $P$
belonged to an ideal $\goth I$ of $\n Q[X_1, ..., X_n]$ then the calculations would be
much simpler. Since in this situation nature is benevolent for us, I concluded that $P$
{\it must belong} to $\goth I$. It was absolutely not seen beforehand. I used a computer
program; yes, computer shows that $P \in \goth I$. Later it was possible to prove
(without computer) that really in all cases $P \in \goth I$.

\medskip

The reader can tell that if a problem has a simple answer then formulas that appear in
its solution do simplify. But

\medskip

{\it We do not know beforehand that our problem has a simple solution!}

\medskip

Really, when we consider a problem and try to generalize it (increasing dimension of
objects for example) it can happen that formulas become more and more complicated, and no
general idea is seen. So, the most important information obtained in the present paper is
the following: the problem of finding of generic components of $\goth T_p(V)$ is solvable
(to find a solution --- if we know that it exists --- is an easier task).

\medskip

And what happens with non-generic components? I do not know. For the case $\goth T_p=T_p$
there is Conjecture 1.5.4 which gives us the complete answer. Without doubt, for any
fixed genus $g$ and the rank $j$ of some auxiliary matrix, there exists a proof of
conjecture 1.5.4 similar to the proof of the theorem 3.3.4 (case $g=3$, $j=2$). But does
exist a uniform proof for all $g$, $j$?

\medskip

If $\goth T_p=T_{p,\goth i}$ then the situation is worse. Practically, for non-generic
components of $\goth T_p(V)$ I can work only ``modulo $p$'' but not ``modulo $p^2$'', see
Section 5.3. Does a simple description of equivalence classes for this case exist? If
not, what is possible to tell about asymptotics of the quantity of irreducible
components, etc.?

\medskip

It is interesting that situations for different Hecke correspondences are quite
unsimilar. Because of applications to (*), we are interested to consider irreducible
components of $\goth T_p(V)$ which are defined over the $p$-ring class field of $K$ (good
components) while components which are defined over $K$ (bad components) generate
obstacles. For the case $\goth T_p=T_p$ the generic components are good and the special
ones are bad, for $\goth T_p=T_{p,1}$ ($g=3$) the situation is (roughly speaking)
inverse, and for $\goth T_p=T_{p,2}$ there is no good components at all. What is a
non-formal explanation of this phenomenon?

\medskip

Finally, the subjects of [L04.1], [L04.2] and of the present paper open a large area of
investigation, see introduction to [L04.2] for a discussion of possibilities of
investigation of high level, and Subsection 1.8 of the present paper. It would be very
important to generalize these results to the functional case.

\medskip

The rest of Section 1 contains definitions (subsections 1.3, 1.4) and statements of the
main theorems and conjectures (subsections 1.5, 1.6).

\medskip

{\bf 1.2. Kolyvagin's trace relations.} \nopagebreak

\medskip

Let us give definitions for the original Kolyvagin's case. Let $N$ be a level,
$X_0(N)=\overline{\Gamma_0(N)\backslash\Cal H}$ the compactification of the modular curve
of level $N$, $K=\n Q(\rd)$ an imaginary quadratic field. In order to simplify the
notations and proofs, we consider only the case when $h(K)=1$, although this restriction
is not essential. Sometimes when we use notation $K^1$ --- the Hilbert class field of $K$
--- this means that the corresponding result is valid for any $K$.

Recall the definition of Heegner point (see for example [GZ86], [K89] for the details).
Points $t$ of the open part of $X_0(N)$ are in one-to-one correspondence with the
isogenies of elliptic curves $\psi_t: A_t \to A'_t$ such that $\Ker \psi_t=\n Z/N\n Z$. A
point $t\in X_0(N)$ is called a Heegner point with respect to $K$ if both $A_t$, $A'_t$
have complex multiplication by the same order of $K$. A Heegner divisor is a Galois orbit
of a Heegner point; Heegner divisors are exactly 0-dimensional Shimura subvarieties of
$X_0(N)$ in the sense of Deligne ([D71]). If all prime factors of $N$ split in $K$ then
for a given $K$ there exists a ``principal'' Heegner point $x_1\in X_0(N)(K^1)$.

The main object of the present paper is a prime $p$. We restrict ourselves to the case
$p\ne 2$ is inert in $K$, i.e. $\left(\frac{-\d}p\right)=-1$. The case $p$ splits in $K$
is much more complicated technically, see for example [L01] where for $g=2$ it is treated
completely. Further, we consider only $p$ such that $p$ does not divide $N$. Recall that
$K^p$ --- the ring class field of $K$ of conductor $p$ --- is an abelian extension of $K$
corresponding to some subgroup in the idele group $I_K$ of $K$, see for example [K89] for
the exact definition. We have: $K^p/K$ is ramified only at $p$, and $\Gal(K^p/K^1)=\n
Z/(p+1)$.

There are 2 objects associated to $p$: a Heegner point $x_p\in X_0(N)(K^p)$ and the
$p$-Hecke correspondence $T_p$ on $X_0(N)$. We have a formula (equality of divisors on
$X$, case $h(K)=1$):

$$T_p(x_1)=\Tr_{K^p/K}(x_p)\eqno{(1.2.1)}$$

(Kolyvagin's trace relation for $X_0(N)$ ([K89])).

We can consider (1.2.1) as a description of the action of $T_p$ on $x_1$ where $x_1$ is
an irreducible component of a 0-dimensional Shimura subvariety of $X$. So, for a triple
$X$, $V$, $\goth T_p$ of (1.1) a high-dimensional analog of (1.2.1) is the description of
the cycle $\goth T_p(V)$. Namely, this cycle is a formal sum of irreducible components of
some other Shimura subvarieties of $X$. The problem is to describe the structure of
$\goth T_p(V)$, particularly its representation as a formal sum of Shimura subvarieties,
their irreducible components and Galois action on them.

\medskip

{\bf 1.3. Definitions related to Hecke correspondences.} \nopagebreak

\medskip

We consider the case $X$ is a Siegel variety of genus $g$ of any level $N$, so the
corresponding algebraic group $G$ is $GSp_{2g}(\n Q)$. The algebra of $p$-Hecke
correspondences on $X$ is the ring of polynomials with $g$ generators denoted by $T_p$,
$T_{p,1}, \dots, T_{p,g-1}$. They are double cosets corresponding to the diagonal
matrices $\tau_p$, $\tau_{p,1}, \dots, \tau_{p,g-1}$ respectively, where

$$\tau_p=\diag(\underbrace{1,\dots,1}_{g\hbox{ times}}, \underbrace{p,\dots,p}_{g\hbox{
times}})$$

$$\hbox{and }\tau_{p,\goth i}=\diag(\underbrace{p,\dots,p}_{\goth i\hbox{ times}},
\underbrace{1,\dots,1}_{g-\goth i\hbox{ times}}, \underbrace{p,\dots,p}_{\goth i\hbox{
times}}, \underbrace{p^2,\dots,p^2}_{g-\goth i\hbox{ times}})\eqno{(1.3.1)}$$ $\goth i=0,
\dots, g$. (Caution: here and below notations of the present paper are other than in
[L04.1], [L04.2]). We have: $T_{p,g}$ is the trivial correspondence, and $T_{p,0}\in \n
Z[T_p, T_{p,1}, \dots, T_{p,g-1}]$.

Let us recall various interpretations of $T_p$, $T_{p,\goth i}$. Let us consider the
double cosets $\Gamma \tau_{p}\Gamma$, $\Gamma \tau_{p,\goth i}\Gamma$ (where
$\Gamma=GSp_{2g}(\n Z_p)$) and their decomposition as a union of ordinary cosets $$\Gamma
\tau_{p}\Gamma =\bigcup_{j\in \goth S(g)}\Gamma \sigma_j, \ \ \ \Gamma \tau_{p,\goth
i}\Gamma =\bigcup_{j\in \goth S(g,\goth i)}\Gamma \sigma_j\eqno{(1.3.2)}$$ where $\goth
S(g)$, $\goth S(g,\goth i)$ are the (abstract) sets of indices and $\sigma_j\in
GSp_{2g}(\n Q)\cap M_{2g}(\n Z)$ are representatives. We use for our calculations the
explicit sets of these representatives described in [L04.1], Section 2.3 for $\goth
T_p=T_p$ and in Section 2.4 for $\goth T_p=T_{p,\goth i}$. These sets will be denoted by
$S(g)$, $S(g,\goth i)$ respectively, i.e. $j \mapsto \sigma_j$ is an isomorphism $\goth
S(g)\to S(g)$, $\goth S(g,\goth i)\to S(g,\goth i)$.

For a generic $t\in X$ the set $T_{p}(t)$ (resp. $T_{p,\goth i}(t)$) is in 1 -- 1
correspondence with $\goth S(g)$ (resp. $\goth S(g,\goth i)$). Let $A_t$ be the abelian
$g$-fold corresponding to $t$. $(A_t)_p$ --- the group of $p$-torsion poins of $A_t$ ---
is a $2g$-dimensional vector space over $\n F_p$ endowed with a non-degenerated skew form
$B$ coming from a Riemann form on $A_t$.

\medskip

{\bf (1.3.3)} The set $T_{p}(t)$ (and hence $\goth S(g)$) is in 1 -- 1 correspondence
with the set of $B$-isotropic $g$-dimensional subspaces in $(A_t)_p$. We denote this
variety by $G_I(g,2g)(\n F_p)=G_I(g,(A_t)_p)$ ($G_I$ means the isotropic Grassmannian).
So, for the case $\goth T_{p}=T_{p}$ we have identifications

$$T_{p}(t)=\goth S(g)= S(g)= G_I(g,2g)(\n F_p)=G_I(g,(A_t)_p)\eqno{(1.3.4)}$$

For $\s\in S(g)$ we denote by

$$W_{\s}\subset (A_t)_p\eqno{(1.3.5)}$$

the corresponding isotropic $\n F_p$-subspace.

Analogously, for $t\in X$, $\goth T_{p}=T_{p,\goth i}$, the set $T_{p,\goth i}(t)$ is
isomorphic to a generalized Grassmannian $G_I(\goth i,g, 2g)$ over $\n F_p$ (see 1.6.1
for the definition); it is a scheme over $\n Z/p^2\n Z$. The analog of (1.3.4) is

$$T_{p,\goth i}(t)=\goth S(g,\goth i)= S(g,\goth i)=G_I(\goth i,g, 2g)\eqno{(1.3.6)}$$

{\bf 1.4. Definitions related to partitions.}  \nopagebreak

\medskip

We consider a generic point $t\in V$. In this case $\goth T_p(t)\subset \goth T_p(V)$,
and the structure of $\goth T_p(V)$ as a union of Shimura subvarieties defines various
partitions of ( = equivalence relations on) $\goth T_p(t)$ and isomorphic sets given in
(1.3.4), (1.3.6). Roughly speaking, 2 points of $\goth T_p(t)$ belong to the same
partition set ( = equivalent) iff they belong to the same Shimura subvariety. So, we must
describe these partitions. It turns out that (roughly speaking) partition sets are
subschemes of $\goth T_p(t)$.

The problem of finding partitions on $\goth T_p(t)$ can be formulated in terms of
algebraic groups $G$, $G_V$ over $\n Q$ associated to $X$, $V$ respectively. Definitions
of 3 different types of partitions are given in [L01], Section 4 in terms of equivalence
relations on $\goth T_p(t)$ ($S$ of [L01] is $\goth T_p(t)$). We shall treat in the
present paper only 2 partitions:

\medskip

(a) the strong equivalence of [L01]: 2 points of $\goth T_p(t)$ are strongly equivalent
($\iff$ belong to the same partition set of strong partition) iff they belong to the same
irreducible component of $\goth T_p(V)$. See [L01], 4.1 for the explicit formulas. For
brevity, we call strong partition by I-partition (I of irreducible).

\medskip

(b) D-equivalence: 2 points of $\goth T_p(t)$ are D-equivalent ($\iff$ belong to the same
partition set of D-partition) iff (roughly speaking) there exists a Shimura subvariety
(in the sence of [D71]) of $\goth T_p(V)$ containing both these points (see 2.8 for the
exact definition). This partition was not considered in [L01].

\medskip

Theorem 2.16 gives us a criterion of D-equivalence.

\medskip

{\bf (1.4.1)} We introduce some notations related to I- (resp. D-)partition of $S(g)$
(resp. $S(g,\goth i)$). To avoid overuse of the word ``respectively'', we give these
notations only for the Hecke correspondence $T_p$; notations for $T_{p, \goth i}$ are
similar. The set of parts of the partition will be denoted by $\Cal S_I(g)$ (resp. $\Cal
S_D(g)$), it is the quotient set of $S(g)$ by the I- (resp. D-)equivalence relation. For
$\s\in S(g)$ we denote the corresponding element of the quotient set by $I(\s)$ (resp.
$D(\s)$). For $k \in \Cal S_I(g)$ (resp. $k \in \Cal S_D(g)$) we denote the corresponding
subset of $S(g)$ by $S(g)_k$ (i.e. $S(g)_{I(\s)}$ (resp. $S(g)_{D(\s)}$) is the set of
elements of $S(g)$ which are I- (resp. D-)equivalent to $\s$). The irreducible component
of $T_p(V)$ (resp. the Shimura subvariety of $T_p(V)$, see Section 2 for the exact
description of this object) that corresponds to $k$ will be denoted by $T_p(V)_k$.
Particularly, for $\s\in S(g)$ the irreducible component of $T_p(V)$ (resp. the Shimura
subvariety of $T_p(V)$) that corresponds to $\s$ will be denoted by $T_p(V)_{I(\s)}$
(resp. $T_p(V)_{D(\s)})$.

\medskip

{\bf (1.4.2)} Clearly I-partition is stronger than D-partition. For $k\in\Cal S_D(g)$ we
denote by $k_I$ the set of parts of I-partition that are contained in $S(g)_k$. Let
$k_1\in k_I \subset \Cal S_I(g)$ be an element, $L_{k_1}$ the field of definition of
$T_p(V)_{k_1}$ and $L_k$ the reflex field of $T_p(V)_k$. According [D71], $\Gal(L_k)$
acts simply transitively on $k_I$ (this action comes from the action of $\Gal(L_k)$ on
the set of irreducible components of $T_p(V)_k$), and $\Gal(L_{k_1})$ is the stabiliser
of $k_1$.

\medskip

We shall study mainly the D-partition, and only 3 conjectures (1.5.6, 1.5.7 and 4.2.20)
are concerned to I-partition and the above Galois action.

\medskip

Because of applications to [L04.2], we shall be interested mainly by the case of $T_p$
for any $g$ and of $T_{p,1}$ for $g=3$.

\medskip


\medskip

{\bf 1.5. Contents of Section 3: Case $\goth T_p=T_p$.}  \nopagebreak

\medskip

Here we formulate theorems and conjectures for the case $\goth T_p=T_p$. We introduce a
$\goth D$-partition of the sets of (1.3.4) ($\goth D$ because of dimension) which has a
simple geometric definition, and we describe relations between $\goth D$-partition and I-
and D-partitions.

The set of parts of $\goth D$-partition is indexed by even numbers $j$, $g\le j \le 2g$.
The definition of $\goth D_j$ --- the $j$-th part of $\goth D$-partition of $S(g)=
G_I(g,(A_t)_p)$ is the following. Recall that $V$ is a Shimura variety of PEL-type
parametrizing abelian $g$-folds with multiplication by $K$. The exact definition of $V$
is given in Subsection 3.1. Since $p$ inert in $K$, we have $O_K/p=\n F_{p^2}$. Let $t$
be a generic point of $V$. Action of $O_K$ on $A_t$ endows $(A_t)_p$ by the structure of
$\n F_{p^2}$-module. Let $\s\in S(g)$ be any element and $W_{\s}\subset (A_t)_p$ from
1.3.5. $\n F_{p^2} W_{\s}$ is a $\n F_{p^2}$-subspace of $(A_t)_p$ whose $\n
F_{p}$-dimension can take any even value from $g$ to $2g$. We define

$$\goth D_j=\{\s\in S(g)\ | \ \dim_{\n F_{p}}(\n F_{p^2} W_{\s})=j\}\eqno{(1.5.1)}$$

Most of below theorems are proved only on the open part of $G_I(g,(A_t)_p)$ (see (3.2.1))
denoted by $G_I(g,(A_t)_p)^{open}$; the corresponding open subset of $S(g)$ is denoted by
$S(g)^{open}$. The exact statement of what is proved is given in the body of the paper;
here we give the statements of these theorems in their complete form. We mark these
theorems with (*) and we leave the number of this theorem in the body of the paper.
Further, we give here statements of some theorems whose proof is not given at all,
because it is elementary. We mark these theorems with (**).

\medskip

{\bf Remark 1.5.2.} Since the below theorems are not proved completely, formally they
should be formulated as conjectures. {\it Absolute} evidence that they are true ({\it and
moreover that their complete proof is a routine work}) comes from consideration of some
particular cases. For example, I made complete calculations for the case $g=2$; they are
so large and so elementary, that there were no meaning to include them in the text of
[L01] (see Proposition 6.3.7, (d)).

\medskip

So, the problem of complete proof of the theorems (*) and (**) can be considered as an
exercise for those who will continue these investigation. See Remark 3.3.2 and 1.8, (1).

\medskip

Let $\s\in S(g)^{open}$. We can associate it a matrix $T(\s)\in M_g(\n Z)$ (see 3.2.2)
and $T(\tilde \s)\in M_g(\n F_p)$ --- its reduction modulo $p$. It has the following
property:

\medskip

{\bf (*) Lemma 3.2.3.} If $\s\in \goth D_{j}$ then rank $T(\tilde \s)=j-g$.

\medskip

{\bf (*) Theorems 3.3.1, 3.3.3.} If $g$ is odd then $\goth D_{2g}$ is a D-partition set.
If $g$ is even $\ne 2$ then $\goth D_{2g}$ is a union of two D-partition sets: 2 elements
$\sigma_i$, $\sigma_j\in\goth D_{2g}$ are D-equivalent iff the ratio $\det T(\tilde
\sigma_i)/ \det T(\tilde \sigma_j)$ is a square in $\n F_p^*$.

\medskip

{\bf Theorem 3.3.5.} If $g$ is odd then all components of $T_p(V)$ do not coincide with
$V$ itself.

\medskip

{\bf (**) Theorem 1.5.3.} If $g$ is even then $\goth D_{g}$ is a part of D-partition, and
moreover it is the only part such that the corresponding component of $T_p(V)$ is $V$
itself.

\medskip

For $j\ne g, 2g$ we have only a

\medskip

{\bf Conjecture 1.5.4.} D-partition of $S(g)$ is a subpartition of its $\goth
D$-partition. 2 elements $\sigma_i$, $\sigma_j\in\goth D_{j}$ are D-equivalent iff
quadratic forms over $\n F_p$ defined by the matrices $T(\tilde \sigma_i)$, $T(\tilde
\sigma_j)$ are isomorphic up to multiplication by a scalar.

\medskip
{\bf Remark.} The theory of quadratic forms over finite fields shows that Conjecture
1.5.4 implies that if $g$ is odd then D-partition coincide with $\goth D$-partition,
while if $g$ is even $\ne 2$ then each $\goth D_{j}$ consists of two parts of D-partition.
\medskip
Parts of D-partition that are contained in $\goth D_{2g}$ are called the good parts, the
only part consisting of $\goth D_{g}$ is called the trivial part, and parts that are
contained in $\goth D_{j}$, $j\ne g, 2g$, are called the bad parts. The corresponding
components of $T_p(V)$ have the same names.

\medskip

{\bf (*) Theorem 3.3.4.} Conjecture 1.5.4 is true for $g=3$, $(\goth r, \goth s) = (2,1)$.

\medskip

{\bf (*) Corollary.} For $g=3$, $(\goth r, \goth s) = (2,1)$ there exist 1 good component
of $T_p(V)$, 1 bad component and no trivial components (i.e. equal to $V$ itself).

\medskip

{\bf (*) Theorem 3.3.6.} For any $g$, any $\goth r\ne \goth s$ the field of definition of
any irreducible component of a good component of $T_p(V)$ is $K^p$.

\medskip

{\bf (**) Theorem 1.5.5.} For all $\goth r, \goth s$ the restriction of I-partition on
all bad parts of D-partition is trivial (i.e. consists of one part). If $\goth r=\goth s$
then the same is true for good parts as well. The field of definition of all mentioned
components is $K^1$.

\medskip

Now let us formulate conjectures on I-partition. If an irreducible component of $T_p(V)$
is defined over $K^1$ then the corresponding part of I-partition is also a part of
D-partition. So, we need only to describe I-partition on the set $\goth D_{2g}$.

There exists the Pl\"ucker embedding $\goth P$: $G_I(g,2g)(\n F_p) \hookrightarrow
G(g,2g)(\n F_p) \hookrightarrow P^m(\n F_p)$ where $m=\left(\matrix 2g\\g\endmatrix
\right)-1$. $\goth P(\goth D_i)\subset\goth P(S(g))\subset\goth P(G_I(g,2g)(\n
F_p))\subset P^m$ are some kind of determinantal varieties in $P^m$.

\medskip

{\bf Remark 3.2.5.} $\goth P(\cup_{j=g}^{2g-2}\goth D_j)$ --- the union of all bad parts
(and the trivial part for even $g$) --- is the intersection of $\goth P(G_I(g,2g))$ with
a codimension 2 linear subspace $P_2$ of $P^m$.

\medskip

This intersection is an analog of the conic line $C$ of the case $g=2$ (see [L01],
Theorem 0.7; Section 3, Figure 1).

\medskip

We consider the set of hyperplanes in $P^m$ containing $P_2$; this set is isomorphic to
$P^1(\n F_p)$. Intersections of $\goth D_{2g}$ with these hyperplanes form a partition of
$\goth D_{2g}$ which we denote by $\goth H$ (partition of hyperplane sections).

\medskip

{\bf Conjecture 1.5.6.} The restriction of I-partition on $\goth D_{2g}$ is the
intersection of $\goth H$-partition with the restriction of D-partition on $\goth D_{2g}$.

\medskip

Particularly, for $g \ne 2$ any part of D-partition in $\goth D_{2g}$ contains $p+1$
parts of I-partition which are indexed by elements of $P^1(\n F_p)$.

\medskip

{\bf Conjecture 1.5.7.} The action of Galois group $\Gal(K^p/K^1)=\n F_{p^2}^*/\n F_p^*$
on the set of parts of I-partition that are contained in any fixed part of D-partition on
$\goth D_{2g}$ coincides with the natural action of $\Gal(K^p/K^1)$ on $P^1(\n F_p)$.

\medskip

{\bf Remark 1.5.8.} I think that the proof of Conjectures 1.5.6, 1.5.7 can be got easily
using the same ideas as the ones used in proofs of other theorems of the present paper.
For the conjecture 1.5.6 this opinion is based on the analogy with the case $g=2$ (see
[L01]). Lemma 3.2.8 and Remark 3.2.9 give an argument in favor of Conjecture 1.5.7 for
odd $g$. I think that the conjecture 1.5.7 is true also for even $g$ because we cannot
imagine any other Galois action.

\medskip

{\bf (1.5.9)} To formulate the remaining theorems of Section 3, we need to quote some
definitions of [L01], Section 3 (slightly changing notations):

\medskip

A restriction of a Hecke correspondence $\goth T_{p}$ to $V$ can be considered as a
linear combination of correspondences $\goth T_{V,i}=\goth T_{p,V,i}$ between $V$ and
subvarieties $V_i \subset X$:

$$\goth T_{p}|_V=\sum_{i} \rho_i\goth T_{V,i} \eqno{(\hbox{[L01]}, 3.1)}$$

where $i$ runs over the set of classes of I-equivalence, $\rho_i$ are multiplicities.

Some subvarieties $V_i$ can coincide but all $\goth T_{V,i}$ are, by definition,

different. We denote the graph of $\goth T_{V,i}$ by $\Gamma_{V,i} \subset
V \times V_i$ and projections $ \Gamma_{V,i}\to V$, $ \Gamma_{V,i}\to V_i$

by $\pi_{1}=\pi_{V,i,1}$, $\pi_{2}=\pi_{V,i,2}$.

\medskip

The degree of $\pi_{2}$ is an important invariant of the corresponding component of
$\goth T_p(V)$, because it enters as the multiplicity in the Abel-Jacobi map. For $\goth
T_p=T_p$ we prove:

\medskip

{\bf (*) Proposition 3.4.1.} For any $g$ the degree of $\pi_{2}$ for any good component
is 1.

\medskip

{\bf Proposition 3.4.3.} For $g=3$ the degree of $\pi_{2}$ for the bad component is $p+1$.

\medskip

{\bf 1.6. Contents of Section 4: Case $\goth T_p=T_{p,\goth i}$.}  \nopagebreak

\medskip

Let $W$ be a module over $\n Z/p^2$, $W_p$ its submodule of $p$-torsion. $W$ has
invariants $d_1=\dim W/W_p$, $d_2=\dim W_p$.

{\bf (1.6.1)} The analog of (1.3.3) for the present case is the following. The set
$T_{p,\goth i}(t)$ is isomorphic to the generalized Grassmannian $G_I(\goth i,g, 2g)$,
i.e. the set of isotropic $W \subset (A_t)_{p^2}$ such that $d_1(W)=g-\goth i$,
$d_2(W)=g+\goth i$. Like in Subsection 1.5, we consider only the open part $S(g,\goth
i)^{open}$ of $S(g,\goth i)$.

Let $\s \in S(g,\goth i)^{open}$. The analog of $T(\s)$ of 3.2.2 is a pair of matrices
$\mu_1(\s),\mu_2(\s)\in M_g(K)$ (see 4.1.2.5). The analog of $\goth D$-partition for this
case is defined in 4.2.9, 4.2.10 (see Remark 4.2.11 explaining why this definition is
distinct from its analog for $\goth T_p=T_p$). Namely, we have a disjoint
union$$S(g,\goth i)=\goth D_*\cup(\bigcup_{j=0}^{\goth i}\goth D_j)\eqno{(1.6.2)}$$
(elements of $\goth D_*$ are the most special elements, elements of $\goth D_{\goth i}$
are the most general elements).

\medskip

{\bf Conjecture 1.6.3.} D-partition of $S(g,\goth i)$ is a subpartition of $\goth
D$-partition.

\medskip

We know the answer for the general part $\goth D_{\goth i}$.

\medskip

{\bf (*) Theorem 4.2.15.} Let $g$, $\goth i$ be arbitrary. If $g-\goth i$ is odd then $\goth D_{\goth i}$ is a D-partition set.
If $g-\goth i$ is even then $\goth D_{\goth i}$ is a union of two D-partition sets. 
\medskip

These sets are described as follows. Let $\sigma_i$, $\sigma_j\in S(g,\goth i)^{open}$ be 
2 elements. We associate to them $(g-\goth i)\times (g-\goth i)$-matrices $G_{22i}, G_{22j}$ with entries in $\n Z$ (see 4.1.5 and the
above formulas). 
$\sigma_i$, $\sigma_j$ are D-equivalent iff the ratio $\det \tilde G_{22i}
/ \det \tilde G_{22j}$ is a square in $\n F_p^*$.
\medskip
From now we consider the case $g=3$, $\goth i=1$.

\medskip

{\bf (*) Theorems 4.2.16, 4.2.17.} Conjecture 1.6.3 is true for $g=3$, $\goth i=1$.

\medskip

We refer to $\goth D_*$, $\goth D_0$, $\goth D_1$ as of special, intermediate, general
type respectively. Further, there exists a disjoint union $\goth D_*=\goth D_{0,*}\cup
\goth D_{1,*}$ defined in Theorem 1.6.5, b, c.

\medskip

{\bf (*) Theorem 4.2.18.} $\goth D_{0,*}$ is a part of D-partition, and moreover it is
the only part such that the corresponding component of $T_p(V)$ is $V$ itself.

\medskip

Theorems 4.2.15 --- 4.2.18 show that the D-partition of $S(3,1)$ is a subpartition of the
partition $S(3,1)=\goth D_{0,*}\cup \goth D_{1,*}\cup \goth D_0\cup \goth D_1$. I do not
know what is the restriction of D-partition on $\goth D_0$ and $\goth D_{1,*}$. See
Section 5.3 for propositions which are the first step to the solution of problem of
finding of this partition. Computer calculations for the case $p=3$ show that there are
more than one part of D-partition sets in $\goth D_{1,*}$. Apparently, if $p$ tends to
infinity, the quantity of these parts should be one.

\medskip

{\bf Theorem 4.2.19.}  Irreducible components of Shimura subvarieties of $T_{p,1}(V)$
that correspond to the parts of D-partition that are contained in $\goth D_0$ are defined
over $K^p$.

\medskip

{\bf Remark.} All other irreducible components of $T_{p,1}(V)$ are defined over ``small''
extensions of $K^1$, see 4.2.21. Unlike the case of $\goth T_p=T_p$, for $g=3$, $\goth
T_p=T_{p,1}$ the Galois action is non-trivial for the intermediate case.

\medskip

{\bf (1.6.4)} Analogs of Conjectures 1.5.6, 1.5.7 (i.e. the description of I-partition on
$\goth D_0$ and the Galois action on the set of these parts) is given in 4.2.20. It is
necessary to emphasize that we have 2 natural options for the partition; I do not know
which of these 2 options holds.

\medskip

Let us give now an analog of the Lemma 3.2.3 --- the interpretation of the sets $\goth
D_*$ in terms of $d_1$, $d_2$ of $\n F_{p^2}W_{\s}$. Their possible values are the
following (dimensions over $\n F_{p}$): (4,4), (4,6), (2,4), (2,6).

\medskip

{\bf (**) Theorem 1.6.5.} a) $\s\in \goth D_0 \cup \goth D_1 \iff d_1, d_2(\n
F_{p^2}W_{\s})= (4,4)$ or $(4,6)$;

b) $\s\in \goth D_{1,*} \iff d_1, d_2(\n F_{p^2}W_{\s})= (2,6)$;

c) $\s\in \goth D_{0,*} \iff d_1, d_2(\n F_{p^2}W_{\s})= (2,4)$. $\square$

\medskip

Finally, we can consider the following geometric construction. Factorization of $W_{\s}$
by $pW_{\s}$ is a projection $pr$ from $G_I(\goth i,g, 2g)$ to $G_I(g-\goth i,2g)$ (with
a fibre $\n F_p^{(g-\goth i)(g-\goth i+1)/2}$). The projection of $\goth D$-partition on
$G_I(2,6)$ is described as follows:

\medskip

{\bf (**) Theorem 1.6.6.} a) $pr(\goth D_*)$ is an intersection of $G_I(2,6)$ (in its
Pl\"ucker embedding) and a codimension 4 linear subspace.

b) $pr(\goth D_* \cup \goth D_0)$ is a section of $G_I(2,6)$ by a hypersurface.

\medskip

{\bf 1.7. Contents of Section 5: Miscellaneous.} \nopagebreak

\medskip

Section 5 contains some separate results that form neither a complete theory nor a
logical chain, and at the moment they have no direct application to the program of
generalization of Kolyvagin's theorem, although maybe will acquire this application in
future. These results can be treated as an introduction to the future investigations.
Some of them are not proved completely and are formulated as conjectures.

\medskip

{\bf Section 5.1. Case of the Hecke correspondence $T_{p,2}$.} This Hecke correspondence
is not as interesting for us as $T_{p,1}$, because no one irreducible component of
$T_{p,2}(V)$ is defined over $K^p$ (see (5.1.3), (5.1.4)), and hence we cannot apply this
correspondence in order to generalize Kolyvagin's method. We find an explicit description
of $\goth D_*$ (Theorem 5.1.1; it turns out that it is non-empty only if $p \equiv 3 \mod
4$) and show that $\goth D_0=\emptyset$ (Lemma 5.3.5.4). Further, we get evidence that
Conjecture 1.6.3 is true for this case (Theorems 4.2.16, 5.1.2).

\medskip

{\bf Section 5.2.} Here we investigate the action of $T_p$ on Shimura subvarieties of
$T_p(V)$. This could be useful for the program of generalization of Kolyvagin's theorem,
because maybe consideration of the action of $T_p$ on various Shimura subvarieties of
$T_p(V)$ will permit us to eliminate by elementary way the unpleasant coefficient $\goth
x_p$ of [L04.2], (1.8), (1.9).

\medskip

{\bf Section 5.3.} The restriction of D-partition to some sets $\goth D_i$ is unknown.
Nevertheless, for $g=3$ it is possible to prove some weakened conditions of
``equivalence'' of $\s_1$, $\s_2\in \goth D_j$ for the cases $\goth i=1$, $j=*$ (Theorem
5.3.2) and $j=0$ (Theorem 5.3.1), and for the cases $\goth i=2$, $j=*$ (Theorem 5.3.3)
and $j=1$ (Conjecture 5.3.5).

\medskip

{\bf Section 5.4.} Here it is proved that there is no coincidences between irreducible
components of $T_{p,\goth i}(V)$ for different $\goth i$.

\medskip

\medskip

{\bf 1.8. Conjectures and possibilities of further investigation.} \nopagebreak

\medskip

Let us give some possibilities of further investigation related to the subject of the
present paper in order of increasing complexity. See also possibilities of further
investigation in [L04.2].

\medskip

1. The reader sees that most theorems of the present paper (they are marked by (*)) are
not proved in all generality. So, the first problem is to prove them completely.

The main reason why most theorems are not proved in all generality is the necessity to
treate all open charts of the Grassmannian corresponding to $\goth T_p(t)$, and not only
the one such chart as it is done in the present paper. The reader can compare the length
of the proof of the Proposition 3.4.3 where the consideration of all open charts is
inevitable, with the length of proofs of other theorems, and evaluate the length of the
complete proofs of all theorems for the case of arbitrary $g$. Apparently, it is
necessary to invent a new method of proof if we do not want to have a several
hundred-page paper. See Remark 3.3.2.

Some theorems marked by (**) are not proved at all. These theorems, as well as
conjectures 1.5.6, 1.5.7, 4.2.20, can be proved easily using the same methods.

\medskip

2. What is the degree of the hypersurface of Theorem 1.6.6, b? What is the interpretation
of parts of I-partition of Conjecture 4.2.20 in terms of the geometry of $G_I(2,6)$ (see
Theorems 1.6.5, 1.6.6)?

It is interesting to find the degree of $\pi_2$ (see 1.5.9) for the bad components for
$g=4$, $j=2$. Is it a multiple of $p+1$?

\medskip

3. Prove Conjecture 1.5.4. We can expect that for any fixed pair $(g,j)$ (notations of
(1.5)) there exists a proof of 1.5.4 for these $(g,j)$ similar to the proof of 3.3.4
(case $g=3, j=2$). Does exist a uniform proof for all $(g,j)$?

\medskip

4. Investigate the case $p$ splits in $K$. According [L01], this case is much more
complicated. Maybe it will be possible to find a good description of this situation.

\medskip

5. To study more complicated types of pairs (and more generally, $n$-tuples) of Shimura
varieties $X$, $V$. For example, if the reflex field of $V$ is any $CM$-field, then we
get more interesting examples of Galois action.

\medskip

6. Investigate in more details the case of $T_{p, \goth i}$. I think that the formula
(4.2.9) for $\goth D_*$ is too rough. In order to get a correct definition of $\goth
D$-partition it is necessary to take into consideration the rank of $\mu_i(\s)$.

\medskip

7. To find D-partition of the intermediate and special cases for $g=3$, $\goth i=1$. If
it is impossible to prove a general theorem, then it would be interesting to find
partitions for some small $p$ by means of computer calculations.

\medskip

8. To prove analogous results for the functional case (Drinfeld modules and more generally abelian T-motives of Anderson). We know that the functional case is usually easier. Maybe in some cases the phenomenon of Section 4.4 of [L04.2] (the pseudo-Euler systems for the case $g=3$, $\goth T_p=T_p$ are identically 0) does not occur, or it will be possible to use the methods of [Z85] (see [L04.2], appendix 3) in order to find the Abel-Jacobi image of bad components? 

\medskip

{\bf Section 2.} D-equivalence.  \nopagebreak

\medskip

Let $X_c$ be a connected component of a Shimura variety of level $N$. Recall that a
$p$-Hecke correspondence $\goth T_p$ on $X_c$ is a diagram

$$\matrix & & X_{p,c} & & \\ & \overset{\pi_1}\to{\swarrow}&&
\overset{\pi_2}\to{\searrow}& \\ X_c &&&& X_c \endmatrix$$

where $X_{p,c}$ is a Shimura variety of (roughly speaking) level $pN$. In this section
$X_{p,c}$ and all subsequent objects with subscript $p$ depend on $\goth T_p$ --- the
type of Hecke correspondence. For $t\in X_c$ we have by definition: $\goth
T_p(t)=\pi_2(\pi_1^{-1}(t))$. It turns out that instead of considering equivalences on
the set $\goth T_p(t)$, it is easier to consider equivalences on the set $\pi_1^{-1}(t)$
and then to apply $\pi_2$ to this equivalence relation.

A non-formal definition of D-equivalence on the set $\goth T_p(t)$ is given in the
introduction; its definition on the set $\pi_1^{-1}(t)$ is analogous. For the exact
definition see 2.8.

We use notations of [D71]. Let us consider Deligne data that define a Shimura variety $X$
of a fixed level, namely, a reductive group $G$ over $\n Q$, a map $h: \n S \to G$ over
$\n R$ (where $\n S=\Res_{\n C/\n R}(G_m)$), and a level $N$ subgroup $\Cal K \subset
G(\n A_f)$. We denote by $\Cal D$ the Hermitean symmetric domain that corresponds to $G$,
$h$. So,

$$X=[\Cal D \times \Cal K\backslash G(\n A_f)]/G(\n Q)$$

Let $C(X)$ be the set of irreducible components of $X$; it is isomorphic to the set of
double cosets $\Cal K\backslash G(\n A_f)/G(\n Q)$. The irreducible component that
corresponds to the trivial double coset $\Cal K \cdot 1 \cdot G(\n Q)$ will be denoted by
$X_c$. Further, we denote $\Gamma=\Cal K \cap G(\n Q)$, so $X_c=\Cal D/\Gamma$ (the
action of some elements of $G(\n Q)$, and hence of $\Gamma$, on $\Cal D$ can be trivial).

We use the same notations for $X_p\twoheadrightarrow X$ and $V\hookrightarrow X$ using
subscripts $p$, $V$ respectively. Let $\Cal K_p \subset \Cal K$ be a level $p$ subgroup,
$$X_p=[\Cal D \times \Cal K_p\backslash G(\n A_f)]/G(\n Q)$$ the corresponding Shimura
variety, $\pi=\pi_1: X_p \to X$ and $C(\pi): C(X_p) \to C(X)$ the natural projection
maps, and $\Gamma_p=\Cal K_p \cap G(\n Q)$. Let

$$\Cal K =\bigcup_{i\in R}\Cal K_p r_i\eqno{(2.1)}$$

be a coset decomposition, here $R$ is the set of these cosets. There is an equivalence
relation on $R$ (we call it 1-equivalence and denote by $\sim$):

$$r_i\sim r_j \iff \Cal K_p r_i G(\n Q) = \Cal K_p r_j G(\n Q) \eqno{(2.2)}$$

We denote the quotient set of $R$ by 1-equivalence by $\goth Q$.

\medskip

{\bf Remark.} (2.2) can be written in the form $\Cal K_p r_i \Gamma = \Cal K_p r_j
\Gamma$, because $r_i\sim r_j \iff r_j \in \Cal K_p r_i G(\n Q)$, i.e. $r_j=kr_ig$, $k
\in \Cal K_p$, $g \in G(\n Q)$, $\lw g=(r_i)^{-1}k^{-1}r_j \lw g \in \Gamma$.

\medskip

{\bf Proposition 2.3.} $C(\pi)^{-1}(X_c)$ [i.e. the set of irreducible components of
$\pi^{-1}(X_c)$] is $\goth Q$.

\medskip

{\bf Proof.} Since $C(X_p)=\Cal K_p\backslash G(\n A_f)/G(\n Q)$, we have
$C(\pi)^{-1}(X_c)=\Cal K_p\backslash \Cal K G(\n Q)/G(\n Q)$. When we choose a set of
representatives of these double cosets, we can choose them in $\Cal K$; hence we can
choose a subset of $\{r_i\}_{i\in R}$. It is clear that $r_i\sim r_j$ is exactly
1-equivalence. $\square$

\medskip

We denote by $X_{p,c}$ the irreducible component of $X_p$ that corresponds to the trivial
double coset $\Cal K_p \cdot 1 \cdot G(\n Q)$. Obviously $X_{p,c}=\Cal D/\Gamma_p$,
$X_{p,c} \subset \pi^{-1}(X_c)$, and we denote by $\pi_c: X_{p,c} \to X_c$ the
restriction of $\pi$ to $X_{p,c}$. Further we consider for simplicity only the following
case:

\medskip

If $\gamma\in\Gamma$ acts trivially on $\Cal D$ then $\gamma\in\Gamma_p$.

\medskip

In this case $\Gamma/\Gamma_p \hookrightarrow \Cal K/\Cal K_p=R$. We denote by $R_c$ the
image of $\Gamma/\Gamma_p$ in $R$. Obviously, $\deg \pi_c: X_{p,c} \to X_c$ is
$\#(\Gamma/\Gamma_p)$. Further, for $t \in \Cal D$ we have a formula

$$\pi_c^{-1}(\bar t)=\{\overline{(t^{r_i}\times r_i)}\}_{i\in R_c}\eqno{(2.4)}$$

where $\bar t\in X_c$ is the projection of $t \times 1$, and $(t^{r_i}\times r_i) \in
\Cal D \times G(\n A_f)$.

\medskip

Now let us consider $V\hookrightarrow X$. Let $G_V \subset G$ be an inclusion of groups
commuting with $h_V: \n S \to G_V$;\ \  $\Cal K_V=G_V(\n A_f)\cap \Cal K$ a level
subgroup; $V$ the Shimura subvariety of $X$ corresponding to $(h_V, G_V, \Cal K_V)$ and
$V_{c}$ the irreducible component of $V$ that corresponds to the trivial double coset
$\Cal K_V \cdot 1 \cdot G_V(\n Q)$.

\medskip

{\bf Proposition 2.5.} The set of irreducible components of $\pi_c^{-1}(V_c)$ is the
quotient set of $\{r_i\}_{i \in R_c}$ by the equivalence relation:

$$r_i\sim r_j \iff \Gamma_p r_i G_V = \Gamma_p r_j G_V \eqno{(2.6)}$$

\medskip

{\bf Proof.} A version of this proposition for Hecke correspondences is [L01],
Proposition 4.4, (1). The proof of the present version of this proposition is completely
analogous. $\square$

\medskip

This equivalence relation is the I-equivalence. Slightly changing notations of (1.4.1),
we denote the quotient set of $\Gamma/\Gamma_p$ by this equivalence relation by
$(\Gamma/\Gamma_p)_I$.

\medskip

The above varieties form the diagram:

$$\matrix \pi_c^{-1}(V_c)&&\to&&X_{p,c}
\\ &&&&& \searrow
\\ \downarrow&&&&\downarrow&&X_p
\\
\\ V_c&&\to&&X_c&&\downarrow
\\&\searrow&&&&\searrow
\\&&V&&\to&&X\endmatrix\eqno{(2.7)}$$

Now we must ``unify'' (in partition sets) those irreducible components of
$\pi_c^{-1}(V_c)$, which belong to a Shimura subvariety. We choose and fix $i\in R_c$, we
choose the corresponding representative $r_i$ in $\Gamma$ and we denote it by $s$. For
any group $\goth G$ we denote the conjugate group $s\goth Gs^{-1}$ by $\goth G^s$. We
define a Shimura variety $V_{p,s}$ that corresponds to the subgroup $\Cal K_{V,p,s}=\Cal
K_V\cap \Cal K_p^{s^{-1}}=G_V(\n A_f)\cap \Cal K_p^{s^{-1}}$. The projection $V_{p,s}\to
V$ is denoted by $\pi_{V,s}$. There exists an $s$-conjugate variety $V_{p,s}^s=s(\Cal
D_V) \times \Cal K_V^s\backslash G_V(\n A_f)^s]/G_V(\n Q)^s$ which enters in the diagram
(2.7) as follows:

$$\matrix V_{p,s}^s&\overset{i_{p,s}}\to{\hookrightarrow}&X_p
\\ \downarrow&&\downarrow
\\V&\to&X\endmatrix$$

where $i_{p,s}$ is induced by identical inclusions of corresponding subgroups of $G(\n
A_f)$, and the left vertical arrow (denoted by $\pi_{V,s}^s$) is induced by the map
$t\mapsto s^{-1}(t)$ for $t\in s(\Cal D_V)$, $s^{-1}(t)\in \Cal D_V$.

\medskip

{\bf Definition 2.8.} 2 elements $\bar s_1$, $\bar s_2$ of $(\Gamma/\Gamma_p)_I$ are
called D-equivalent if there exists $s\in \Gamma$ such that both irreducible components
of $\pi_c^{-1}(V_c)$ that correspond to $\bar s_1$, $\bar s_2$ are contained in
$V_{p,s}^s$. 2 elements $s_1$, $s_2$ of $\Gamma/\Gamma_p$ are called D-equivalent if
their images in $(\Gamma/\Gamma_p)_I$ are D-equivalent.

\medskip

Now let us find an algebraic criterion of D-equivalence. Let

$$\Cal K_V =\bigcup_{i\in R_{V,s}}\Cal K_{V,p,s} r_{V,i,s}$$

be a coset decomposition. There is the relation of 1-equivalence on $R_{V,s}$ defined
like in (2.2), and the corresponding quotient set $\goth Q_{V,s}$.  Let $J_{V,s}$ be a
subset of $R_{V,s}$ defined as follows:

$$j\in J_{V,s} \iff sr_{V,j,s}s^{-1} \in \Cal K_pG(\n Q)$$

It is obvious that $J_{V,s}$ is a union of partition sets of 1-equivalence on $R_{V,s}$;
we denote the corresponding quotient set by $\goth Q_{V,0,s}\subset \goth Q_{V,s}$. It is
easy to see that there exists an inclusion $\alpha_s: \goth Q_{V,0,s} \to
(\Gamma/\Gamma_p)_I$ is defined at the level of representatives as follows. Let
$r_{V,i,s}\in \goth Q_{V,0,s}$, $sr_{V,i,s}s^{-1}=kg$, where $k\in \Cal K_p$, $g\in G(\n
Q)$. Then $g\in \Gamma$. By definition, $gs$ is a representative of
$\alpha_s(r_{V,i,s})$. Particularly, $\alpha_s(1)=s$.

\medskip

{\bf Proposition 2.9.} $\goth Q_{V,0,s}$ is isomorphic to the set of irreducible
components of $(\pi_{V,s}^s)^{-1}(V_c)\cap X_{p,c}=\pi_c^{-1}(V_c)\cap V_{p,s}^s$ (both
intersections are in $X_p$). Moreover, the corresponding inclusion is exactly $\alpha_s$.

\medskip

{\bf Proof.} An easy calculation in groups and cosets, with use of (2.4). $\square$

\medskip

{\bf Proposition 2.10.} 2 elements $r_1$, $r_2\in \Gamma/\Gamma_p$ are D-equivalent iff
the double cosets $\Cal K_pr_i\Cal K_V$ coincide ($i=1,2$), i.e. iff

$$r_2\in\Cal K_pr_1\Cal K_V\eqno{(2.11)}$$

{\bf Proof.} We see that 2 elements $r_1$, $r_2\in \Gamma/\Gamma_p$ are D-equivalent iff
$r_2\in \im(\alpha_{r_1})$. This means that $\exists i\in R_{V,r_1}$ such that $r_1
r_{V,i,r_1}r_1^{-1}=kg$, where $k\in \Cal K_p$, $g\in G(\n Q)$, and such that $r_2=gr_1$.
This is equivalent to (2.11). $\square$

\medskip

Now we must replace the coset decomposition (2.1) by the coset decomposition (1.3.2). In
order to simplify notations, until the end of this section $\tau_p$ (resp. $S(g)$) will
mean either $\tau_p$ (resp. $S(g)$ ) itself (for $\goth T_p=T_p$) or $\tau_{p,\goth i}$
(resp. $S(g,\goth i)$ ) of (1.3.1) (for $\goth T_p=T_{p,\goth i}$). In both cases we have
$\Cal K_p=\tau_p \Cal K \tau_p^{-1}\cap\Cal K $, $\Gamma _p=\tau_p \Gamma \tau_p^{-1}\cap
\Gamma$. We need a lemma on relation between 2 sets of representatives of cosets.

\medskip

{\bf Lemma 2.12.} Let $\Gamma \tau_{p}\Gamma =\bigcup_{j\in \goth S(g)}\Gamma \sigma_j'$
be a coset decomposition such that all representatives $\sigma_j'$ satisfy
$$\s_j'\in\tau_{p}\Gamma\eqno{(2.13)}$$ We denote $$r_j=\tau_{p} \s_j^{\prime
t-1}\eqno{(2.14)}$$ Then $\Gamma = \cup_j \Gamma _p r_j$, and vice versa. $\square$

\medskip

{\bf Remark.} We use notation $\sigma_j'$ instead of $\sigma_j$, because not all
$\sigma\in S(g)$ satisfy (2.13). Namely,

$$\hbox{$\sigma\in S(g)$ satisfy (2.13) iff $\sigma\in S(g)^{open}$.} \eqno{(2.15)}$$

Substituting (2.14) in (2.11) we get immediately the following

\medskip

{\bf Theorem 2.16.} 2 elements $\s_1$, $\s_2\in S(g)$ are D-equivalent iff the double
cosets $\Cal K\s_i\Cal K_V$ coincide ($i=1,2$). This condition can be rewritten as
follows: $$\exists \goth g\in \Cal K_V\hbox{ such that }\sigma_1\alpha(\goth g)
\sigma_2^{-1}\in\Cal K \eqno{(2.17)}$$

{\bf Remark 2.18.} We need to consider $\pi_2$-images of irreducible components of the
above subvarieties of $V_{p,c}$. In principle, $\pi_2$ can glue some of them, see [L01]
for the detailed consideration of this phenomenon. But if $p$ is split in $K$ then it is
easy to check that (for the case of $V$, $X$ defined in 3.1) there is no such glueing.
So, Theorem 2.16 describes us also D-equivalence on the set $\goth T_p(t)$, and not only
on the set $\pi_1^{-1}(t)$.

\medskip

Now we need a lemma that describes $L_{k_1}$ of (1.4.2) --- the field of definition of
irreducible components of $T_p(V)_k$. Let $s\in\Gamma$ be an element such that $s\in
S(g)_{k_1}$, i.e. $T_p(V)_k\subset V_{p,s}^s$. For simplicity, we formulate it only for
the case of $V$, $X$ defined in 3.1, and moreover for the case when $\goth r\ne \goth s$,
i.e. when the reflex field of $V$ is $K$:

\medskip

{\bf Lemma 2.19.} $L_{k_1}$ is the abelian extension of $K$ that corresponds (by the
global class field theory) to the group $\det \Cal K_{V,p,s}$, where $\det: G_V(\n A_{\n
Q})\to\Res_{K/\n Q}(G_m)(\n A_{\n Q})=I_K$ is the determinant map ($G_m$ is the
multiplicative group, $I_K$ is the idele group). $\square$

\medskip

{\bf Section 3.} Case $\goth T_p=T_p$.  \nopagebreak

\medskip

{\bf Subsection 3.1.} Definition of $V$.  \nopagebreak

\medskip

Here we define the inclusion $V \subset X$. $X$ is the Siegel variety of genus $g$ and
level $N$, i.e. it corresponds to the following Deligne data:

\medskip

1) $G=GSp_{2g}(\n Q)$;

\medskip

2) The map $h_X:\Res_{\n C/\n R}G_m\to G$ and the subgroup $\Cal K=\Cal K_N \subset G(\n
A_f)$ of level $N$ are the standard ones: they are defined like in one-dimensional case,
but instead numbers we have in mind $g\times g$-blocks of matrices.

\medskip

Recall the explicit formulas for the unitary group $G_V=GU(\goth r, \goth s)$ where
$\goth r+\goth s=g$, $K=\n Q(\rd)$. If $L$ is a field such that $\rd \not\in L$ (we shall
consider only cases $L=\n F_p, \ \n Q, \ \n Q_p$) and $\goth g\in G_V(L)$ then $\goth
g=h+k\sqrt{-\Delta}$, where $\goth g$ satisfy

$$\goth gE_{\goth r \goth s}\bar \goth g^t=\lambda E_{\goth r \goth s}\eqno{(3.1.1)}$$

i.e. $h,k\in M_g(L)$ satisfy

$$hE_{\goth r \goth s}h^t+\Delta kE_{\goth r \goth s}k^t=\lambda E_{\goth r \goth
s}\eqno{(3.1.2)}$$

$$kE_{\goth r \goth s}h^t=hE_{\goth r \goth s}k^t\eqno{(3.1.3)}$$

(here $E_{\goth r \goth s}=\left(\matrix 1 & 0 \\0&-1\endmatrix \right)$ where sizes of
diagonal blocks are $r$, $s$).

There exists an inclusion $\alpha:G_V \to G$ defined as follows:

$$\alpha(\goth g)=\left(\matrix h & kE_{\goth r \goth s}\\-\Delta E_{\goth r \goth s}k &
E_{\goth r \goth s}hE_{\goth r \goth s}\endmatrix \right)\eqno{(3.1.4)}$$ Let all prime
divisors of $N$ split in $K$. The trivial $\n Q$-structure on $G_V$, the standard map
$h_V:\Res_{\n C/\n R}G_m\to GU(\goth r, \goth s)$ and the level subgroup $\Cal K_V=\Cal
K\cap \alpha(G_V(\n A_f))$ are Deligne data for $G_V$. We denote the corresponding
Shimura variety by $V$. The inclusion $\alpha$ defines an inclusion $V \subset X$.

It is well-known that $\dim V=\goth r \goth s$, the reflex field of $V$ is $K$ if $\goth
r\ne \goth s$ and $\n Q$ otherwise, and if $\goth r\ne \goth s$ then the field of
definition of irreducible components of $V$ is $K^1$ --- the Hilbert class field of $K$.
We shall be interested mainly in the first non-trivial case $g=3$, $\goth r=2$, $\goth
s=1$. In this case $V$ is a Picard modular surface.

\medskip

{\bf Subsection 3.2.} Preliminary notations and lemmas.  \nopagebreak

\medskip

Recall that a description of $S(g)$ can be found in [L04.1], Section 2.3. We find
D-equivalence only on the open part of $S(g)$. This open part $S(g)^{open}$ corresponds
to the case $I=\{1,\dots ,g\}$ in notations of [L04.1]. Elements of this open part have
the form $$\s =\s(s)=\left( \matrix E_g & s\\0&pE_g \endmatrix \right)\eqno{(3.2.1)}$$
where $s\in M_g(\n Z)^{symm}$ and its entries belong to a fixed system of residues modulo
$p$ (see Remark 3.3.2 for the case when $\s\not\in S(g)^{open}$).

Since $p$ is inert in $K$, we have: $K\otimes\n Q_p =\Bbb Q_{p^2}$. We denote the residue
map by res: $\Bbb Z_{p^2}\to\Bbb F_{p^2}$ or simply by tilde, and the non-trivial
automorphism of $\Bbb Z_{p^2},\Bbb F_{p^2}$ by bar.

We need the following objects related to $\s=\s(s)$ ($E_g=1_g=1$ is the unit matrix):

$$t=t(\s)=E_g-\sqrt{-\Delta}sE_{\goth r \goth s}, \ \ T=T(\s)=\bar t(\s)E_{\goth r \goth
s}t(\s)^t=E_{\goth r \goth s}+\Delta sE_{\goth r \goth s}s\eqno{(3.2.2)}$$

Recall that $W_{\s }\subset (A_t)_p$ is defined in (1.3.5).

\medskip

{\bf Lemma 3.2.3.} $\rank T(\tilde \s ) = \dim_{\n F_{p}}(\n F_{p^2} W_{\s})-g$.

\medskip

{\bf Proof. }Straightforward. Let $\{e_i\}$ ($i=1, \dots, g)$ be an $O_K/p$-basis of
$(A_t)_p$, so $\{e_i, \rd e_i\}$ form an $\n F_p$-basis of $(A_t)_p$. $W_{\s }$ has an
$\n F_p$-basis whose matrix in $\{e_i, \rd e_i\}$ is $(E_g,\tilde s)$. The matrix of an
$\n F_p$-basis of $\rd W_{\s}$ is $(-\Delta \tilde sE_{\goth r \goth s}, E_{\goth r \goth
s})$. So, $\dim_{\n F_{p}}(\n F_{p^2} W_{\s })$ is the rank of $\left(\matrix E_g &
\tilde s\\ -\Delta \tilde sE_{\goth r \goth s}& E_{\goth r \goth s}\endmatrix \right)$.
Elementary transformations transform it to $\left(\matrix E_{\goth r \goth s} & \tilde
s\\ 0& E_{\goth r \goth s}+\Delta \tilde sE_{\goth r \goth s}\tilde s\endmatrix
\right)=\left(\matrix E_g & \tilde s\\ 0&T(\tilde \s )\endmatrix \right)$.  $\square$

\medskip

Recall that the $\goth D$-partition is defined by (1.5.1). We have a

\medskip

{\bf Corollary 3.2.4.} $\s_i$ is $\goth D$-equivalent to $\s_j$ $\iff$ ranks of $T(\tilde
\s_i)$, $T(\tilde \s_j)$ are equal. Particularly, $\s \in \goth D_g \iff T(\tilde \s ) =
0$ (for even $g$), and $\s \in \goth D_{2g} \iff t(\tilde \s ),T(\tilde \s )$ are
invertible. $\square$

\medskip

{\bf Remark 3.2.5.} $\goth D_{2g}$ and its complement in $S(g)$ have a geometric
interpretation in the Pl\"ucker embedding of $S(g)=G_I(g,2g)(\n F_p)$. To simplify
notations, we consider the case $g=3$.  We denote coordinates in the space $P^{m}(\n
F_p)$ ($m=\left(\matrix 2g\\g\endmatrix \right)-1$) by $v_{i_1 i_2 i_3}$ where $i_j\in
\{1,2,3,1',2',3'\}$ (for example, $v_{123}$-th coordinate of any $\s $ is 1 and
$v_{1'2'3'}$-th coordinate of any $\s $ is $\det \tilde s$).

It is clear that $$\det \tilde t=\goth f_1(\tilde s)+\rd\goth f_2(\tilde s)
\eqno{(3.2.6)}$$ where $$\matrix \goth f_1(\tilde s)=\ve_{1'23}v_{1'23}+
\ve_{12'3}v_{12'3}+ \ve_{123'}v_{123'}+ \ve_{1'2'3'}v_{1'2'3'} \\ \goth f_2(\tilde
s)=\ve_{123}v_{123}+ \ve_{1'2'3}v_{1'2'3}+ \ve_{1'23'}v_{1'23'}+
\ve_{12'3'}v_{12'3'}\endmatrix\eqno{(3.2.7)}$$ and $\ve_{***}$ are some constants ($\pm$
some power of $\d$). The same formulas hold for any $g$: indices run over all subsets of
the set of elements of the matrix $\left(\matrix 1,&\dots ,&g\\ 1',& \dots ,&g'\endmatrix
\right)$ such that each column contains 1 element of this subset; indices of the first
(resp. second) formula of (3.2.7) correspond to subsets containing even (resp. odd)
elements of the first line of this matrix.

So, $S(g)-\goth D_{2g}=\cup_{j=g}^{2g-2}\goth D_j$ is the intersection of the Pl\"ucker
embedding of $\goth P(G_I(g,2g))$ with a codimension 2 linear subspace $P_2$ of $P^m$.
$P_2$ is given by equations $\goth f_i(\tilde s)=0$, $i=1,2$. This is an analog of the
conic line $C$ of the case $g=2$ (see [L01], Theorem 0.7; Section 3, Figure 1).

\medskip

Recall that $\goth H$ is the partition of hyperplane sections of Conjecture 1.5.6.
$\s_1$, $\s_2$ belong to the same part of $\goth H$-partition iff $(\goth f_1(\tilde
s_1):\goth f_2(\tilde s_1))=(\goth f_1(\tilde s_2):\goth f_2(\tilde s_2))$. Action of
$\Bbb F_{p^2}$ on $(A_t)_p$ induces the action of $\Bbb F_{p^2}^*/\Bbb F_{p}^*$ on
$S(g)=G_I(g,2g)(\n F_p)$. We call it $\Cal F$-action. We have the following lemma which
for odd $g$ gives evidence in favor of Conjecture 1.5.7:

\medskip

{\bf Lemma 3.2.8.} Let $\s_1\in\goth D_{2g}^{open}$. If $g$ is even then all points of
the $\Cal F$-orbit of $\s_1$ belong to the same part of $\goth H$-partition. If $g$ is
odd then there is exactly 1 point of the $\Cal F$-orbit of $\s_1$ in each part of $\goth
H$-partition.

\medskip

{\bf Proof.} Let $\s_2=\rd(\s_1)$ respectively the $\Cal F$-action. We use notations of
Lemma 3.2.3 omitting for simplicity non-essential factors $E_{\goth r \goth s}$. We have:

$$\det (E_g - \rd\tilde s_1) = \goth f_1(\tilde s_1)+\rd\goth f_2(\tilde s_1)$$

$$\det (\Delta\tilde s_1 -\rd E_g) = \goth f_1(\tilde s_2)+\rd\goth f_2(\tilde s_2)$$
hence $\goth f_1(\tilde s_2)+\rd\goth f_2(\tilde s_2)=(\rd)^g(\goth f_1(\tilde
s_1)+\rd\goth f_2(\tilde s_1))$. The lemma follows immediately from this formula.
$\square$

\medskip

{\bf Remark 3.2.9.} The meaning of the Lemma 3.2.8 is the following. We can only guess
what is the action of Galois group on the set of parts of I-partition in a given good
part of D-partition. There are 2 natural actions: $\Cal F$-action and the action of $\n
Z/(p+1)$ on the set of parts of $\goth H$-partition. Lemma 3.2.8 shows that they coincide
for odd $g$, hence most likely it is the desired Galois action.

\medskip

{\bf Subsection 3.3.} Finding of D-equivalence.  \nopagebreak

\medskip

Let $i=1,2$, $s_i\in M_g(\Bbb Z)^{symm}$, $\s_i=\s(s_i)\in S(g)^{open}$ be as in (3.2.1).
We denote $t(\s_i)$, $T(\s_i)$ simply by $t_i$, $T_i$.

\medskip

{\bf Theorem 3.3.1.} $\s_1, \s_2\in\goth D_{2g}\cap S(g)^{open}$ are D-equivalent iff the
following condition (3.3.1.1) holds:

\medskip

{\bf (3.3.1.1)} Either $g$ is odd, or $g$ is even and the ratio $\det T(\tilde \sigma_i)/
\det T(\tilde \sigma_j)$ is a square in $\n F_p^*$.

\medskip

{\bf Proof.} We consider $\goth g$ of (2.17), $\goth g=h+k\sqrt{-\Delta}$, $h,k\in
M_g(\Bbb Z_p)$. The explicit calculation shows that the $g\times g$-block structure of
$\sigma_1\alpha(\goth g) \sigma_2^{-1}$ is the following:

$$\sigma_1\alpha(\goth g) \sigma_2^{-1}=\left(\matrix A_{11} & p^{-1}A_{12} \\
A_{21}&A_{22}\endmatrix\right)\eqno{(3.3.1.2)}$$

where entries of all $A_{ij}=A_{ij}(h,k,s_1,s_2)$ are some polynomials with integer
coefficients in entries of $h,k,s_1,s_2$ and hence $\in M_g(\n Z_p)$. This means that

\medskip

$\s_1$ is D-equivalent to $\s_2 \iff \exists \goth g \in G_V(\n Z_p)$ such that $\tilde
A_{12}=0$.

\medskip

The proof is based on the following crucial formula:

$$\bar t_1\goth g t_2=B+A_{12} E_{\goth r \goth s}\sqrt{-\Delta}\eqno{(3.3.1.3)}$$

where $B\in M_g(\Bbb Z_p)$ is some matrix. This means that we can rewrite (2.17) as
follows:

\medskip

$\sigma_1$ is $D$-equivalent to $\sigma_2 \iff \exists \goth g\in \Cal K_V\hbox{ such
that }$

$$\bar {\tilde t}_1\tilde \goth g\tilde t_2\in M_g(\Bbb F_p)\eqno{(3.3.1.4)}$$

We shall prove:

\medskip

{\bf A.} If $\exists \goth g \in G_V(\n Z_p)$ such that (3.3.1.4) holds (i.e. $\tilde
A_{12}=0$) then (3.3.1.1) holds;

\medskip

{\bf B.} If (3.3.1.1) holds then $\exists \goth g \in G_V(\n Z_p)$ such that $A_{12}=0$.

\medskip

To prove (A) we denote $\bar {\tilde t}_1\tilde \goth g\tilde t_2$ by $\tilde B$ (i.e. we
reduce (3.3.1.3) modulo $p$). Since both $\tilde t_i$ are invertible we get $\tilde \goth
g=\bar {\tilde t}_1^{-1}\tilde B\tilde t_2^{-1}$. Substituting this value of $\tilde
\goth g$ in (3.1.1), we get

$$\tilde B(E_{\goth r \goth s}\tilde T_2E_{\goth r \goth s})^{-1}\tilde B^t=\lambda
\tilde T_1\eqno{(3.3.1.5)}$$

The theory of quadratic forms over $\n F_p$ shows that (3.3.1.5) implies  (3.3.1.1).

Inversely, if (3.3.1.1) holds, then $\exists B \in GL_g(\n Z_p)$ such that $B(E_{\goth r
\goth s}T_2E_{\goth r \goth s})^{-1}B^t=\lambda T_1$. So, $\goth g=\bar
t_1^{-1}Bt_2^{-1}\in G_V(\n Z_p)$, and $A_{12}=0$.

Finally, we must check that all possible types of $T(\tilde \s)$ can be realized. This is
obvious (it is sufficient to consider diagonal matrices $s$). $\square$

\medskip

{\bf Remark 3.3.2.} The first step of reformulation of the Theorem 3.3.1 for the case of
any $\s\in\goth D_{2g}$ (not necessarily $\s\in S(g)^{open}$) is the following. If $
\s=\left( \matrix A & B\\0&D\endmatrix \right)$ where $A$, $B$, $D$ are $g\times
g$-matrices described in [L04.1], Section 2.3, (2), then we can set
$t(\s)=A+{\sqrt{-\Delta }}(B+D)$ up to possible factors $-1$ and/or $E_{\goth r \goth s}$.

\medskip

{\bf Theorem 3.3.3.} If $\s_1, \s_2\in S(g)^{open}$, $\s_1\in\goth D_{2g}$ and
$\s_2\not\in\goth D_{2g}$ then $\s_1$ and $\s_2$ are not D-equivalent.

\medskip

{\bf Proof.} We consider here only the case when $\s_1$ is D-equivalent to the matrix
$\left( \matrix E_g & 0\\0&pE_g\endmatrix \right)$, i.e. the signature of $\tilde T_1$ is
trivial (proofs for the cases of another signatures of $\tilde T_1$ are analogous). So,
we can take $\s_1=\left( \matrix E_g & 0\\0&pE_g\endmatrix \right)$, $s_1=0$ and hence
$t_1=E_g$. Assume that there exists $\goth g$ from (2.17). Then (3.3.1.3) implies that
$\tilde \goth g\tilde t_2\in M_g(\n F_p)$. This means that $\tilde k=\tilde h \tilde s_2 
E_{\goth r \goth s}$. Substituting this formula in (3.1.2) we get $\tilde h\tilde T_2
\tilde h^t=\tilde \lambda E_{\goth r \goth s}$. Since det $\tilde T_2=0$, this is
contradicts to the condition $\lambda \in \n Z_p^*$. $\square$

\medskip

{\bf Theorem 3.3.4.} If $g=3$, $r=2$, $s=1$ and 2 elements $\s_1$, $\s_2\in\goth
D_{4}\cap S(g)^{open}$ then $\s_1$ and $\s_2$ are D-equivalent.

\medskip

{\bf Proof.} det $ \tilde t_i=0$ implies that 2 eigenvalues of $\tilde s_iE_{21}$ are
$\pm\frac{1}{\sqrt{-\Delta }}$, and another one --- denoted by $r_i$ --- is in $\Bbb
F_p$. Let $w_i=(w_{i1},w_{i2},w_{i3})\in (\Bbb F_p)^3$, $v_i=(v_{i1},v_{i2},v_{i3})\in
(\Bbb F_{p^2})^3$ and $\bar v_i=(\bar v_{i1}, \bar v_{i2},\bar v _{i3})$ be eigenvectors
of $\tilde s_iE_{21}$ with eigenvalues $r_i,\frac{1}{\sqrt{-\Delta}}$ and
$\frac{-1}{\sqrt{-\Delta}}$ respectively. We denote $ L_i=\left(\matrix
w_{i1}&v_{i1}&\bar v_{i1}\\
w_{i2}&v_{i2}&\bar v_{i2}\\
w_{i3}&v_{i3}&\bar v_{i3}\endmatrix\right)$,

so det $L_i \ne 0$ and

$$ \tilde t_i=L_iD_iL_i ^{-1}\eqno{(3.3.4.1)}$$

where $D_i=\diag (1-r_i\sqrt{-\Delta},0,2)$. We assume existence of $\goth g$ satisfying
(2.17), and we denote as earlier $B=\bar {\tilde t}_1\tilde \goth g\tilde t_2$. (3.3.1.3)
implies $B \in M_3(\Bbb F_p)$.

Taking $\tilde t_i$ from (3.3.4.1), we get:

$$ \bar D_1 \bar L_1^{-1} \tilde \goth gL_2D_2=\bar L_1^{-1}BL_2\eqno{(3.3.4.2)}$$For any
matrix $X$ we have: the second line (resp. column) of $\bar D_1 X $ (resp. $X D_2$) is
$0$. Hence, (3.3.4.2) implies that both second line and column of $\bar L_1 ^{-1}BL_2$
are $0$. Further, we denote by $T_{132}$ the permutation matrix $\left(\matrix 1 & 0 &
0\\0 & 0 & 1\\
0 & 1 & 0\endmatrix\right)$, so we have equalities

$$\bar L_i=L_iT_{132}\eqno{(3.3.4.3)}$$

and (because $\bar B=B$)

$$\overline{\bar L_1^{-1}BL_2}=L_1^{-1}B\bar L_2=T_{132} \bar L_1^{-1}
BL_2T_{132}\eqno{(3.3.4.4)}$$

(3.3.4.4) implies that the third line and column of $\bar L_1 ^{-1}BL_2$ are also $0$.
So, $$\bar D_1 \bar L_1^{-1} \tilde \goth gL_2D_2=b\cdot \diag (1,0,0)\eqno{(3.3.4.5)}$$

where $b \in \Bbb F_p$.

Solving the equation $XD_2=b\cdot \diag(1,0,0)$ we see that $\bar D_1\bar L_1^{-1} \tilde
\goth gL_2=\left(\matrix \frac{b}{1-r_2\sqrt{-\Delta }}&*&0
\\0&*&0\\0&*&0\endmatrix\right)$. Analogously, solving the equation $\bar
D_1X=\left(\matrix\frac{b}{1-r_2\sqrt{-\Delta}}
&*&0\\0&*&0\\0&*&0\endmatrix\right)$ we see that

$$\bar L_1^{-1} \tilde \goth gL_2=\left(\matrix \frac{b}{(1-r_2\sqrt{-\Delta
})(1+r_1\sqrt{-\Delta })}&*&0 \\*&*&*\\0&*&0\endmatrix\right)\eqno{(3.3.4.6)}$$ All these
considerations are invertible, i.e. (3.3.4.6) is equivalent to (3.3.1.4). We denote $\bar
L_1^{-1} \tilde \goth gL_2 $ by $ \goth g_0$. Now we shall show that (3.3.1.4) and

$$\tilde \goth gE_{21}\bar {\tilde \goth g}^t=\lambda E_{21}\eqno{(3.3.4.7)}$$

(the $\Bbb F_{p^2}$-residue of (3.1.1)) can be satisfied simultaneously for any $s_1,
s_2$ of the statement of the theorem. Really, substituting $\goth g_0$ to (3.3.4.7) we get

$$\goth g_0(L_2^{-1}E_{21}(\bar L_2^{-1})^t)\bar \goth g_0^t=\lambda(\bar
L_1^{-1}E_{21}(L_1^{-1})^t)\eqno{(3.3.4.8)}$$

Analog for the present case of a theorem that eigenvectors of a symmetric matrix are

orthogonal is the following. We denote

$$D_{0i}\overset{\opr}\to{=}L_i^{-1}s_iE_{21}L_i=\diag
(r_i,\frac{1}{\sqrt{-\Delta}},-\frac {1}{\sqrt{-\Delta}})$$ so we have:
$$s_i=L_iD_{0i}L_i^{-1}
E_{21}=s_i^t=E_{21}(L_i^{-1})^tD_{0i}L_i^t$$ i.e. $$(L_i^{-1}E_{21}(L_i^{-1})^t)D_{0i}
(L_i^tE_{21}L_i)=D_{0i}$$ Since all diagonal entries of $D_{0i}$ are different we get

that $L_i^tE_{21}L_i$ is a diagonal matrix. It is invertible, we denote entries of the

inverse matrix as follows: $$(L_i^tE_{21}L_i)^{-1}=diag (d_{i1}, d_{i2}, \bar d_{i2})
\eqno{(3.3.4.9)}$$

where $d_{i1}\in\Bbb F_{p}$, $d_{i2}\in\Bbb F_{p^2}$. Taking into consideration
(3.3.4.3), (3.3.4.6) and (3.3.4.8) give us the following equations:

$$\left(\matrix b' & x_{12} & 0\\ x_{21} & x_{22} & x_{23}\\0 & x_{32} & 0
\endmatrix \right) \left(\matrix d_{21} & 0 & 0 \\ 0 & 0 & d_{22}\\ 0 & \bar d_{22} & 0
\endmatrix\right)\left(\matrix \bar b' & \bar x_{21} & 0\\ \bar x_{12} & \bar x_{22} &
\bar x_{32}\\0 & \bar x_{23} & 0 \endmatrix \right)= \lambda  \left(\matrix d_{11} & 0 &
0 \\ 0 & 0 & d_{12} \\ 0 &
\bar d_{12} & 0\endmatrix \right)$$

where $b'= \frac {b} {(1-r_2 \sqrt{-\Delta })(1+r_1 \sqrt{-\Delta })}$. For any $s_1$,
$s_2$ a solution can be found immediately, for example $x_{12}$= $x_{21}$= $x_{22}$= $0$,
$b\in\Bbb F_{p}^*$ arbitrary, $\lambda= \frac{b'\bar b'd_{21}}{d_{11}}$, $x_{23}$=$1$,
$x_{32}$=$\frac{\lambda \bar d_{12}}{d_{22}}$. So, we have got $\tilde \goth g$
satisfying (3.3.1.4), (3.3.4.7). Further, we have

\medskip

{\bf Lemma 3.3.4.10.} The reduction map $G_V(\n Z_p) \to G_V(\n F_p)$ is surjective.
$\square$

\medskip

Application of this lemma gives us the theorem. $\square$

\medskip

{\bf Theorem 3.3.5.} If $g$ is odd then all components of $T_p(V)$ do not coincide with
$V$ itself.

\medskip

{\bf Proof.} Obvious. Really, we have the following criterion of coincidence ([L01],
Section 4a, Proposition 4.4, (4)): $$ \Gamma_V\s^{-1} \cap \Gamma \ne
\emptyset\eqno{(3.3.5.1)}$$

(the inclusion $G_V \to G$ is $\alpha$). If $\goth g\in \Gamma_V$ satisfies $\alpha(\goth
g)\s^{-1} \in \Gamma$ then $\lambda(\goth g)=p$. This is impossible by a trivial reason
(in this case we have $N_{K/\n Q}(\det \goth g)=p^g$ which for odd $g$ contradicts to a
condition $p$ inert in $K$). $\square$

\medskip

{\bf Theorem 3.3.6.} For any $g$, any $\goth r\ne \goth s$ the field of definition of any
good component is $K^p$.

\medskip

{\bf Proof.} We fix $\s\in S(g)^{open}$. According 2.15, we can take $r=\tau_p\s^{t-1}$
(see 2.14). We apply Lemma 2.19 ($s$ of 2.19 is our $r$). It is easy to see that for
points $\goth l$ of $K$, $\goth l\ne p$, the $\goth l$-component of $\det \Cal K_{V,p,r}$
contains $O_{K_{\goth l}}^*$, so we must find only the $p$-component of $\det \Cal
K_{V,p,r}$. Let $\goth g=h+k\rd\in G_V(\n Z_p)$ be as earlier ($G_V(\n Z_p)$ is the
$p$-component of $\Cal K_V$). The $p$-component of $\Cal K_p$ is the set of matrices
whose (2,1)-block is $\equiv 0 \mod p$.

\medskip

{\bf 3.3.6.1.} The following 3 conditions are equivalent:

$$\hbox{$\goth g$ belongs to the $p$-component of $\Cal K_{V,p,r}$}\eqno{(3.3.6.2)}$$

$$\hbox{$\iff r \alpha(\goth g) r^{-1}$ belongs to the $p$-component of $\Cal
K_p$}\eqno{(3.3.6.3)}$$

$$\hbox{$\iff $ the (2,1)-block of $r \alpha(\goth g) r^{-1}$ is $\equiv 0 \mod
p$}\eqno{(3.3.6.4)}$$

We denote the (2,1)-block of $r \alpha(\goth g) r^{-1}$ by $\goth B_{21}$. There is a
formula for it (practically, this is 3.3.1.3): $$\bar t E_{\goth r \goth s}\goth g t^t =
\goth A-\goth B_{21}\rd\eqno{(3.3.6.5)}$$ where $t=t(\s)$ and $\goth A\in M_g(\n Z_p)$.
If $\det \tilde t\ne 0$ (i.e. $\s\in \goth D_{2g}$) then $\tilde \goth B_{21} = 0
\implies \det \goth g \in \n Z_p+p\rd\n Z_p$. It is easy to see that $\det \goth g$ can
take any values in $\n Z_p^*+p\rd\n Z_p$, i.e. the field of definition of the good
component is $K^p$. $\square$

\medskip

Analogously, we get the following

\medskip

{\bf Theorem 3.3.7.} The field of definition of all bad components is $K$. $\square$

\medskip

{\bf Subsection 3.4. Finding of the degree of $\pi_2$. } \nopagebreak

\medskip

We use notations of 1.5.9.

\medskip

{\bf Proposition 3.4.1.} For any $g$ the degree of $\pi_{2}$ for any good component is 1.

\medskip

{\bf Proof.} We restrict ourselves by the case of the good component with thivial
signature of $T(\tilde \s)$; proof for other components is analogous. So, we can take
$s=0$ and $\s=\left( \matrix E_g&0 \\0&pE_g\endmatrix \right)$. We apply formula 4.10 of
[L01] for the degree of $\pi_2$. We set $\goth g=h+k\rd$ like above ($g$ of [L01], 4.10
is $\alpha(\goth g)$ in the notations of the present paper). So, the degree of $\pi_2$ is
the quantity of $\s'\in S(g)$ such that there exists $\goth g$ such that

$$\s'\s\alpha(\goth g) \in p \cdot GSp_{2g}(\n Z)\eqno{(3.4.1.1)}$$

Particularly, for $\goth g$ satisfying (3.4.1.1) we have

$$\ord_p(\det \goth g)=0\eqno{(3.4.1.2)}$$

We use description of $S(g)$ from [L04.1], Section 2.3 ($\s'$ is denoted in [L04.1],
Section 2.3 by $\gamma$). Let us fix $\s'$. Attached to $\s'$ is a subset $I$ of $\{1,
\dots, g\}$.

Let $j\in \{1, \dots, g\}$, $j \not\in I$. Multiplying the $j$-th line of $\s'\s$ by the
$r$-th (resp. $g+r$-th) column of $\alpha(\goth g)$ we get that (3.4.1.1) implies
$h_{jr}\in \n Z$ (resp. $k_{jr}\in \n Z$). Analogously, for $i \in I$, multiplying the
$g+i$-th line of $\s'\s$ by columns of $\alpha(\goth g)$, we get that $h_{ir}, k_{ir} \in
\frac 1p \n Z$. Further, multiplying the $i$-th line of $\s'\s$ by columns of
$\alpha(\goth g)$, we see that all terms of the corresponding scalar product --- except
one --- are integer, hence this only remaining term --- which is $h_{ir}$ or $k_{ir}$ ---
is also integer.

Let us consider the case $I\ne \emptyset$, and let $i$ be the minimal element of $I$.
Multiplying the $i$-th line of $\s'\s$ by all columns of $\alpha(\goth g)$ we get that
(3.4.1.1) implies that the first $i$ lines of $\goth g$ are linearly dependent mod $p$.
This contradicts to (3.4.1.2). So, we get that if for the given $\s'$ there exists $\goth
g$ such that (3.4.1.1) holds then $I=\emptyset$. There exists only one such $\s'=\left(
\matrix pE_g&0 \\0&E_g\endmatrix \right)$, and we can take $\goth g=1$. $\square$

\medskip

Now let us consider the case of the bad component for $g=3$.

\medskip

{\bf Remark 3.4.2.} To use the formula 4.10 of [L01], we need to check all elements of
$S(g)$, not only of its open part. This explains why the proof of the following
Proposition 3.4.3 is so long: I am forced to treat all 8 open components of $S(g)$
separately. I do not see a method to find a uniform method of theating. For higher $g$
and a partition set contained in $\goth D_i$ it is possible to treat some cases
uniformly, so the quantity of cases is a polynomial of $g$, $i$.

We shall see (Remark 3.4.4) that the set of $\s'\in S(3)$ satisfying (3.4.1.1) is an
irreducible subvariety of $S(3)=G_I(3,6)$. We can expect that this is true for all $g$,
$i$. Proof of this fact will reduce the quantity of cases that we shall have to consider
in order to find the degree of $\pi_{2}$.

\medskip

{\bf Proposition 3.4.3.} For $g=3$ the degree of $\pi_{2}$ for the bad component is $p+1$.

\medskip

{\bf Proof.} We use the same formula 4.10 of [L01]. We fix $\s=\left( \matrix E_g&W
\\0&pE_g\endmatrix \right)$, where $s=W=\left( \matrix 0&0&0
\\0&w_1&w_2\\0&w_2&w_1\endmatrix \right)$, $w_1,w_2\in \n Z$ satisfy $w_1^2-w_2^2\equiv
\frac{-1}{\d} \mod p^2$ (this implies that $s$ corresponds to the bad component). It is
more convenient to take $\goth g=\frac 1p (h+k\rd)$. Let $\s'=\left( \matrix
A&B\\0&D\endmatrix \right)$ and $I$ be as in [L04.1], Section 2.3. We denote
$P=\s'\s\alpha(\goth g)$. Although we prove this theorem for $g=3$, some arguments are
valid for any $g$, so we shall use sometimes $g$ instead of 3.

The entries of $i$-th line of $\s'\s$ are multiples of $p$ for $i \not\in I$,
$i\in\{1,\dots,g\}$, and for all $i\in\{g+1,\dots,2g\}$. For this set of $i$'s conditions
$P_{ij}\in p\n Z$ can be treated as linear congruences modulo $p$ on elements $h_{\alpha
j}$, $k_{\beta j}$:

$$P_{ij}\in p\n Z \iff \sum_{\alpha=1}^g (C_{11})_{i \alpha} \ h_{\alpha j} +
\sum_{\alpha=1}^g (C_{12})_{i \alpha} \ k_{\alpha j} \equiv 0 \hbox{ (for $i\not\in I$,
$i\in\{1,\dots,g\}$)}$$

$$P_{ij}\in p\n Z \iff \sum_{\alpha=1}^g (C_{21})_{i \alpha} \ h_{\alpha j} +
\sum_{\alpha=1}^g (C_{22})_{i \alpha} \ k_{\alpha j} \equiv 0 \hbox{ (for
$i\in\{g+1,\dots,2g\}$)}$$

$$P_{i,g+j}\in p\n Z \iff \sum_{\alpha=1}^g (C_{31})_{i \alpha} \ h_{\alpha j} +
\sum_{\alpha=1}^g (C_{32})_{i \alpha} \ k_{\alpha j} \equiv 0 \hbox{ (for $i\not\in I$,
$i\in\{1,\dots,g\}$)}$$

$$P_{i,g+j}\in p\n Z \iff \sum_{\alpha=1}^g (C_{41})_{i \alpha} \ h_{\alpha j} +
\sum_{\alpha=1}^g (C_{42})_{i \alpha} \ k_{\alpha j} \equiv 0 \hbox{ (for
$i\in\{g+1,\dots,2g\}$)}$$

where $j\in\{1,\dots,g\}$ and the coefficient matrix $C$ has the $4\times 2$-block
structure, $(C_{uv})_{xy}$ is the $(x,y)$-th element of the $(u,v)$-block of $C$. It is
easy --- but tedious --- to write down expressions for all 8 blocks of $C$ (for example,
$C_{21}=C_{42}=0$, $C_{22}=-\d D E_{\goth r \goth s}$, $C_{41}=D E_{\goth r \goth s}$
etc).

Looking at the second (resp. fourth) block line of $C$ we get immediately that the
corresponding congruences imply that all entries of $k$ (resp. $h$) are $p$-integer. This
means that

\medskip

{\bf (3.4.3.1)} If the rank of $C$ modulo $p$ is $2g$ then all entries of $h$, $k$ are 0
modulo $p$.

\medskip

{\bf Case $I=\{1,2,3\}$.} $\s'=\left( \matrix 1&B_3 \\0&p\endmatrix \right)$, where
$s=B_3=\{b_{ij}\}$ is any symmetric $3\times 3$-matrix. We shall use also the (1,2)-block
partitions of $B_3$, $W$, $h$, $k$: $B_3=\left( \matrix
B_{11}&B_{12}\\B_{12}^t&B_{22}\endmatrix \right)$, $W=\left( \matrix 0&0\\0&W_2\endmatrix
\right)$, $h=\left( \matrix \goth h_{11} & \goth h_{12}
                        \\ \goth h_{21} & \goth h_{22} \endmatrix \right)$ and
analogously for $k$. Calculating explicitly entries of $P$ we see that the condition
(3.4.1.1) gives us immediately congruences:

\medskip

3.4.3.2) $h\equiv \Delta(W+pB_3)E_{21} k \mod p^2$

\medskip

3.4.3.3) $k\equiv -(W+pB_3)E_{21} h \mod p^2$

\medskip

Substituting $h$ from (3.4.3.2) to (3.4.3.3) we get

\medskip

3.4.3.4)  $[1_3+\Delta WE_{21}WE_{21}+\Delta p(B_3E_{21}W+WE_{21}B_3)E_{21}]k\equiv 0
\mod p^2$ which is in the above notations

\medskip

3.4.3.5) $\goth k_{11}+\Delta p B_{12}E_{11} W_2E_{11} \goth k_{21}\equiv 0 \mod p^2$

\medskip

3.4.3.6) $\goth k_{12}+\Delta p B_{12}E_{11} W_2E_{11} \goth k_{22}\equiv 0 \mod p^2$

\medskip

3.4.3.7) $\Delta p W_2E_{11} B_{12}^t \goth k_{11}+\Delta p (B_{22}E_{11} W_2 + W_2E_{11}
B_{22})E_{11}\goth k_{21}\equiv 0 \mod p^2$

\medskip

3.4.3.8) $\Delta p W_2E_{11} B_{12}^t \goth k_{12}+\Delta p (B_{22}E_{11} W_2 + W_2E_{11}
B_{22})E_{11}\goth k_{22}\equiv 0 \mod p^2$

\medskip

Substituting value of $\goth k_{11}$ from (3.4.3.5) to (3.4.3.7) (resp. $\goth k_{12}$
from (3.4.3.6) to (3.4.3.8)), we get that (3.4.3.7), (3.4.3.8) become

\medskip

3.4.3.9) $(B_{22}E_{11} W_2 + W_2E_{11} B_{22})E_{11}\goth k_{2i}\equiv 0 \mod p$,
$i=1,2$.

\medskip

If $\det (B_{22}E_{11} W_2 + W_2E_{11} B_{22})\not\equiv 0 \mod p$, then (3.4.3.9) imply
that $\goth k_{2i}\equiv 0 \mod p$. Further, (3.4.3.5), (3.4.3.6) imply that $\goth
k_{1i}\equiv 0 \mod p^2$, (3.4.3.2) implies that $\goth h_{1i}\equiv 0 \mod p^2$,
$h\equiv 0 \mod p$. This contradicts to (3.4.1.1).

\medskip

Now we need a

\medskip

{\bf Lemma 3.4.3.10.} If $\det (B_{22}E_{11} W_2 + W_2E_{11} B_{22})\equiv 0 \mod p$ then
$B_{22}E_{11} W_2 + W_2E_{11} B_{22}\equiv 0 \mod p$ and $B_{22}=\gamma \left( \matrix
w_2&w_1\\w_1&w_2\endmatrix \right)$ where $\gamma$ is a scalar factor.

\medskip

{\bf Proof.} Direct calculation. Really, condition $\det (B_{22}E_{11} W_2 + W_2E_{11}
B_{22})\equiv 0 \mod p$ in terms of $b_{ij}$ is

$$[- \Delta w_2^2 b_{22} +(2 \Delta w_1w_2 +2 w_2 \rd)b_{23} +(\Delta w_2^2-2\Delta w_1^2
-2 w_1 \rd)b_{33}]\cdot$$

$$[- \Delta w_2^2 b_{22} +(2 \Delta w_1w_2 -2 w_2 \rd)b_{23} +(\Delta w_2^2-2\Delta w_1^2
+2 w_1 \rd)b_{33}]=0$$

this is the union of 2 imaginary lines on the projective plane $(b_{22}:b_{23}:b_{33})$
whose real intersection point corresponds to the case $B_{22}E_{11} W_2 + W_2E_{11}
B_{22}\equiv 0 \mod p$. Solving the corresponding system, we get the above expression for
$B_{22}$. $\square$

\medskip

In this case (3.4.3.7), (3.4.3.8) are always satisfied, substituting values of $\goth
k_{1i}$ from (3.4.3.5), (3.4.3.6) to (3.4.3.2) we get (congruences are modulo $p^2$):

\medskip

$h\equiv \left( \matrix \Delta p B_{12}E_{11} \goth k_{21} & \Delta p B_{12}E_{11} \goth
k_{22}
\\ \Delta (W_2 + p B_{22})E_{11}\goth k_{21} & \Delta (W_2 + p B_{22})E_{11}\goth k_{22}
\endmatrix \right)$ and hence

$$p g\equiv \left( \matrix \Delta p B_{12}E_{11}(1_2-\rd W_2E_{11}) \goth k_{21} & \Delta
p B_{12} E_{11}(1_2-\rd W_2E_{11}) \goth k_{22}
\\ (\Delta (W_2 + p B_{22})E_{11} + \rd \cdot 1_2) \goth k_{21} & (\Delta (W_2 + p
B_{22})E_{11} + \rd \cdot 1_2) \goth k_{22}
\endmatrix \right)\eqno{(3.4.3.11)}$$

Now we check the condition $(p\goth g)E_{21} (p\bar \goth g^t) = \lambda E_{21}$, where
$\ord_p(\lambda)=2$. Recall that $b_{ij}$, $k_{ij}$ are entries of $B_3$, $k$
respectively. Diagonal entries of $(p\goth g)E_{21} (p\bar \goth g^t)$ are elements of
the ring of polynomials $\n Z[w_1, w_2, \d, b_{ij}, k_{ij}]$ factored by
$\d(w_1^2-w_2^2)+1=0$. A calculation in this ring shows that $$\hbox{if $B_{22}=0$ then
$((p\goth g) E_{21} (p\bar \goth g^t))_{11} = Q p^2 ((p\goth g) E_{21} (p\bar \goth
g^t))_{33}$} \eqno{(3.4.11')}$$ for some $Q\in\n Z[w_1, w_2, \d, b_{ij},
k_{ij}]/(\d(w_1^2-w_2^2)+1)$. 

Since $B_{22}$ enters in (3.4.3.11) with a coefficient $p$, we see that the condition
$\ord_p((p\goth g) E_{21} (p\bar \goth g^t)_{11})=2$ together with ($3.4.11'$)
contradicts to the condition $\ord_p((p\goth g) E_{21} (p\bar \goth g^t)_{33})=2$. This
means that for $\s'$ of type $I=\{1,2,3\}$ there is no $\goth g$ satisfying (3.4.1.1).

\medskip

{\bf Case $I=\{2,3\}$.} $\s'=\left( \matrix p & 0 & 0 & 0
             \\ -D^t & 1 & 0 & B
             \\    0 & 0 & 1 & D
             \\    0 & 0 & 0 & p \endmatrix \right)$ where $D$, $B$ are respectively
$1\times 2$, $2\times 2$-matrices, $B$ is symmetric. The condition that the entries of
the first and the third block lines of $\s'\s\alpha(\goth g)$ are in $p \n Z$ gives us
immediately that entries of $\goth h_{1i}$, $\goth k_{1i}$ are in $p \n Z$ ($i=1,2$).

\medskip

{\bf Subcase 1.} $D=0$. The condition that the entries of the second block line of
$\s'\s\alpha(\goth g)$ are in $p \n Z$ becomes

$$\goth h_{2i} \equiv \Delta (W_2+pB)E_{11}\goth k_{2i} \mod p^2\eqno{(3.4.3.12)}$$

$$\goth k_{2i} \equiv - (W_2+pB)E_{11}\goth h_{2i} \mod p^2\eqno{(3.4.3.13)}$$

or $[E_{11} + \Delta W_2E_{11}W_2 + p \Delta (W_2E_{11}B+BE_{11}W_2))]E_{11} \goth h_{2i}
\equiv 0$  mod $p^2$, i.e. $$(W_2E_{11}B+BE_{11}W_2)E_{11} \goth h_{2i} \equiv 0 \mod
p\eqno{(3.4.3.13')}$$

{\bf Subcase 1a.} $B\ne \gamma \left( \matrix w_2&w_1\\w_1&w_2\endmatrix \right)$. Using
Lemma 3.4.3.10 we get that ($3.4.3.13'$) implies $\goth h_{2i}, \goth k_{2i} \equiv 0$ 
mod $p$. (3.4.3.12), (3.4.3.13) become $\goth k_{2i} \equiv - W_2E_{11} \goth h_{2i}$,
$\goth h_{2i} \equiv \d W_2E_{11} \goth k_{2i}$ mod $p^2$. This means that
$(1_2+W_2E_{11}\d)\goth g_{2i}\equiv 0 \mod p$, i.e. that the second and the third lines
of $\goth g$ are linearly dependent. This contradicts to (3.4.1.1).

\medskip

{\bf Subcase 1b.} $B =\gamma \left( \matrix w_2&w_1\\w_1&w_2\endmatrix \right)$. For
these $\s'$ there exist $\goth g$ satisfying (3.4.1.1). Really, if we take any $\goth
h_{1i}$, $\goth k_{1i}$ with entries in $p \n Z$, any $\goth h_{2i}$ with entries in $\n
Z$

and $\goth k_{2i}$ satisfying (3.4.3.13) then entries of $\s'\s\alpha(\goth g)$ are in $p
\n Z$. Now we must find $\goth h_{1i}$, $\goth k_{1i}$, $\goth h_{2i}$ such that $\goth
g$ satisfies (3.1.1). This can be done by many ways. For example, we can take $\goth
g_{12}=\goth g_{21}=0$ and $\goth h_{22}$ satisfying $\goth h_{22}E_{11}\goth
h_{22}^t=\mu E_{11}$ for some $\mu$. Simple calculations show that for all $\gamma$
(3.1.1) can be satisfied.

\medskip

{\bf Subcase 2.} $D\ne 0$. As earlier, elementary transformations of congruences
$(\s'\s\alpha(\goth g))_{ij}\equiv 0 \mod p$ give us the desired. Firstly, the condition
that the entries of the third block line of $\s'\s\alpha(\goth g)$ are in $p \n Z$
implies that

$$D E_{11}\goth h_{2i} \equiv 0 \mod p, \ \ D E_{11}\goth k_{2i} \equiv 0 \mod
p\eqno{(3.4.3.14)}$$

These conditions together with the condition that the entries of the second block line of
$\s'\s\alpha(\goth g)$ are in $p \n Z$ implies that

$$D E_{11}W_2E_{11} \goth h_{2i} \equiv 0 \mod p, \ \ D E_{11}W_2E_{11} \goth k_{2i}
\equiv 0 \mod p\eqno{(3.4.3.15)}$$

It is easy to check that $W_2E_{11}$ has no eigenvectors in $\n F_p$ (because $\left(
\frac {-\Delta} p \right)=-1$), so (3.4.3.14), (3.4.3.15) imply that entries of $\goth
h_{2i}$, $\goth k_{2i}$ are in $p \n Z$.

Now the condition that the entries of the third block line of $\s'\s\alpha(\goth g)$ are
in $p \n Z$ becomes

$$-D^t \goth h_{1i} + \goth h_{2i} - \Delta W_2E_{11} \goth k_{2i} \equiv 0 \mod p^2
\eqno{(3.4.3.16)}$$

$$D^t \goth k_{1i} - \goth k_{2i} - W_2E_{11} \goth h_{2i} \equiv 0 \mod p^2
\eqno{(3.4.3.17)}$$

Elementary transformations reduce this system to the system

$$\left( \matrix -w_1 d_1 + w_2d_2 & d_1 \\ -w_2 d_1 + w_1d_2 & d_2 \endmatrix \right)
\left( \matrix \goth h_{1i}\\ \goth k_{1i}\endmatrix \right) \equiv 0 \mod p^2$$ (here
$d_i$ are entries of $D$).

Since det  $\left( \matrix -w_1 d_1 + w_2d_2 & d_1 \\ -w_2 d_1 + w_1d_2 & d_2 \endmatrix
\right)$ is never 0, we get that $\goth h_{1i}, \ \goth k_{1i} \equiv 0 \mod p^2$ --- a
contradiction to (3.4.1.2). So, for this case there is no $\goth g$ satisfying (3.4.1.1).

\medskip

{\bf Case $I=\{1\}$.} $\s'=\left( \matrix 1 & 0 & 0 & b_{11} & 0 & 0
                          \\ 0 & p & 0 & 0 & 0 & 0
                          \\ 0 & 0 & p & 0 & 0 & 0
                          \\ 0 & 0 & 0 & p & 0 & 0
                          \\ 0 & 0 & 0 & 0 & 1 & 0
                          \\ 0 & 0 & 0 & 0 & 0 & 1 \endmatrix \right)$. (3.4.1.1) gives
us immediately that the second and the third lines of $h,k$ are $\equiv 0 \mod p$, and
their first lines are $\equiv 0 \mod p^2$. This contradicts to (3.4.1.2).

\medskip

{\bf Case $I=\{2\}$.} $\s'=\left( \matrix p & 0 & 0 & 0 & 0 & 0
                         \\ -d_{12} & 1 & 0 & 0 & b_{22} & 0
                          \\ 0 & 0 & p & 0 & 0 & 0
                          \\ 0 & 0 & 0 & 1 & d_{12} & 0
                          \\ 0 & 0 & 0 & 0 & p & 0
                          \\ 0 & 0 & 0 & 0 & 0 & 1 \endmatrix \right)$. Writing
explicitly the matrix $C$ for this case, we see immediately that its rank is 6, i.e.
(3.4.3.1) implies that $h,k\equiv 0 \mod p$. Further, the condition $[\s'\s\alpha(\goth
g)]_{2i}\in p\n Z$ implies that

$$\left( \matrix d_{12} & -1- w_1 \rd & w_2 \rd\endmatrix \right)\goth g\equiv 0 \mod p$$
This contradicts to (3.4.1.2).

\medskip

{\bf Case $I=\{3\}$.} $\s'=\left( \matrix    p &    0 & 0 & 0 & 0 & 0
                          \\    0 &    p & 0 & 0 & 0 & 0
                          \\ -d_{13} & -d_{23} & 1 & 0 & 0 & b_{33}
                          \\    0 &    0 & 0 & 1 & 0 & d_{13}
                          \\    0 &    0 & 0 & 0 & 1 & d_{23}
                          \\    0 &    0 & 0 & 0 & 0 & p \endmatrix \right)$. The non-0
lines of $C$ modulo $p$ form the following matrix $C'$:

$C'=\left( \matrix              1 &   0 &    0 &   0 &       0 & 0
                          \\    0 &   1 &    0 &   0 & -\d w_1 & \d w_2
                          \\    0 &   0 &    0 & -\d &       0 & \d d_{13}
                          \\    0 &   0 &    0 &   0 &     -\d & \d d_{23}
                          \\    0 &   0 &    0 &   1 &       0 & 0
                          \\    0 & w_1 & -w_2 &   0 &       1 & 0
                          \\    1 &   0 & -d_{13} &   0 &       0 & 0
                          \\    0 &   1 & -d_{23} &   0 &       0 & 0
\endmatrix \right)$

If $d_{13}\not\equiv 0 \mod p$ then it is seen immediately that the rank of $C'$ is 6; in
any case, the determinant of the submatrix of $C'$ formed by removing of its third and
seventh lines is $\pm w_2(w_2d_{23}^2+2w_1d_{23}-w_2)$ which is never 0 because of
$w_1^2-w_2^2\equiv \frac{-1}{\d} \mod p^2$ and because $p$ is inert in $K$.

We denote the third line of (1,1)-block (resp. of (1,2)-block) of $\s'\s$ by $l_1$ (resp.
by $l_2$). It is easy to check that the condition $[\s'\s\alpha(\goth g)]_{3i}\in p\n Z$
implies that

$$(l_1 + \rd l_2)E_{21}\goth g\equiv 0 \mod p$$ (like in the case $I=\{2\}$). This
contradicts to (3.4.1.2).

\medskip

{\bf Cases $I=\{1,3\}$ and $I=\{1,2\}$.} $\s'=\left( \matrix             1 &    0 & 0 &
b_{11} & 0 & b_{13}
                          \\    0 &    p & 0 & 0      & 0 & 0
                          \\    0 &   -d_{23} & 1 & b_{13} & 0 & b_{33}
                          \\    0 &    0 & 0 & p      & 0 & 0
                          \\    0 &    0 & 0 & 0      & 1 & d_{23}
                          \\    0 &    0 & 0 & 0      & 0 & p \endmatrix \right)$ and
$\s'=\left( \matrix        1 &    0 & 0 & b_{11} & b_{12} & 0
                          \\    0 &    1 & 0 & b_{12} & b_{22} & 0
                          \\    0 &    0 & p & 0      & 0      & 0
                          \\    0 &    0 & 0 & p      & 0      & 0
                          \\    0 &    0 & 0 & 0      & p      & 0
                          \\    0 &    0 & 0 & 0      & 0      & 1 \endmatrix \right)$
respectively. The non-0 lines of $C$ are

$\left( \matrix 0&1&0&0&-\d w_1&\d w_2 \\
                0&w_1&-w_2&0&1&0\\
                0&0&0&0&-\d&\d d_{23}\\
                0&1&-d_{23}&0&0&0\endmatrix \right)$,

$\left( \matrix 0&0&1&0&-\d w_2&\d w_1 \\
                0&w_2&w_1&0&0&1\\
                0&0&0&0&0&\d\\
                0&0&-1&0&0&0\endmatrix \right)$ respectively. Both $C$ are of rank 4 for
all $\s'$, i.e. $h_{2i}$, $h_{3i}$, $k_{2i}$, $k_{3i}$ are in $p\n Z$. Considering the
condition that the $(1,i)$-th and the $(1,i+g)$-th element of $\s'\s\alpha(\goth g)$ are
integer, we get immediately that $h_{1i}$, $k_{1i}$ are in $p\n Z$. Considering again the
same condition, we get that moreover $h_{1i}$, $k_{1i}$ are in $p^2\n Z$. This is
sufficient to get a contradiction.

\medskip

{\bf Case $I=\emptyset$.} $\s'=\left( \matrix p & 0
                          \\ 0 & 1 \endmatrix \right)$. This is the trivial case, $\goth
g=1$ satisfies (3.4.1.1).

\medskip

We see that in $S(3)$ there are $p$ matrices of type $I=\{2,3\}$, Subcase 1b, and 1
matrix of type $I=\emptyset$ such that there exists $\goth g$ satisfying (3.4.1.1). This
gives us our result. $\square$

\medskip

{\bf Remark 3.4.4.} It is easy to check that the above $p+1$ matrices form a projective
line in $S(3)=G_I(3,6)\subset P^{19}$.

\medskip

{\bf Section 4. Case of Hecke correspondence $\goth T_p=T_{p,\goth i}$. } \nopagebreak

\medskip

{\bf Subsection 4.1.} Some matrix equalities.  \nopagebreak

\medskip

Here we formulate and prove some matrix equalities that will be necessary for the next
subsection. The open part $S(g, \goth i)^{open}\subset S(g, \goth i)$ is the set of the
following block matrices (diagonal blocks have sizes $\goth i$, $g-\goth i$, $\goth i$,
$g-\goth i$):

$$\s=\left (\matrix p& 0&0&p C^t\\
-A&1&C&U\\ 0&0&p&pA^t\\ 0&0&0&p^2 \endmatrix\right),\eqno{(4.1.1)}$$

where $A,C$ are $\goth i\times(g-\goth i)$-matrices, $U$ is a $(g-\goth i)\times(g-\goth
i)$-matrix with entries in $\n Z$. Moreover, entries of $A,C$ run over a fixed system of
representatives of $\n Z$ modulo $p\n Z$, diagonal and upper-triangular entries of $U$
run over a fixed system of representatives of $\n Z$ modulo $p^2\n Z$, and
lower-triangular entries of $U$ are defined uniquely by the relation

$$CA^t - AC^t = U-U^t\eqno{(4.1.2)}$$

which is equivalent to the condition $\s\in GSp_{2g}$.

Let us define some objects associated to $\s$. Firstly, we define $\mu_1=\mu_1(\s)$,
$\mu_2=\mu_2(\s)\in M_g(\n Z(\rd))\hookrightarrow M_g(\n Z_{p^2})$ (here and below
diagonal blocks have sizes $\goth i$, $g-\goth i$):

$$\mu_1=\left (\matrix -1&-C^t\rd \\ A+C\rd&-1-U\rd \endmatrix\right),\mu_2=\left
(\matrix \rd &-A^t\rd \\ A+C\rd&-1-U\rd \endmatrix\right)\eqno{(4.1.2.5)}$$

To simplify notations, we consider only the case $(\goth r,\goth s)=(\goth i, g-\goth
i)$, and we denote $E_{\goth r \goth s}=E$. Further, we define
$G=\mu_{1}E\bar\mu_{1}^{t}$. Let the block structure of $G$ be $\left (\matrix G_{11} &
G_{12} \\ G_{21} & G_{22}  \endmatrix\right)$. Finally, if $\det \mu_1(\s)$ is invertible
in $\n Z_{p^2}$ then we define matrices $X_1$, $X_2$ by the formula $$\mu_{2}=\left
(\matrix X_{1} & X_{2}\\0&1 \endmatrix\right)\mu_{1}\eqno{(4.1.3)}$$

and $F=G^{-1}=(\mu_{1}E\bar\mu_{1}^{t})^{-1}=\left (\matrix F_{11} & F_{12} \\ F_{21} &
F_{22}  \endmatrix\right)$.

\medskip

{\bf Remark 4.1.4.} Matrices $\mu_1$, $\mu_2$ play symmetric roles in the contents of the
present subsection, so it is possible to rewrite it for the case $\det \mu_2(\s)$ is
invertible in $\n Z_{p^2}$. Recall that if $\s\not\in\goth D_*$ (see (4.2.9)) then either
$\det \mu_1(\tilde \s)$, or $\det \mu_2(\tilde \s)$, or both of them are $\ne 0$. So,
when we shall consider in the sequel an element $\s\not\in\goth D_*$, we shall assume
always that $\det \mu_1(\tilde \s)\ne0$. If not, then we interchange roles of
$\mu_1(\s)$, $\mu_2(\s)$ and we get the proof for the case $\det \mu_1(\tilde \s)=0$,
$\det \mu_2(\tilde \s)\ne0$.

\medskip
We shall apply terminology "real", "imaginary" etc. to the extensions $\n Q_p \hookrightarrow \n Q_{p^2}$, $\n F_p \hookrightarrow \n F_{p^2}$ in the obvious sense. 
\medskip
{\bf Proposition 4.1.5.} $G$ is real (and hence symmetric) matrix.

\medskip

{\bf Proof.} Follows immediately from the definition of $G$ and (4.1.2). $\square$

\medskip

Clearly the same is true for $F$.

\medskip

{\bf Proposition 4.1.6.} $\im X_1=F_{11}$, $\im X_2=F_{12}$ (and hence $F_{21}=\im
X_2^t$).

\medskip

{\bf Proof.} Writing down the real and the imaginary parts of the (1,1)- and (1,2)-block
entries of the equality (4.1.3) we get 4 equalities:

$$- \re X_1 + \re X_2 A - \Delta \im X_2 C=0 \eqno{(4.1.6.1r)}$$

$$-\im X_{1}+\re X_{2} C + \im X_{2} A=1\eqno{(4.1.6.1i)}$$

$$\Delta \im X_{1} C^t -\re X_{2} + \Delta \im X_{2} U=0\eqno{(4.1.6.2r)}$$

$$- \re X_1 C^t - \re X_2 U - \im X_2 =-A^t \eqno{(4.1.6.2i)}$$

where the last 2 symbols $r$ (resp. $i$); 1 (resp. 2) in the number of the equality mean
that it comes from the equality for the real (resp. imaginary) part of (1,1) (resp.
(1,2))-blocks of the equality (4.1.3).

Now we eliminate $\re X_1$, $\re X_2$ from these equalities (firstly we find $\re X_2$
from $(4.1.6.2r)$ and secondly $\re X_1$ from $(4.1.6.1r)$; we get

$$\im X_1 (-1 + \Delta C^t C) + \im X_2 (A+\Delta UC)=1\eqno{(4.1.6.3)}$$

$$\im X_1 (-\Delta C^t U - \Delta C^t AC) + \im X_2 (-1-\Delta UAC^t+\Delta CC^t - \Delta
U^2)=-A^t\eqno{(4.1.6.4)}$$

It is sufficient to check that

$$\im X_{1}\hbox{ } G_{11} + \im X_{2}\hbox{ } G_{21}=1\eqno{(4.1.6.5)}$$

$$\im X_{1}\hbox{ } G_{12} + \im X_{2}\hbox{ } G_{22}=0\eqno{(4.1.6.6)}$$

Writing down explicitly $G_{jk}$ we get immediately that (4.1.6.5) = (4.1.6.3) and
(4.1.6.6) = (4.1.6.3) $\cdot A^t$ + (4.1.6.4). $\square$

\medskip

{\bf Corollary 4.1.7.} $\det\im \tilde X_{1}\ne0\iff \det \tilde G_{22}\ne0$. $\square$

\medskip

{\bf 4.1.8.} If $\det \im \tilde X_{1}\ne0$ then we define $W=-F_{11}^{-1}F_{12}=-(\im
X_{1})^{-1}\im X_{2}$.

\medskip

{\bf Corollary 4.1.9.} In this case we have for any $g$, $\goth i$:

$$F_{11}W + F_{12}=0;\eqno{(4.1.9.1)}$$

$$WG_{22}=G_{12}\eqno{(4.1.9.2)}$$

$$F_{21}W+F_{22}=G_{22}^{-1}\eqno{(4.1.9.3)}$$

$$\det G=\det G_{22}\det (G_{11}-WG_{22}W^t)\eqno{(4.1.9.4)}$$

{\bf Proof.} (4.1.6.6) (resp. (4.1.6)) and the definition of $W$ imply immediately
(4.1.9.2) (resp. (4.1.9.1)). Further, (4.1.9.3) follows from a general matrix formula:
let $\goth F=\left (\matrix \goth F_{11} & \goth F_{12} \\ \goth F_{21} & \goth F_{22} 
\endmatrix\right)$, $\goth G=\goth F^{-1}=\left (\matrix \goth G_{11} & \goth G_{12} \\
\goth G_{21} & \goth G_{22}  \endmatrix\right)$ be arbitrary block matrices, $\goth
F_{11}$ is invertible and $\goth W=-\goth F_{11}^{-1}\goth F_{12}$. So, $\goth
F_{21}\goth W+\goth F_{22}=\goth G_{22}^{-1}$. Analogously, (4.1.9.4) also holds for any
$\goth F$, $\goth G=\goth F^{-1}$. $\square$

\medskip

Now we consider the case when $\det \im \tilde X_{1}=0$.

\medskip

{\bf Lemma 4.1.10.} Let $D\ne0$ be a matrix with entries in $\n F_p$ such that $D\im
\tilde X_{1}=0$. Then $D\im \tilde X_{2}\ne0$. Particularly, $\det \im \tilde
X_{1}=0\hbox{ implies } \im \tilde X_{2}\ne0$.

\medskip

{\bf Proof.} We reduce $(4.1.6.1i)$, $(4.1.6.2r)$ modulo $p$, multiply them by $D$ from
the left and eliminate $D\re \tilde X_2$. $D\im \tilde X_{1}=0$ implies $D\im \tilde
X_{2}(\tilde A+\Delta \tilde U \tilde C))=D$, this contradicts to conditions $D \im
\tilde X_{2}=0$, $D\ne0$. $\square$

\medskip

For a vector row $X=(x_1,x_2)$ we denote by $X^O$ the orthogonal vector $(-x_2, x_1)$,
and for a vector column $X=\left (\matrix x_1\\ x_2\endmatrix\right)$ $X^O$ will mean
$\left (\matrix -x_2\\ x_1\endmatrix\right)$.

\medskip

{\bf Proposition 4.1.11.} Let $g=3$, $\goth i=1$, $\tilde F_{11}=0$. Then

$$\tilde F_{12}^O \tilde F_{22} \tilde F_{12}^{Ot}=-\det \tilde F\eqno{(4.1.11.1)}$$

(we identify a $1\times 1$-matrix with a number);

$$\tilde F_{21}^O \tilde F_{21}^{Ot}=\frac{-1}{\det \tilde F} \tilde
G_{22}\eqno{(4.1.11.2)}$$

$$\tilde F_{22} \tilde F_{12}^{Ot}=-\det \tilde F \ \tilde G_{21}^{O}\eqno{(4.1.11.3)}$$

$$\hbox{The determinant of the matrix formed by vectors $\tilde F_{21}$, $\tilde
G_{21}^O$ is 1.}\eqno{(4.1.11.4)}$$

\medskip

{\bf Proof.} All these equalities hold for any symmetric invertible $3\times 3$-matrix
$F$ having $F_{11}=0$. $\square$

\medskip

Now we need a lemma for the case when both $\det \tilde \mu_{1}=\det \tilde \mu_{2}= 0$.
The following simple proof of this lemma is due to A. Zelevinskij.

\medskip

{\bf Lemma 4.1.12.} Let $g$ is arbitrary, $\goth i=1$. If $\det \tilde \mu_{1}=\det
\tilde \mu_{2}= 0$ then the $(g-1)\times g$-matrix $\left (\matrix \tilde A+\tilde
C\rd&-1-\tilde U\rd \endmatrix\right)$ has rank $\le g-2$.

\medskip

{\bf Proof.} We consider the following $(g+1)\times g$-matrix $\mu$:

$\left (\matrix -1&-\tilde C^t\rd \\  \rd&-\tilde A^t\rd \\ \tilde A+\tilde
C\rd&-1-\tilde U\rd \endmatrix\right)$

For each $i=1,\dots, g-1$ we consider 2 minors $m_{i1}$, $m_{i2}$ of this matrix:

$m_{i1}$ is formed by all lines of $\mu$, except the $i$-th line of the block $\left
(\matrix \tilde A+\tilde C\rd&-1-\tilde U\rd \endmatrix\right)$, and by all columns of
$\mu$;

$m_{i2}$ is formed by all lines of the block $\left (\matrix \tilde A+\tilde
C\rd&-1-\tilde U\rd \endmatrix\right)$, and by all columns of $\mu$, except the $i$-th
column of the block $-1-\tilde U\rd$.

Elementary transformations show that (taking into consideration 4.1.2) for all $i$ we
have $\det m_{i1}= \det m_{i2}$. This implies immediately the lemma. $\square$

\medskip

{\bf Subsection 4.2.} Finding of D-equivalences of some $\s_k$. \nopagebreak

\medskip

Let $\s_1$, $\s_2\in S(g,\goth i)^{open}$, i.e. they are matrices of type (4.1.1). Let us
consider objects associated to $\s$ defined in Section 4.1, namely $A$, $U$, $C$,
$\mu_{i}$, $G$, $G_{ij}$, $F$, $F_{ij}$, $X_{i}$ ($i,j=1,2$). These objects for the above
$\s_k$ ($k=1,2$) will be denoted by $A_k$, $U_k$, $C_k$, $\mu_{ik}$, $G_k$, $G_{ijk}$,
$F_k$, $F_{ijk}$, $X_{ik}$ respectively. Recall that $\s_1$, $\s_2$ are D-equvalent iff
$\exists \goth g=h+k\rd\in G_V(\n Z_p)$ such that $\s_1\alpha(\goth g)\s_2^{-1}$ has
integer entries.

Analog of (3.3.1.2) is the following:

$$\s_1\alpha (\goth g)\s_2^{-1}=\left(\matrix * & * & * & p^{-1} A_{14}\\p^{-1} A_{21} &
* & p^{-1} A_{23} & p^{-2} A_{24} \\ * & * & * & p^{-1} A_{34}\\ * & * & * & *
\endmatrix\right)\eqno{(4.2.1)}$$

where *'s and $A_{**}$ are some polynomials with integer coefficients in $h,k$, $A_i$,
$U_i$, $C_i$ ($i=1,2$) and hence $\in M_g(\n Z_p)$.

Analog of $t(\s)$ of (3.2.2) is the pair $\mu_{1}(\s)$, $\mu_{2}(\s)$. For $i=1,2$ we
denote

$$B_i=B_i(\s_1,\s_2)=\mu_{i1}\goth gE\bar \mu_{i2}^tE^{\alpha}\eqno{(4.2.2)}$$

where $\alpha=0$ for $i=1$ and $\alpha=1$ for $i=2$ (this factor $E^{\alpha}$ is not
important). Analogs of (3.3.1.3) are

$$\im B_1=\left (\matrix *&A_{14}\\ A_{23}&A_{24} \endmatrix\right) ; \ \ \im B_2=\left
(\matrix *&A_{34}\\ A_{21}&A_{24} \endmatrix\right)\eqno{(4.2.3)}$$ These equalities show
that $B_1, B_2$ have the form $$\left (\matrix \n Z_p + \rd \n Z_p & \n Z_p + p\rd \n Z_p
\\\n Z_p + p\rd \n Z_p &\n Z_p + p^2\rd \n Z_p
\endmatrix\right)\eqno{(4.2.4)}$$ (sets mean that the corresponding entry of $B_1, B_2$
belongs to this set).

\medskip

{\bf (4.2.5)} Now we can formulate an analog of (3.3.1.4):

\medskip

$\s_1$ is $D$-equivalent to $\s_2 \iff \exists g \in G_V(\n Z_p)$ such that $B_1, B_2$
satisfy (4.2.4).

\medskip

In some cases I cannot find classes of D-equivalence. Nevertheless, as a first step of
the future investigations in many cases it is possible to prove that some $\s_1$,
$\s_2\in S(g,\goth i)$ satisfy the following condition that is weaker than (4.2.5):

\medskip

{\bf Condition 4.2.6.} $\exists g \in G_V(\n Z_p)$ such that

$$B_1, B_2\in \left (\matrix \n Z_p + \rd \n Z_p & \n Z_p + p\rd \n Z_p \\\n Z_p + p\rd
\n Z_p &\n Z_p + p\rd \n Z_p
\endmatrix\right)\eqno{(4.2.7)}$$ (the only difference is the first power of $p$ in
(2,2)-block) which is equivalent to the condition

$$\tilde B_1, \tilde B_2\in \left (\matrix \n F_p + \rd \n F_p & \n F_p\\\n F_p&\n
F_p\endmatrix\right)\eqno{(4.2.8)}$$

\medskip

Now we define the $\goth D$-partition of $S(g, \goth i)$ (see 1.6.2). We give formulas
only for $\s\in S(g,\goth i)^{open}$ leaving the definitions for other $\s$ as a subject
of further investigation, see Remark 4.2.11. Firstly,

$$\s\in \goth D_* \iff \det \mu_1(\tilde \s)= \det \mu_2(\tilde \s)=0\eqno{(4.2.9)}$$

Further, if $\det \mu_1(\tilde \s)\ne 0$ then for $j=0,\dots,\goth i$ we define

$$\s\in \goth D_j \iff \rank \ im \ \tilde X_1 = j\eqno{(4.2.10)}$$

{\bf Remark 4.2.11.} We have the following table of the corresponding formulas for $\goth
D$-partition for the cases $\goth T_p=T_p$ and $T_{p,\goth i}$:

\settabs 6 \columns

\medskip

\+ $\goth T_p\ \ =$ & $T_p$ && $T_{p,\goth i}$ \cr

\medskip

\+ & Lemma 3.2.3 (property) && Formulas 4.2.9, 4.2.10 (definition)\cr

\medskip

\+ & Formula 1.5.1 (definition) && Theorems 1.6.5, 1.6.6 (properties)\cr

\medskip

The reason of this non-coincidence is the following: theorems 1.6.5, 1.6.6 do not permit
us to distinguish between all possible cases, so I cannot take their formulas as
definitions. See also 1.8.6.

Particularly, I do not know the exact definition of $\goth D$-partition for the
complement of $S(g,\goth i)^{open}$, hence all subsequent theorems are proved only for
$S(g,\goth i)^{open}$.

\medskip

{\bf Remark 4.2.12.} If $\det \mu_1(\tilde \s)= 0$, $\det \mu_2(\tilde \s)\ne 0$ then
(see Remark 4.1.4) we can define $\goth D_j$, $j=0,\dots,\goth i$, interchanging the
roles of $\mu_1$, $\mu_2$ in (4.2.10). Since for any invertible complex matrix $X$ we
have $\rank \ im \ X=\rank \ im \ (X^{-1})$ the formula (4.2.10) is invariant
respectively the permutation of $\mu_1$, $\mu_2$.

\medskip

{\bf Remark 4.2.13.} For $g=3$, $\goth i=1$ this definition is equivalent to the
following one. We define a $2\times 2$-matrix $\goth R=\goth R(\s)$ with entries in $\n
F_p$ whose $i$-th line $(\goth r_{i1}, \goth r_{i2})$ is the line of coordinates of $\det
\mu_i(\tilde \s)$ in the basis $(1, \rd)$:

$$\goth r_{i1}+\goth r_{i2}\rd=\det \mu_i(\tilde \s),\ \ \ \goth r_{ij}\in \n F_p, \ \
i=1,2 $$

We have: ($k=0,1$):

$$\goth D_*=\{\s\in S(3,1)^{open}| \rank (\goth R(\s))=0\}\eqno{(4.2.13.1)}$$

$$\goth D_k=\{\s\in S(3,1)^{open}| \rank (\goth R(\s))=k+1\}\eqno{(4.2.13.2)}$$

{\bf Remark 4.2.14.} For $g=3$, $\goth i=2$ Lemma 5.3.5.4 shows that $\goth
D_0=\emptyset$.

\medskip

{\bf Theorem 4.2.15.} For any $g$, $\goth i$ we have: if $g-\goth i$ is odd then $\s_1$, $\s_2\in \goth D_{\goth
i}$ are D-equivalent. If $g-\goth i$ is even then $\s_1$ is D-equivalent to $\s_2$ iff the ratio $\det \tilde G_{221}
/ \det \tilde G_{222}$ is a square in $\n F_p^*$.


\medskip

{\bf Proof.} We have $\det \tilde \mu_{1i}\ne 0$, $i=1,2$. (4.2.2) implies that
$$gE=\mu_{11}^{-1}B_1\bar\mu_{12}^{t-1}\eqno{(4.2.15.1)}$$ $$B_2=\left (\matrix X_{11} &
X_{21}\\0&1 \endmatrix\right)B_1\left (\matrix \bar X_{12} & \bar X_{22}\\0&1
\endmatrix\right)^{t}E\eqno{(4.2.15.2)}$$ When we shall D-prove equivalence of $\s_1$,
$\s_2$ in this and subsequent theorems, we shall prove that there exist $B_1$, $B_2$
satisfying the following condition:

$$B_1, B_2\in \left (\matrix \n Z_p + \rd \n Z_p & \n Z_p\\\n Z_p&\n
Z_p\endmatrix\right)\eqno{(4.2.15.2')}$$ which is stronger than (4.2.4), and when we
shall prove non-D-equivalence, we shall prove that do not exist $B_1$, $B_2$ satisfying
(4.2.8). We denote $B_1=\left (\matrix \goth B_{11} & \goth B_{12}\\ \goth B_{21}&\goth
B_{22} \endmatrix\right)$. Let us assume that ($4.2.15.2'$) for $B_1$ is satisfied.
(4.2.15.2) gives us $$B_2=\left (\matrix * & X_{11}\goth B_{12} + X_{21}\goth B_{22}\\
\goth B_{21} \bar X_{12}^t + \goth B_{22} \bar X_{22}^t&\goth B_{22}
\endmatrix\right)E\eqno{(4.2.15.3)}$$

So, if

$$\im X_{11}\goth B_{12} + \im X_{21}\goth B_{22}=0 \eqno{(4.2.15.4)}$$

$$\im X_{12}\goth B_{21}^t + \im X_{22}\goth B_{22}^t=0 \eqno{(4.2.15.5)}$$

then ($4.2.15.2'$) for $B_2$ is satisfied.

Condition $\s_i\in \goth D_{\goth i}$ is equivalent to the condition $\det \im X_{1i} \in
\n Z_p^*$, $i=1,2$. So, if we set $$\goth B_{12} =W_{1}\goth B_{22}, \ \ \goth B_{21}
=\goth B_{22}W_2^t\eqno{(4.2.15.7)}$$

($W_i$ are defined in (4.1.8)) then (4.2.15.4), (4.2.15.5) will be satisfied.
Substituting the expression for $g$ from (4.2.15.1) to (3.1.1) we get that $B_1$ must
satisfy $$B_1F_2\bar B_1^t = \lambda G_1\eqno{(4.2.15.9)}$$ Writing down (4.2.15.9) in
block form and taking into consideration (4.2.15.7), we get that (3.1.1) is equivalent to
the following 3 equalities:

$$\goth B_{11} F_{112} \bar \goth B_{11}^t + W_1 \goth B_{22} F_{212} \bar \goth B_{11}^t
+ \goth B_{11}F_{122}\goth B_{22}^t W_1^t + W_1 \goth B_{22} F_{222} \goth B_{22}^t W_1^t
= \lambda G_{111}\eqno{(4.2.15.10)}$$

$$\goth B_{11}(F_{112}W_2 + F_{122}) \goth B_{22}^t + W_1 \goth B_{22}(F_{212}W_2 +
F_{222}) \goth B_{22}^t =\lambda G_{121}\eqno{(4.2.15.11)}$$

$$\goth B_{22}(W_2^t(F_{112}W_2 + F_{122}) +(F_{212}W_2 + F_{222}) )\goth
B_{22}^t=\lambda G_{221}\eqno{(4.2.15.12)}$$

For those who does not know Corollary 4.1.9 these formulas look rather complicated, does
not it? But using this corollary we get that (4.2.15.12) becomes

$$\goth B_{22}G_{222}^{-1}\goth B_{22}^t=\lambda G_{221}\eqno{(4.2.15.13)}$$

and (4.2.15.11) follows from (4.2.15.12). The theory of quadratic forms over $\n F_p$ shows that the statement of the theorem is necessary, and if it is satisfied, then
(4.2.15.11), (4.2.15.12) are satisfied.

Finally, (4.2.15.10) has the form $\goth B_{11} P \bar \goth B_{11}^t + Q^t \bar \goth
B_{11}^t + \goth B_{11} Q+R=0$ where $P,Q,R$ are some real matrices, $P,R$ are symmetric
and $\det \tilde P\ne 0$. The substitution $\goth B_{11}=\goth B_{11}'-Q^tP^{-1}$
eliminates the linear terms: we get $\goth B_{11}' P \overline {\goth B_{11}'}^t +R' =0$,
$R'=R-Q^tP^{-1}Q$. It is easy to prove that always $\det \tilde R'\ne 0$ (if $\s_2=\left
(\matrix E_g &0 \\ 0&pE_g\endmatrix\right)$ (i.e. if $A_2, U_2, C_2=0$) then this follows
immediately from (4.1.9.4); the proof for other $\s_2$ is similar). So, the theory of
Hermitian forms over $\n Z_p$ shows that we can choose $\goth B_{11}$ such that
(4.2.15.10) can be satisfied. $\square$

\medskip

{\bf Proposition 4.2.16.} For the case $g=3$ we have: if $\goth i=1$, $\s_1\in \goth
D_0$, $\s_2\in \goth D_1$ or $\goth i=2$, $\s_1\in \goth D_1$, $\s_2\in \goth D_2$ then
$\s_1$ is not D-equivalent to $\s_2$.

\medskip

{\bf Proof.} According Remark 4.1.4, we consider only the case when $\det \tilde
\mu_{11}\ne 0$. We have: $\det\im \tilde X_{11}=0$, $\det\im \tilde X_{12} \ne 0$. Let us
assume existence of $B_1$.

\medskip

Case $\goth i=1$. Equalities (4.2.15.1) - (4.2.15.3) and the reductions of (4.2.15.4) -
(4.2.15.5) continue to be true, so (4.2.15.7) becomes $$\tilde \goth B_{21} =\tilde \goth
B_{22}\tilde W_2^t\eqno{(4.2.16.1)}$$ hence $\im \tilde X_{21}\tilde \goth B_{22}=0$.
(4.2.16.1) implies that $\im \tilde X_{21}\tilde \goth B_{21}=0$, i.e. $\im \tilde X_{21}
\cdot (\tilde \goth B_{21} \ \ \tilde \goth B_{22})=0$. According Lemma 4.1.10, we see
that $\tilde B_1$ is degenerate. This contradicts to (4.2.15.1).

\medskip

Case $\goth i=2$. Equalities (4.2.15.1) - (4.2.15.3), the reductions of (4.2.15.4) -
(4.2.15.5), and (4.2.16.1) continue to be true. There exists a $1\times2$-matrix $D\ne0$
such that $D\ \im \tilde X_{11}=0$, hence the reduction of (4.2.15.4) implies that $D\
\im \tilde X_{21} \tilde \goth B_{22}=0$, Lemma 4.1.10 implies that $\tilde \goth
B_{22}=0$ and (4.2.16.1) implies that $\tilde \goth B_{21}=0$, hence $\det \tilde B_1=0$.
This contradicts to (4.2.15.1). $\square$

\medskip

{\bf Proposition 4.2.17.} For the case $g=3$, $\goth i=1$ we have: if $\s_1\in \goth D_0
\cup \goth D_1$, $\s_2\in \goth D_*$ then $\s_1$ and $\s_2$ are not D-equivalent.

\medskip

{\bf Proof.} $\s_1\in \goth D_0 \cup \goth D_1$ implies that there exists $i=1$ or 2 such
that $\det \mu_{i1}\ne 0\mod p$. For a matrix $A$ we denote by $A_{(j)}$ the $j$-th
column of $A$. According Lemma 4.1.12, $(\bar\mu_{i2}^t)_{(3)}=\gamma
(\bar\mu_{i2}^t)_{(2)}\mod p$ where the coefficient of proportionaly $\gamma\in \n
F_p(\rd)$ obviously is not real. This implies that $(B_i)_{(3)}=\gamma (B_i)_{(2)}\mod
p$, hence both $(B_i)_{(2)}= (B_i)_{(3)}=0\mod p$. Since $(B_i)_{(2)}=\mu_{i1}\goth
gE\cdot(\bar\mu_{i2}^t)_{(2)}$ and $(\bar\mu_{i2}^t)_{(2)}\ne0\mod p$, we get that the
columns of $\mu_{i1}\goth gE$ are linearly dependent modulo $p$ --- a contradiction to
the condition $\det \mu_{i1}, \det \goth g\ne 0\mod p$. $\square$

\medskip

{\bf Proposition 4.2.18.} Let $\s\in\goth D_*$. Then $T_{p,1}(V)_{I(\s)}=V \iff
\s\in\goth D_{0,*}$ (see 1.4.1 for the notations). Particularly, $\goth D_{0,*}$ is a
part of D-partition.

\medskip

{\bf Proof.} To simplify calculations, we consider only the case of $\s$ such that (in
notations of 4.1.1) $A,C=0$, and $U\in M_2(\n Z)$ is any symmetric matrix. Calculations
for the remaining $\s$ are similar but more complicated. Like in the proof of the theorem
3.3.5, we use 3.3.5.1.

This means that $T_{p,1}(V)_{I(\s)}=V$ iff there exists $\goth g\in G_V$ such that
$\ord_p(\det \goth g)=3$ and entries of $\alpha(\goth g)\s^{-1}$ are integer.

We write $\goth g=h+k\rd$ and we denote block entries of $\goth g,h,k$ by $\goth g_{ij}$,
$h_{ij}$, $k_{ij}$ respectively. Writing down the entries of $\alpha(\goth g)\s^{-1}$ we
get that they are integer iff

\medskip

{\bf (4.2.18.0)} Entries of $h_{11}$, $h_{21}$, $k_{11}$, $k_{21}$ are $\equiv 0 \mod p$

and

$$\matrix -h_{12}U+k_{12}\equiv 0 \mod p^2 \\ \d k_{12}U+h_{12}\equiv 0 \mod p^2 
\endmatrix\eqno{(4.2.18.1)}$$

$$\matrix -h_{22}U+k_{22}\equiv 0 \mod p^2 \\ \d k_{22}U+h_{22}\equiv 0 \mod p^2 
\endmatrix\eqno{(4.2.18.2)}$$

(4.2.18.1) (resp. (4.2.18.2)) implies that $(1+\d U^2)h_{12}$ (resp. $(1+\d U^2)h_{22}$)
$\equiv 0 \mod p^2$.

Further, we have:$$\ord_p(\det(1+U\rd))=0 \iff \ord_p(\det(1+\d U^2))=0;$$

$$\ord_p(\det(1+U\rd))=1 \iff 1+\d U^2\equiv 0 \mod p, \ \ \ \not\equiv 0 \mod p^2,$$ and
in this case $\ord_p(\det\frac{1+\d U^2}p)=0$;

$$\ord_p(\det(1+U\rd))\ge 2 \iff 1+\d U^2\equiv 0 \mod p^2$$

This means that if $\ord_p(\det(1+U\rd))=0$ then $h_{i2}$, $k_{i2} \equiv 0 \mod p^2$;
together with (4.2.18.0) this means that $\ord_p(\det \goth g)\ge 5$, i.e.
$T_{p,1}(V)_{I(\s)}\ne V$.

If $\ord_p(\det(1+U\rd))=1$ then (4.2.18.0) --- (4.2.18.2) imply that $h_{ij}=ph'_{ij}$,
$k_{ij}=pk'_{ij}$ with integer $h'_{ij}$, $k'_{ij}$ and that $h'_{i2}$, $k'_{i2}$ satisfy

$$\matrix -h'_{12}U+k'_{12}\equiv 0 \mod p \\ \d k'_{12}U+h'_{12}\equiv 0 \mod p 
\endmatrix$$

$$\matrix -h'_{22}U+k'_{22}\equiv 0 \mod p \\ \d k'_{22}U+h'_{22}\equiv 0 \mod p 
\endmatrix$$

These conditions mean that $\left (\matrix h'_{12}+\rd k'_{12}\\ h'_{22}+\rd
k'_{22}\endmatrix\right) (1-\rd U)=0$, i.e. the second and the third columns of $\frac1p
\goth g$ are linearly dependent. This means that $\ord_p(\det \goth g)>3$, i.e.
$T_{p,1}(V)_{I(\s)}\ne V$.

Finally, if $\ord_p(\det(1+U\rd))\ge 2$ then we can take $\goth g_{11}=p$, $\goth
g_{12}=\goth g_{21}=0$, $\goth g_{22}=h_{22}(1+U\rd)$, and the condition $\lambda(\goth
g)=p^2$ becomes $\frac{1+\d U^2}{p^2}=h_{22}^{-1}{h_{22}^{-1}}^t$. This equation can be
solved for any $U$ such that the corresponding $\s$ is in $\goth D_{0,*}$. $\square$

\medskip

{\bf Theorem 4.2.19.} For $g=3$, $\goth i=1$ the field of definition of any intermediate
irreducible component of $T_{p,1}(V)$ is $K^p$.

\medskip

{\bf Proof} is similar to the proof of the theorem 3.3.6. In the present case $\s\in
S(g,\goth i)^{open}$ is given by (4.1.1). It satisfies the analog of (2.13), and we can
take $r=\tau_{p,\goth i}\s^{t-1}$. The $p$-component of $\Cal K_p$ is the set of matrices
$\gamma\in GSp_{2g}(\n Z_p)$ having the block structure $\gamma_{ij}$, $i,j=1,\dots,4$
such that

$$\gamma_{12},\gamma_{32},\gamma_{41},\gamma_{43}\equiv 0 \mod p,\ \ \gamma_{42}\equiv 0
\mod p^2\eqno{(4.2.19.1)}$$

Analogs of (3.3.6.1) - (3.3.6.3) hold for the present situation, and (3.3.6.4) is changed
(according (4.5.2.1)) as follows:

\medskip

{\bf (4.2.19.2)} $\goth g$ belongs to the $p$-component of $\Cal K_{V,p,r}$ iff

\medskip

$[r \alpha(\goth g) r^{-1}]_{12}$, $[r \alpha(\goth g) r^{-1}]_{32}$, $[r \alpha(\goth g)
r^{-1}]_{41}$, $[r \alpha(\goth g) r^{-1}]_{43}$ are $\equiv 0 \mod p$,

\medskip

and $[r \alpha(\goth g) r^{-1}]_{42}$ is $\equiv 0 \mod p^2$

\medskip

To formulate an analog of (3.3.6.5), we denote $B_i=B_i(\s,\s)$ ($i=1,2$) like in
(4.2.2). The analog of (3.3.6.5) is the following:

$$\im B_1=\left( \matrix *&[r \alpha(\goth g) r^{-1}]_{41}\\[r \alpha(\goth g)
r^{-1}]_{32}&[r \alpha(\goth g) r^{-1}]_{42}\endmatrix \right)\ \ \ \im B_2=\left(
\matrix *&[r \alpha(\goth g) r^{-1}]_{43}\\[r \alpha(\goth g) r^{-1}]_{12}&[r
\alpha(\goth g) r^{-1}]_{42}\endmatrix \right)\eqno{(4.2.19.3)} $$

According Remark 4.1.4, we can restrict ourselves by the case when $\det \mu_1(\tilde
\s)\ne 0$. Let us take any $\goth g\in \Cal K_{V,p,r}$. Since $\tilde X_1$ is real, the
reduction of (4.2.15.4) implies that $(\im \tilde X_2) \tilde \goth B_{22}=0$, and lemma
4.1.10 implies that $\det \tilde \goth B_{22}=0$. (4.2.19.2), (4.2.19.3) imply that
$\tilde B_1=\left( \matrix *+\rd * & *\\ *& \tilde \goth B_{22}\endmatrix \right)$

(entries of *'s are in $\n F_p$). Since sizes of blocks are (1,2), this means that $\im
\det \tilde B_1=0$, i.e. $\det B_1\in \n Z_p+p\rd\n Z_p$, and hence $\det \goth g\in \n
Z_p+p\rd\n Z_p$. For any $\s\in\goth D_0$, any $\alpha\in\n Z_p+p\rd\n Z_p$ it is easy to
construct examples of matrices $\goth g\in \Cal K_{V,p,r}$ such that $\det \goth g
=\alpha$. This means that the field of definition of this component is $K^p$. $\square$

\medskip

To formulate the conjecture on the restriction of I-partition on $\goth D_0$ and on the
Galois action on the corresponding set of irreducible components, we must recall that the
projectivization of the set of $2\times 2$-matrices of rank 1 is $P^1 \times P^1$. Recall
that for $\s\in \goth D_0$ the matrix $\goth R(\s)$ (defined in 4.2.13) is of rank 1. Let
$P(\goth R(\s))\in P^1(\n F_p)\times P^1(\n F_p)$ be its projectivization. So, there
exists a map $\goth D_0\to P^1(\n F_p)\times P^1(\n F_p)$ defined on $\s\in \goth D_0$ as
follows: $\s \mapsto P(\goth R(\s))$. Let $\pi_i:\goth D_0\to P^1(\n F_p)$ ($i=1,2$) be
the composition of this map with the projection of $P^1(\n F_p)\times P^1(\n F_p)$ to the
$i$-th factor. Let $\s\in \goth D_0$. We denote the restriction of $\pi_i$ to
$S(3,1)_{D(\s)}$ by $\pi_{\s,i}$ (we use notations of (1.4.1), (1.4.2)).

\medskip

{\bf Conjecture 4.2.20.} There exists $i=1$ or 2 such that parts of I-partition
restricted to $S(3,1)_{D(\s)}$ coincide with the fibers of $\pi_{\s,i}$. Particularly,
$D(\s)_I$ --- the set of these parts --- is indexed by elements of $P^1(\n F_p)$, and the
action of Galois group $\Gal(K^p/K^1)=\n F_{p^2}^*/\n F_p^*$ on $D(\s)_I$ coincides with
the natural action of $\Gal(K^p/K^1)$ on $P^1(\n F_p)$.

\medskip

I do not know whether $i$ is 1 or 2.

\medskip

{\bf Remark 4.2.21.} For $g=3$, $\goth i=1$ there exists a constant $c$ (which does not
depend on $p$) and an extension $L/K^1$ of degree $c$ such that the field of definition
of any special component of $T_{p,1}(V)$ is a subfield of $L$.

\medskip

gde-to eshchyo, no ne znayu, gde].

\medskip

{\bf Idea of the proof.} We use notations of Theorem 4.2.19. We consider only one $\s$
given by (4.1.1) with $A,C=0$, $U=\left( \matrix u_2&u_3\\u_3&-u_2\endmatrix \right)$
where $u_2$, $u_3$ satisfy $\Delta(u_2^2+u_3^2)\equiv -1\mod p$. We denote by $\nu$ the
coefficient of proportionality between the second and the third lines of $\tilde
\mu_1(\s)$. Obviously $\nu$ is not real, hence (4.2.2), (4.2.19.2), (4.2.19.3) imply that
$$\tilde B_1=\left( \matrix *&0&0\\0&0&0 \\0&0&0\endmatrix \right)\eqno{(4.2.21.1)}$$
where $*\in \n F_{p^2}^*$. Moreover, the expression for $\mu_2$ shows that (4.2.21.1) is
equivalent to (4.2.8) for $B_2$.

A direct calculation shows that (4.2.21.1) is equivalent to the following relations
between entries $\goth g_{ij}$ of $\goth g$:

$$\goth g_{22} - \goth g_{33} - \nu \goth g_{23} - \bar\nu \goth
g_{32}=0\eqno{(4.2.21.2)}$$

$$-\nu \goth g_{12}+\goth g_{13}=0\eqno{(4.2.21.3)}$$

$$\nu \goth g_{21}+\goth g_{31}=0\eqno{(4.2.21.4)}$$


Even diagonal matrices $\goth g$ satisfying (4.2.21.2)--(4.2.21.4) can have (almost) any
determinant in $\n F_{p^2}^*$ after reduction; really, their determinants form a subgroup
of index $\le 3$ in $\n F_{p^2}^*$. This implies the desired. $\square$

\medskip

{\bf Section 5. Miscellaneous.} \nopagebreak

\medskip

{\bf Subsection 5.1. Case $g=3$, $\goth i=2$.} \nopagebreak

\medskip

We use notations of Section 4.1 (particularly of 4.1.1), so $E$ will mean $E_{12}$. Let
us describe $\goth D_*$. We denote $A=(a_1 \ a_2)$, $C=(c_1 \ c_2)$ (the $1\times
2$-matrices), $U=(u)$ (the $1\times 1$-matrix). Further, if $p\equiv 3 \mod 4$ then we
denote $\frac1{\sqrt{\Delta}}$ by $r$ (we have $r\in \n F_p$) and $r \rd$ by $i$ (we have
$i\in \n F_{p^2}$ is purely imaginary and $i^2=-1$).

\medskip

{\bf Theorem 5.1.1.} $\goth D_*$ is non-empty iff $p\equiv 3 \mod 4$. In this case $\goth
D_*\cap S(3,2)^{open}$ is the set of matrices of form 4.1.1 such that

$$a_1^2+a_2^2\equiv 1, \ \ U\equiv 0, \ \ C\equiv (\ve r a_2, -\ve r a_1)\mod
p\eqno{(5.1.1.1)}$$ where $\ve=\pm 1$.

\medskip

{\bf Proof. } Direct calculations. Conditions $\det \mu_1(\tilde \s)=\det \mu_2(\tilde
\s)=0$ are the following (equalities in $\n F_p$):

$$u+a_1c_1+a_2c_2=0$$

$$u-a_1c_1-a_2c_2=0$$

$$a_1^2+a_2^2=1$$

$$c_1^2+c_2^2=1/\Delta$$

Solving this system we get the desired. $\square$

\medskip

{\bf Remark.} I am very astonished by this result, because this is the only case (in
[L01] and in the present paper) when something depends on residue of $p$ modulo 4. 

\medskip

{\bf Theorem 5.1.2.} There exist $\s_1\in \goth D_1, \s_2 \in \goth D_*$ which are not
D-equivalent.

\medskip

{\bf Remark.} Conjecturally, any $\s_1\in \goth D_1, \s_2 \in \goth D_*$ are not
D-equivalent.

\medskip

{\bf Proof.} We choose $\s_1$, $\s_2$ having respectively (in notations of 4.1.1)

$$A=(1\ \  0), \ \ U=0, \ \ C=(0\ \ 0)\eqno{(5.1.2.1)}$$

$$A=(1\ \  0), \ \ U=0, \ \ C=(0\ \ r)\eqno{(5.1.2.2)}$$

(i.e. for $\s_2$ we have $a_1=1$, $a_2=0$, $\ve=-1$ of 5.1.1.1). We shall prove that the
condition 4.2.6 is not satisfied for these $\s_i$. gg We shall make all calculations in
$\n F_p$ without indicating tilde, i.e. all elements are reduced. We have

$\mu_{12}= \left( \matrix -1 & 0 & 0
\\ 0 & -1& -i
\\ 1 & i &-1
\endmatrix \right)$,

$\mu_{22}= \left( \matrix \rd & 0 & -\rd
\\ 0 & \rd& 0
\\ 1 & i &-1
\endmatrix \right)$.

We denote the third line of $\mu_{11}$ by $l$. Conditions (4.2.8) for $B_1$ have the
following form:

\medskip

5.1.2.3. $\mu_{11}\goth gE \bar l^t$ is real (condition of reality of $(B_1)_{*3}$);

\medskip

5.1.2.4. $(\mu_{11}\goth gE)_{31}$ is real (condition of reality of $(B_1)_{31}$).

\medskip

Reality of $(B_1)_{32}$ follows from these conditions, because the second line of
$\mu_{12}$ is a lineal combination with real coefficients of its first and third lines.

\medskip

Conditions (4.2.8) for $B_2$ have the following form:

\medskip

5.1.2.5. $\mu_{21}\goth gE \bar l^t$ is real (condition of reality of $(B_2)_{*3}$);

\medskip

5.1.2.6. $(\mu_{21}\goth gE)_{32}$ is imaginary (condition of reality of $(B_2)_{32}$);

\medskip

5.1.2.7. $-r (B_2)_{32} + (B_2)_{33}=0$ (because the first line of $\mu_{22}$ is the
linear combination of its second and third lines with coefficients $-i$, $\rd$).

\medskip

Since $$\mu_{21}=\left( \matrix X_{11}&X_{21} \\ 0 & 1 \endmatrix
\right)\mu_{11}\eqno{(5.1.2.8)}$$

we see that (5.1.2.3), (5.1.2.5) are equivalent to the condition

\medskip

5.1.2.9. $\mu_{11}\goth gE \bar l^t=\left( \matrix *\\ 0 \\0 \endmatrix \right)$ where
$*$ is real. So,  (5.1.2.7) implies

$(B_2)_{32}=0 \Longrightarrow (\mu_{21}\goth gE)_{32}=0 \Longrightarrow (\mu_{11}\goth
gE)_{32}=0$ (because of 5.1.2.8).

\medskip

We see that 5.1.2.7, 5.1.2.9 are equivalent to the condition

$$\mu_{11}\goth g=\left( \matrix v_{11}&v_{12}&v_{11} -i v_{12} +z_1 \\
v_{21}&v_{22}&v_{21} -i v_{22}\\
z_2& 0 & z_2 \endmatrix \right)\eqno{(5.1.2.10)}$$ where $v_{ij}\in \n F_{p^2}$,
$z_{i}\in \n F_{p}$.

\medskip

Now we take $Z=\left( \matrix &0&1\\
\bar l^t&1&0\\
&0&0\endmatrix \right)$, so $\goth g'\overset{\opr}\to{=}\mu_{11}\goth gZ=
\left( \matrix z_1&v_{12}&v_{11}\\
0&v_{22}&v_{21}\\
0& 0 & z_2 \endmatrix \right)$. In terms of ${\goth g'}^{-1}$ the condition (3.1.1) has
the form

$$Z^{-1}E\bar Z^{t-1}=\lambda {\goth g'}^{-1}G_1\bar {\goth g'}^{t-1}\eqno{(5.1.2.11)}$$

($G_1=\mu_{11}E\bar \mu_{11}^t$ from 4.1). We have: $Z^{-1}E\bar Z^{t-1}=\left( \matrix
-1&-i&-1\\ i&0&i\\ -1&-i&0\endmatrix \right)$ and ${\goth g'}^{-1}$ has entries similar
to the entries of $\goth g'$. Writing explicitly the (2,2)-entry of the right hand side
of (5.1.2.11), we get a contradiction. $\square$

\medskip

{\bf Proposition 5.1.3.} The field of definition of the intermediate component that
corresponds to $\s$ given by (5.1.2.1) is $K$.

\medskip

{\bf Proof.} We use notations of Theorem 4.2.19. We denote the $(i,j)$-th element (not
block!) of $B_1$ by $b_{ij}$, $1\le i,j\le 3$. We get that the left hand sides of
(4.2.15.4), (4.2.15.5) are $\left( \matrix b_{33}\\-b_{23}\endmatrix \right)$, $\left(
\matrix \bar b_{33}\\ - \bar b_{32}\endmatrix \right)$, hence the condition that $B_2$
satisfies (4.2.4) is $\tilde B_1=\left( \matrix *&*&b_{13}\\*&*&0\\b_{31}&0&0\endmatrix
\right)$ where $b_{13}$, $b_{31}\in \n F_p$ and $*\in \n F_{p^2}$. Writing down
explicitly condition (3.1.1) 

we get immediately that $\det B_1$ --- and hence $\det \goth g$ --- can take any value.
So, the field of definition of the corresponding component is $K$. $\square$

\medskip

{\bf Remark 5.1.4.} Conjecturally, all special components of $T_{p,2}(V)$ are defined
over some field $L$ of Remark 4.2.21.

\medskip

{\bf Subsection 5.2.} Action of $T_p$ on components of $T_p(V)$.  \nopagebreak

\medskip

We use notations of (1.4.1). Let $\s\in S(g)$ be an element. For $\goth i=0,\dots,g$ we
denote

$$S(g,\goth i,\s)=\{\s'\in S(g)|\s'\s\in \Gamma\tau_{p,\goth i}\Gamma\}$$ Therefore
$S(g)=\bigcup_{\goth i=0}^g S(g,\goth i,\s)$ (all unions in this subsection are
disjoint), and for $k \in \Cal S_D(g,\goth i)$ we denote

$$S(g,\goth i,\s)_k=\{\s'\in S(g,\goth i,\s)|\s'\s\in \bigcup_{\s_j\in S(g,\goth i)_k}
\Gamma\s_j\}.$$

It is easy to check that if $\s_1, \s_2 \in S(g)$ are D-equivalent then generally
$S(g,\goth i,\s_1)_k\ne S(g,\goth i,\s_2)_k$, but $$S(g,\goth i,\s_1)_k\ne\emptyset \iff
S(g,\goth i,\s_2)_k\ne\emptyset\eqno{(5.2.1)}$$

Let us fix some $\goth k\in \Cal S_D(g)$ (i.e. a class of D-equivalence in $S(g)$ = a
Shimura subvariety $T_p(V)_{\goth k}$ of $T_p(V)$).The answer on the problem of
description of $T_p(T_p(V)_{\goth k})$ is given by the following two propositions:

\medskip

{\bf Proposition 5.2.2.} Let $\goth i\in\{0,\dots,g\}$, $k \in \Cal S_D(g,\goth i)$. Then
$T_{p,\goth i}(V)_k$ is a component of $T_p(T_p(V)_{\goth k})$ iff for $\s\in S(g)_{\goth
k}$ we have: $S(g,\goth i,\s)_k\ne\emptyset$ (according (5.2.1), this condition does not
depend on a choice of $\s$). $\square$

\medskip

Now let us find miltiplicities. We have a formula $T_p^2=\sum_{\goth i=0}^gW_p(\goth
i)T_{p,\goth i}$ where $W_p(\goth i)=\prod_{j=1}^\goth i(p^j+1)$. Let us fix $\s''\in
S(g,\goth i)_k$. We can represent the coset $\Gamma\s''$ by $W_p(\goth i)$ ways as the
following product:

$$\Gamma\s'' = \Gamma\s'_j \s_j\eqno{(5.2.3)}$$

where $\s'_j$, $\s_j\in S(g)$, and $j=1, \dots, W_p(\goth i)$ is the number of the way of
representation as the product.

\medskip

{\bf Proposition 5.2.4.} The multiplicity of $T_{p,\goth i}(V)_k$ in $T_p(T_p(V)_{\goth
k})$ is the quantity of $j$ in (5.2.3) such that we have $\s'_j\in S(g)_{\goth k}$.
$\square$

\medskip

Now let us apply propositions 5.2.2, 5.2.4 to some $\goth i, \goth k, k$. We denote 2
good elements of $\Cal S_D(3)$ by $k_{g1}$, $k_{g2}$ and the bad element by $k_b$.
Analogously, we denote 2 general elements of $\Cal S_D(3,1)$ by $k_{1g1}$, $k_{1g2}$, the
intermediate elements by $k_{1i\alpha}$ ($\alpha$ runs over the set of components of
$T_{p,1}(V)$ of intermediate type), the special element corresponding to $\goth D_{0,*}$
by $k_{1s0}$ and other special elements by $k_{1s\alpha}$ ($\alpha$ runs over the set of
components of $T_{p,1}(V)$ of special type).

\medskip

{\bf 5.2.5. Case $\goth i=1$, $\goth k=k_{g1}$.} We take $\s=\left( \matrix E_3 &
0\\0&pE_3\endmatrix \right)\in S(3)_{k_{g1}}$, and let $\s'\in S(3)$ be a generic element
such that $\s'\s\in S(3,1)$. This $\s'$ corresponds to the case $I=\{2,3\}$, where $I$ is
from [L04.1], Section 2.3. We take $\s'$ from the proof of the Proposition 3.4.3, Case
$I=\{2,3\}$, $D=(d_1 \ \ d_2)$. In these notations $\s'\s$ is a matrix described in
(4.1.1), with $A=D^t$, $C=0$, and $U$ of (4.1.1) is $\equiv 0 \mod p$. This means that
$\mu_1(\s'\s)= \left( \matrix -E_1&0\\D^t&-E_2 \endmatrix \right)$, $\mu_2(\s'\s)= \left(
\matrix \rd E_1&-D\rd \\D^t&-E_2 \endmatrix \right)$. So, $\goth R(\s'\s) = \left(
\matrix -1 & 0 \\ 0 & 1-d_1^2 - d_2^2 \endmatrix \right)$ ($\goth R$ is defined in
4.2.13). This means that $\s'\s$ never belongs to $\goth D_*$, and belongs to $\goth D_0$
iff $d_1^2 + d_2^2 =1$. So, we have got a

\medskip

{\bf Proposition 5.2.6.} $T_p(T_p(V)_{k_{g1}})$ (i.e. $T_p$ of a good component of
$T_p(V)$) contains components of $T_{p,1}(V)$ of general and intermediate types. $\square$

\medskip

{\bf Remark 5.2.7.} (a) Consideration of other (non-generic) components of $S(3)$, and of
$k_{g2}$, shows that $T_p$ of any good component of $T_p(V)$ does not contain components
of $T_{p,1}(V)$ of special type.

\medskip

(b) In terms of the Pl\"ucker coordinates on $S(3)$ (see notations of Remark 3.2.5) the
condition $d_1^2 + d_2^2 =1$ is equivalent to the condition $v_{1'23}(\s') -
v_{12'3}(\s') - v_{123'}(\s')=0$.

\medskip

{\bf 5.2.8. Case $\goth i=1$, $\goth k=k_{b}$.} We take $\s$ from the proof of the
Proposition 3.4.3 (see the first line of the proof), and $\s'$ as in 5.2.5. We have
$$\goth R(\s'\s) = \left( \matrix - \Delta w_1 (-d_1^2 w_1 + 2 d_1 w_2 -w_1) &
-w_1(1-d_1^2) \\ \Delta  (-d_1^2 w_1 + 2 d_1 w_2 -w_1) & 1-d_1^2 \endmatrix
\right)\eqno{(5.2.9)}$$ We have $\det\goth R(\s'\s) =0$  
computer, file logachev19, calculations on pages 21.1.2, 45.1-3)]. So, we have got a

\medskip

{\bf Proposition 5.2.10.} $T_p$ of the bad component of $T_p(V)$ contains the component
of $T_{p,1}(V)$ of intermediate type. $\square$

\medskip

As earlier we have a

\medskip

{\bf Remark 5.2.11.} Consideration of other (non-generic) components of $S(3)$ shows that
$T_p$ of the bad component of $T_p(V)$ does not contain components of $T_{p,1}(V)$ of
general type, and does contain components of special type (although the above matrices
(5.2.9) are all non-0).

\medskip

Now we use Proposition 5.2.4 in order to find multiplicity for the case when $k=k_{1i1}$
corresponds to a component of intermediate type. We take the simplest value of $\s''$
defined by (4.1.1) with $U=0$, $C=0$, $A^t=(1 \ \ 0)$, this $\s''$ belongs to
$S(3,1)_{k_{1i1}}$. To describe $\s_j$, $\s'_j$ ($j=0, \dots, p$; $W_p(1)=p+1$) we need
the following matrices:

$M=\left( \matrix p & 0 & 0 & 0& 0 & 0
                \\ -1 & 1 &  0 & 0& 0&0
                \\ 0 & 0 & 1 & 0 & 0&0
                                \\ 0 & 0 & 0 &1&1&0
                                \\ 0 & 0 & 0 &0&p&0
                \\ 0 & 0 & 0 & 0 & 0&p \endmatrix \right)$,

$M_j= \left( \matrix 1 & 0 & 0 & j&j&0
                \\ 0 & 1 &  0 &j&j&0
                \\ 0 & 0 & 1 & 0&0&0
                                \\ 0 & 0 & 0 &p&0&0
                                \\ 0 & 0 & 0 &0&p&0
                \\ 0 & 0 & 0 & 0 & 0&p\endmatrix \right)$, $j=0, \dots, p-1$.

We have $\s_j=M$, $\s'_j=M_j$, $j=0, \dots, p-1$ and $\s_p=M_0$, $\s'_p=M$. Using for
example (3.2.6), (3.2.7) we get that $M\in S(3)_{k_b}$ and all $M_j\in S(3)_{k_{g1}}\cup
S(3)_{k_{g2}}$. So, we get

\medskip

{\bf Proposition 5.2.12. (a)} The multiplicity of the simplest intermediate component of
$T_{p,1}(V)$ in $T_p(T_p(V)_{k_b})$ is $1$.

\medskip

{\bf (b)} Let $\goth m_i$ ($i=1,2$) be the multiplicity of the simplest intermediate
component of $T_{p,1}(V)$ in $T_p(T_p(V)_{k_{gi}})$. Then $\goth m_i$ is the quantity of
$j\in\{0, \dots, p-1\}$ such that $M_j\in S(3)_{k_{gi}}$ (i.e. the signature of $\tilde
T(M_j)$ is of type $i$). Particularly, $\goth m_1+\goth m_2=p$. $\square$

\medskip

{\bf Remark 5.2.12$'$.} We can expect that are similar results hold for any intermediate
component of $T_{p,1}(V)$.

\medskip

{\bf 5.2.13. Case $\goth i=2$, $\goth k=k_{g1}$.} We follow 5.2.5. Namely, we take the
same $\s$ as in 5.2.5, and let $\s'\in S(3)$ be a generic element such that $\s'\s\in
S(3,2)$. This $\s'$ corresponds to the case $I=\{3\}$, where $I$ is from [L04.1], Section
2.3. We take $\s'$ from the proof of the Proposition 3.4.3, Case $I=\{3\}$, $D^t=(d_{13}
\ \ d_{23})$. In these notations $\s'\s$ is a matrix described in (4.1.1), with $A=D^t$,
$C=0$, and $U$ of (4.1.1) is $\equiv 0 \mod p$. We have $\mu_1(\s' s)= \left( \matrix
-E_2&0\\A&-E_1 \endmatrix \right)$, $\mu_2(\s' s)= \left( \matrix \rd E_2&-A^t\rd
\\A&-E_1 \endmatrix \right)$. We get that always $\det \mu_1\ne 0$, $X_1=\rd \left(
\matrix 1+d_{13}^2&d_{13} d_{23}\\d_{13} d_{23}&1+d_{23}^2\endmatrix \right)$, $\det \im
X_1= 1+d_{13}^2+d_{23}^2$. So, we have got a

\medskip

{\bf Proposition 5.2.14.} $T_p(T_p(V)_{k_{g1}})$ (i.e. $T_p$ of a good component of
$T_p(V)$) contains components of $T_{p,2}(V)$ of general and intermediate types. $\square$

\medskip

{\bf Remark 5.2.15.} Apparently, $T_p$ of any good component of $T_p(V)$ does not contain
components of $T_{p,2}(V)$ of special type. To check it carefully, we must consider $\s'$
of non-generic type.

\medskip

{\bf 5.2.16. Case $\goth i=2$, $\goth k=k_{b}$.} It is easy to check that $T_p$ of the
bad component of $T_p(V)$ contains all components of $T_{p,2}(V)$.

\medskip

{\bf Subsection 5.3.} A theorem on a weak ``equivalence'' of components of $\goth
D_{1,*}$. \nopagebreak

\medskip

Here we prove that in many cases the condition that $\s_1$, $\s_2$ are $\goth
D$-equivalent implies that $\s_1$, $\s_2$ satisfy 4.2.6.

Throughout all this subsection all calculations will be made in $\n F_p$, i.e. we shall
consider reductions of all objects. To simplify notations we shall not indicate tilde.

We use here notations of the beginning of Subsection 4.2.

\medskip

{\bf Proposition 5.3.1.} For the case $g=3$, $\goth i=1$ we have: any $\s_1$, $\s_2\in
\goth D_0\cap S(3,1)^{open}$ satisfy 4.2.6.

\medskip

{\bf Proof.} $\s\in \goth D_0 \iff \{\det \mu_1(\s)=0$ and $\det \mu_2(\s)\ne 0\}$ or
$\{\det \mu_1(\s)\ne 0$ and $\det\im X_{1}=0\}$.

\medskip

{\bf Case 1.} Both $\s_1$, $\s_2$ satisfy $\{\det \mu_1\ne 0$ and $\im X_{1}=0\}$.

\medskip

According (4.1.6), $F_{11i}=0$. (4.2.15.4), (4.2.15.5) and (4.1.6) give us:

$$F_{121} \goth B_{22}=F_{122} \goth B_{22}^t=0\eqno{(5.3.1.1)}$$

This means that if

$$\goth B_{22} = k F_{211}^O F_{122}^O\eqno{(5.3.1.2)}$$

where $k\in \n F_p^*$ (recall that all $F$, $G$ are real symmetric) then $B_2$ satisfy
(4.2.8).

Analogs of (4.2.15.10) - (4.2.15.12) for the present case (taking into consideration
(5.3.1.1)) are

$$\goth B_{12} F_{212} \bar \goth B_{11}^t + \goth B_{11}F_{122}\goth B_{12}^t + \goth
B_{12} F_{222} \goth B_{12}^t = \lambda G_{111}\eqno{(5.3.1.3)}$$

$$\goth B_{12}(F_{212}\goth B_{21}^t + F_{222} \goth B_{22}^t) =\lambda
G_{121}\eqno{(5.3.1.4)}$$

$$\goth B_{22} F_{222} \goth B_{22}^t=\lambda G_{221}\eqno{(5.3.1.5)}$$

Substituting (5.3.1.2) to (5.3.1.5) we get 

$$k^2 F_{211}^O (F_{122}^O F_{222} F_{122}^{Ot}) F_{211}^{Ot}= \lambda G_{221}$$ Applying
(4.1.11.1) for $F_2$ and (4.1.11.2) for $F_1$ we get

$$k^2\frac{\det F_2}{\det F_1}G_{221}=\lambda G_{221}$$

i.e. if we take any $k\in \n Z_p^*$, $\lambda=k^2\frac{\det F_2}{\det F_1}$ and $\goth
B_{22}$ from (5.3.1.2) then (5.3.1.5) will be satisfied.

Substituting (5.3.1.2) to (5.3.1.4) we get (taking into consideration (4.1.11.3) for
$F_2$)

$$\goth B_{12}Z =\lambda G_{121}\eqno{(5.3.1.6)}$$ where $Z=F_{212}\goth B_{21}^t -k
(\det F_2)G_{212}^O F_{121}^O$.

Since vectors $F_{212}$, $G_{212}^O$ are never linearly dependent (because of (4.1.11.4))
and $F_{121}\ne 0$, we get that for almost all $\goth B_{21}$ $\det Z\ne 0$. We choose
any such $\goth B_{21}$, so for $\goth B_{12}= \lambda G_{121} Z^{-1}$ (5.3.1.4) is
satisfied.

Finally, in order to show that (5.3.1.3) has a solution respectively $\goth B_{11}$ it is
sufficient to check that $\goth B_{12} F_{212}\ne 0$. Let us assume the contrary: $\goth
B_{12} F_{212}=0$. This means that $\goth B_{12} = \alpha F_{122}^O$. Substituting this
value to (5.3.1.6) we get $ \alpha k (\det F_2) F_{122}^O G_{212}^O F_{121}^O = \lambda
G_{121}$. This implies that vectors $F_{121}^O$, $G_{121}$ are linearly dependent. But
their determinant is 1 (because of (4.1.11.4)) --- a contradiction.

\medskip

{\bf Case 2.} $\s_1$ satisfy $\det \mu_{1}(\s_1)\ne 0$, $\det \mu_{2}(\s_1)=0$, and
$\s_2$ satisfy $\det \mu_{1}(\s_2)=0$, $\det \mu_{2}(\s_2)\ne 0$.

\medskip

Here we change slightly the definitions of the previous cases. Namely, we set:
$\mu_{21}=\left (\matrix 0&X_{21} \\ 0&1\endmatrix\right)\mu_{11}$, $\mu_{12}=\left
(\matrix 0&Y_{22} \\ 0&1\endmatrix\right)\mu_{22}$,  $B_0=\mu_{11}gE\bar \mu_{22}^t$,
$B_0=\left (\matrix \goth B_{11} & \goth B_{12}\\ \goth B_{21}&\goth B_{22}
\endmatrix\right)$. Further, let $G_1$, $F_1$ be as earlier, and $G_2=F_2^{-1}= \mu_{22}
E \bar \mu_{22}^t$, with the same partition on blocks as earlier. As earlier we have:
$\im X_{21}$ is proportional to $F_{121}$, $\im Y_{22}$ is proportional to $F_{122}$, and
all are non-zero.

Further, we have $$B_2=\left (\matrix 0&X_{21} \\ 0&1\endmatrix\right)B_0, \ \ \ B_1=B_0
\left (\matrix 0&0 \\ \bar Y_{22}^t &1 \endmatrix\right) \eqno{(5.3.1.7)}$$

(4.2.7) and (5.3.1.7) imply that $\im X_{21} \goth B_{22}=\goth B_{22}\im Y_{22}^t=0$.
Substituting the expression for $\goth g$: $\goth g=\mu_{11}^{-1}B_0 \bar \mu_{22}^{t-1}$
in (3.1.1) we get formulas analogous to (5.3.1.3) - (5.3.1.5). The end of the proof is
similar to the one of the Case 1.

Proof for the remaining cases when for one or both $\s_1$, $\s_2$ satisfy $\{\det
\mu_1=0$ and $\det \mu_2\ne 0\}$ is analogous (we use Remark 4.1.4). $\square$

\medskip

{\bf Proposition 5.3.2.} For the case $g=3$, $\goth i=1$ we have: if $\s_1$, $\s_2\in
\goth D_*\cap S(3,1)^{open}$ and $\s_2$ satisfy the condition $A_2=C_2=0$ then $\s_1$,
$\s_2$ satisfy 4.2.6.

\medskip

{\bf Remark.} I think that this is true for any $\s_1$, $\s_2\in \goth D_*$, without the
restriction $A_2=C_2=0$.

\medskip

{\bf Proof.} We denote the scalar product of vectors by $<.,.>$. In order to simplify
notations, sometimes we does not distinguish between the row and column vectors. We
denote $A_1=\left (\matrix a_{1} \\ a_{2} \endmatrix\right)$, $C_1=\left (\matrix c_{1}
\\ c_{2} \endmatrix\right)$, $U_i=\left (\matrix u_{22i} & u_{23i} \\ u_{23i} & u_{33i}
\endmatrix\right)$. We assume the existence of $\goth g$ satisfying (4.2.6), transform
the corresponding formulas and show that these transformations are invertible. We shall
use notations $\goth g_i$ for the $i$-th column of $\goth g$. (4.1.2.5), Lemma 4.1.12,
condition $\s_1$, $\s_2\in \goth D_*\cap S(3,1)^{open}$ and (4.2.7) imply that

$$\left (\matrix -1&-c_1\rd&-c_2\rd\\ 1&-a_1&-a_2\endmatrix\right)
(\goth g_2|\goth g_3)\left (\matrix -1+u_{222}\rd\\ u_{232}\rd\endmatrix\right)= \left
(\matrix 0\\ 0\endmatrix\right)\eqno{(5.3.2.1)}$$

$$\left (\matrix a_1 + c_1 \rd & -1 - u_{221}\rd & -u_{231}\rd \endmatrix\right) \goth
g_1=0\eqno{(5.3.2.2)}$$

We denote $M=\left (\matrix a_1 + c_1 \rd & -1 - u_{221}\rd & -u_{231}\rd \\
-1&-c_1\rd&-c_2\rd\\ 1&-a_1&-a_2\endmatrix\right)$. A direct calculation shows that

$$\det \mu_{11}=\det \mu_{21}= 0 \iff \det M =0\eqno{(5.3.2.3)}$$

If $A=C=0$ then the same arguments as in the proof of 4.2.18 show us that the proposition
holds. If $A\ne 0$ or $C\ne 0$ then the second and the third lines of $M$ are linearly
independent and in this case (5.3.2.3) means that the first line of $M$ is a linear
combination of its second and third line.

Let us find a solution to the above equations. We denote the vector product of the lines
of $\left (\matrix -1&-c_1\rd&-c_2\rd\\ 1&-a_1&-a_2\endmatrix\right)$ by $\goth A$,
$\goth A=(c_1a_2-c_2a_1)\rd, -a_2 -c_2\rd, a_1+c_1\rd$. (5.3.2.1) implies that $(\goth
g_2|\goth g_3)\left (\matrix -1+u_{222}\rd\\ u_{232}\rd\endmatrix\right)$ is the
$\eta'\goth A$ for some coefficient $\eta'$. We have:

\medskip

$(-1+u_{222}\rd)\goth g_2+u_{232}\rd \goth g_3=\eta'\goth A$, or

$$\goth g_3=\eta\goth A+\gamma \goth g_2\eqno{(5.3.2.4)}$$

where $\gamma=\frac{1-u_{222}\rd}{u_{232}\rd}$ ($\left (\frac{-\d}p\right)=-1$ implies
$u_{232}\ne 0$) and $\eta$ is a coefficient.

Conditions $<E_{21}\goth g_2,\bar \goth g_1>=0$, $<E_{21}\goth g_3,\bar \goth g_1>=0$ and
(5.3.2.4) imply that $<E_{21}\goth A, \bar \goth g_1>=0$. Condition $<E_{21}\goth
g_3,\bar \goth g_2>=0$ and (5.3.2.4) imply that

$$\eta=-\frac{\gamma <E_{21}\goth g_2,\bar \goth g_2>}{<E_{21}\goth A, \bar \goth
g_2>}\eqno{(5.3.2.5)}$$

Further, Lemma 4.1.12 implies that $<E_{21}\goth A, \bar \goth A>=0$. Condition $\s_2\in
\goth D_*$ implies that $\gamma \bar \gamma =-1$. So, we get that (5.3.2.4), (5.3.2.5)
imply that always $<E_{21}\goth g_2,\bar \goth g_2>= <E_{21}\goth g_3,\bar \goth g_3>$.
Very well. Further, $\goth g_1$ must be orthogonal to both vectors $a_1 + c_1 \rd , -1 -
u_{221}\rd , -u_{231}\rd$ (condition (5.3.2.2)) and $\goth A$. But again Lemma 4.1.12
implies that they are proportional.

All above considerations are invertible, so in order to find $\goth g$ it is enough to
choose a vector $\goth g_1$ orthogonal to $\goth A$ such that $<E_{21}\goth g_1,\bar
\goth g_1>=1$, and a vector $\goth g_2$ orthogonal to $E_{21}\goth g_1$ such that
$<E_{21}\goth g_2,\bar \goth g_2>=1$ and such that $<E_{21}\goth a,\bar \goth g_2>\ne 0$.
It is clear that it is possible. Finally, we take $\goth g_3$ from 5.3.2.4, and all
conditions are satisfied. $\square$

\medskip

{\bf Theorem 5.3.3. } For the case $g=3$, $\goth i=2$ we have: there exists $\s_1\in
\goth D_*\cap S(3,2)^{open}$ such that for any $\s_2\in \goth D_*\cap S(3,2)^{open}$ we
have: $\s_1$, $\s_2$ satisfy 4.2.6.

\medskip

{\bf Proof.} We use the notations of the beginning of 5.1, we take $\s_1$ given by
(5.1.2.2), (in 5.1.2 this element is denoted by $\s_2$), and $\s_2$ is given by (5.1.1.1)
with any $a_1$, $a_2$, $\ve$. Formula for all $\goth g$ satisfying (4.2.8) is the
following: $$\goth g=k\left( \matrix \omega +\frac i2 (\ve a_2 -\rho) & \frac 12 (a_2
+\ve\rho) + i \ve (\omega - a_1) & 0
\\ \frac 12 (\ve a_2 +\rho) + i (\omega - a_1) & -\ve \omega +\frac i2 (-a_2 +\ve\rho) &
0 \\ 0&0&1\endmatrix \right)\eqno{(5.3.3.1)}$$ where $k$, $\omega$ are parameters in $\n
F_p^*$ such that $\rho=\sqrt{1-(a_1-2\omega)^2}\in \n F_p^*$ (clearly this condition can
be satisfied for many values of $\omega$).  $\square$

\medskip

{\bf Remark 5.3.4.} Method of finding of the formula (5.3.3.1).

\medskip

We denote by $l_i$ the third line of $\mu_{1i}$ ($i=1,2$). We can treat the free
3-dimensional module $\n F_{p^2}^3$ as the free 6-dimensional module $\n F_p^6$, and the
imaginary part of the Hermitian form $H(v_1, v_2)\overset{\opr}\to{=}v_1 \bar v_2^t$ as a
bilinear form on $\n F_p^6$. More exactly, there exists an isomorphism $\iota: \n
F_{p^2}^3 \to \n F_p^6$ such that $$\im H(v_1, v_2) = <\iota(v_1),
\iota(v_2)>\eqno{(5.3.4.1)}$$ where $<.,.>$ is the ordinary scalar product on $\n F_p^6$
(finding of formulas for $\iota$ is an elementary exercise).

Conditions (4.2.8) are corollaries of the following conditions: $$\im H(l_1 \goth gE,
m'_i)=0, \ \ \ \im H(m_i, l_2 \bar{\goth g}\bar E^t)=0\eqno{(5.3.4.2)}$$ where $m_i$
(resp. $m'_i$), $i=1,\dots,5$, are 5 non-equal lines of matrices $\mu_{11}, \mu_{21}$
(resp. $\mu_{12}, \mu_{22}$). (5.3.4.1), (5.3.4.2) imply that $\iota(l_1 \goth gE_{21}) =
k' v'$, $\iota(l_2 \bar \goth g\bar E_{21}^t) = k v$, where $v$, $v'$ are 6-vectors whose
coordinates are 5-minors of matrices formed by $\iota(m_i)$, $\iota(m'_i)$ respectively,
$i=1,\dots,5$, and $k,k'\in \n F_p^*$. An easy calculation shows that $v=\iota(l_1)$,
$v'=\iota(l_2)$, hence (4.2.8) is a corollary of the conditions

$$l_1 \goth gE_{21}=k'l_2 \eqno{(5.3.4.3)}$$

$$\goth gE_{21}\bar l_2^t=k \bar l_1^t\eqno{(5.3.4.4)}$$

Further, multiplying (5.3.4.3) by $\bar l_2^t$ from the right and (5.3.4.4) by $\bar l_1$
from the left, taking into consideration that $H(l_1,l_1)=H(l_2,l_2)=3$ (this follows
from the explicit formulas for them, see below), we get that if $p>3$ then $k=k'$.

Now we represent $\goth gE_{21}$ in the block form: $\goth g= \left( \matrix \goth
g_{11}&\goth g_{12}\\ \goth g_{21}&\goth g_{22}\endmatrix \right)$ with sizes of diagonal
blocks 2,1. We denote $l_i=(l_{i1},-1)$ ($i=1,2$), where $l_{i1}$ is a 2-vector formed by
2 first coordinates of $l_i$, and $-1$ is the third coordinate of both $l_1$, $l_2$. In
these notations formulas (5.3.4.3), (5.3.4.4) are equivalent to the following expressions
for $\goth g_{ij}$ in terms of $\goth g_{11}$:

$$\goth g_{12}=-k \bar l_{11}^t + \goth g_{11} \bar l_{21}^t$$

$$\goth g_{21}=-k l_{21} + l_{11} \goth g_{11}$$

$$\goth g_{22}=-k + l_{11} \goth g_{11} \bar l_{21}^t$$

Now we substitute these formulas in (3.1.1). We get:

$$\goth g_{11}\goth A\bar\goth g_{11}^t + \goth g_{11}\goth B+\bar \goth B^t\bar\goth
g_{11}^t=\goth C\eqno{(5.3.4.5)}$$

$$(\goth g_{11}\goth A\bar\goth g_{11}^t +\bar \goth B^t\bar\goth g_{11}^t -k^2 E_2)\bar
l_{11}^t=0 \eqno{(5.3.4.6)}$$

$$l_{11}\goth g_{11}\goth A\bar\goth g_{11}^t\bar l_{11}^t=-\lambda-k^2\eqno{(5.3.4.7)}$$

where $\goth A=E_2-\bar l_{21}^tl_{21}$, $\goth B=k\bar l_{21}^tl_{11}$, $\goth C=\lambda
E_2-\bar l_{11}^tl_{11}=\left( \matrix \lambda+k^2 & ik^2\\-ik^2&\lambda+k^2 \endmatrix
\right)$.

Now we make a change of variables in order to eliminate the linear terms in (5.3.4.5):
$\goth g_{11}=\gamma-\bar \goth B^t \goth A^{-1}$, and we substitute the square part of
(5.3.4.5) to (5.3.4.6), (5.3.4.7). We get:

$$\gamma\goth A\bar \gamma^t= \goth C+ \bar \goth B^t \goth A^{-1}\goth B$$

$$[(\lambda+k^2)E_2 - \gamma \goth B + \bar \goth B^t \goth A^{-1}\goth B]\bar
l_{11}^t=0$$

$$-l_{11} \gamma \goth B \bar l_{11}^t-l_{11}\bar \goth B^t \bar \gamma^t \bar l_{11}^t +
l_{11}(\goth C+ 2\bar \goth B^t \goth A^{-1}\goth B )\bar l_{11}^t + (\lambda+k^2) =0$$

Substituting values of all objects in the above equations, we get that (5.3.4.6),
(5.3.4.7) become respectively

$$-2k\gamma\bar l_{21}^t + (\lambda-3k^2)\bar l_{11}^t=0\eqno{(5.3.4.8)}$$

$$-2kl_{11}\gamma\bar l_{21}^t-2kl_{21}\gamma\bar l_{11}^t+3\lambda-11k^2=0
\eqno{(5.3.4.9)}$$

Substituting the expression for $\gamma\bar l_{21}^t$ from (5.3.4.8) to (5.3.4.9) we get
that (5.3.4.9) becomes $\lambda=k^2$ and (5.3.4.8) becomes

$$\gamma\bar l_{21}^t + k\bar l_{11}^t=0\eqno{(5.3.4.10)}$$

Now we write $\gamma=\left( \matrix \gamma_{11} & \gamma_{12} \\ \gamma_{21} &
\gamma_{22} \endmatrix \right)$ and $\gamma_{j1}=x_j+iy_j$. (5.3.4.10) gives us
expressions of $\gamma_{j2}$ in terms of $x_j, y_j$. Substituting in (3.1.1) we get a
system of equations with unknowns $x_j, y_j$. Formula (5.3.3.1) comes directly from these
solutions. 

\medskip

{\bf Conjecture 5.3.5.} For the case $g=3$, $\goth i=2$ we have: any $\s_1$, $\s_2\in
\goth D_1\cap S(3,2)^{open}$ satisfy 4.2.6.

\medskip

{\bf Idea of the proof.} We use notations of Subsection 4.1. According Remark 4.1.4, we
consider only the case $\det m_{11}\ne 0$, $\det m_{12}\ne 0$. We have $G_{22i}=0$
(Corollary 4.1.7) and hence $$\im X_{1i}G_{12i}=0\eqno{(5.3.5.1)}$$ (from 4.1.6.6).
Further, $F_{11i}G_{12i}=0$ (definition of $F$, $G$). Further, $$\goth
B_{22}=0\eqno{(5.3.5.2)}$$ --- like in the Proposition 4.2.16, case $\goth i=2$, hence
(4.2.15.5) implies that $$\im X_{12}\goth B_{21}^t=0\eqno{(5.3.5.3)}$$

\medskip

{\bf Lemma 5.3.5.4.} $\im X_{1i}\ne 0$ (recall that all equalities are in $\n F_p$).

\medskip

{\bf Proof.} If not, then eliminating $\re X_{2i}$ from ($4.1.6.1i$), ($4.1.6.2r$), we
get $\im X_{2i}(\Delta B_iC_i+A_i)=1$ --- a contradiction, because $\im X_{2i}(\Delta
B_iC_i+A_i)$ is a rank 1 matrix. $\square$

\medskip

So, (5.3.5.1) and (5.3.5.3) imply that $$\goth B_{21}=kG_{212}\eqno{(5.3.5.5)}$$

Analogs of (5.3.1.3) - (5.3.1.5) for the present case (taking into consideration
(5.3.5.2)) are

$$\goth B_{11} F_{112} \bar \goth B_{11}^t + \goth B_{12}F_{212}\bar\goth B_{11}^t+ \goth
B_{11}F_{122}\goth B_{12}^t + \goth B_{12} F_{222} \goth B_{12}^t = \lambda
G_{111}\eqno{(5.3.5.6)}$$

$$\goth B_{21}F_{112}\bar\goth B_{11}^t + \goth B_{21}F_{122} \goth B_{12}^t =\lambda
G_{211}\eqno{(5.3.5.7)}$$

$$\goth B_{21} F_{112} \goth B_{21}^t=\lambda G_{221}\eqno{(5.3.5.8)}$$

Since $G_{21i}F_{11i}=0$ and $G_{21i}F_{12i}=1$ (definition of $F$, $G$), we get from
(5.3.5.5) that (5.3.5.8) is always satisfied and (5.3.5.7) becomes $\goth
B_{12}=\frac{\lambda}k G_{121}$. So, we must only to substitute these values in (5.3.5.6)
and to find $\goth B_{11}$. I think that this is always possible.
\medskip

{\bf Subsection 5.4.} Non-coincidence of components of $T_*(V)$ for different types of
Hecke correspondences.  \nopagebreak

\medskip

Arguments like in 3.3.5 show that components of $T_p(V)$ do not coincide with components
of $T_{p,\goth i}(V)$ for any $\goth i$.

\medskip

For the case $\goth T_p=T_{p,0}$ the theory is similar to the one of Section 3. Analog of
(3.2.1) is $\s =\s(s)=\left( \matrix E_g & s\\0&p^2E_g \endmatrix \right)$ where $s\in
M_g(\n Z)^{symm}$ and its entries belong to a fixed system of residues modulo $p^2$,
$T=T(\s)$ is defined like in (3.2.2).

Components of $T_{p,0}(V)$ are defined over $K^{p^2}$, $K^p$ and $K^1$. For $g=3$ the
component $T_{p,0}(V)_{I(\s)}$ corresponding to $\s$ is defined over $K^1$ iff $\det T
\equiv 0 \mod p^2$. We call it the special component.

\medskip

{\bf Proposition 5.4.1.} For $g=3$ the special component of $T_{p,0}(V)$ does not
coincide with $V$ itself.

\medskip

{\bf Proof.} We use 3.3.5.1. Really, we shall prove that even the condition $G_V(\n
Q_p)\s^{-1} \cap G(\n Q_p) \ne \emptyset$ is not satisfied.

This condition is equivalent to the following one: there exists $\goth g = h + k \sqrt
{-\Delta}
\in G_V(\n Q_p), h, k \in M_3(\n Z_p)$, such that $$\alpha(\goth g)\s^{-1} \hbox{ has
integer entries} \eqno{(5.4.1.1)}$$ and $\ord_p(\det(\goth g))=3$.

We have: $\alpha (\goth g)\s^{-1} = \left (\matrix h & p^{-2}(-hs + kE_{21}) \\ -\Delta
E_{21}k & p^{-2}(\Delta E_{21}ks + E_{21}hE_{21})\endmatrix\right)$ hence (5.4.1.1) is equivalent
to the following congruences: $$k \equiv hsE_{21}  \mod  p^2 \eqno{(5.4.1.3)}$$ $$\Delta
ks + hE_{21} \equiv 0  \mod  p^2 \eqno{(5.4.1.4)}$$

Substituting $k$ from (5.4.1.3) to (5.4.1.4) we get: $$hT \equiv 0 \mod  p^2
\eqno{(5.4.1.5)}$$

{\bf (5.4.1.6)} To simplify calculations, we consider a particular value of $\s$, namely
$\s=\s(s)$, where $s=W$ is from the beginning of the proof of 3.4.3, $w_1,w_2\in \n Z$
satisfy $w_1^2-w_2^2\equiv \frac{-1}{\d} \mod p^4$. We have $T\equiv \left (\matrix 1 & 0
& 0 \\ 0 & 0 & 0\\0 & 0 & 0\endmatrix \right) \mod p^4$, hence (5.4.1.5) is equivalent to

$$ h \equiv \left (\matrix 0 & * & * \\ 0 & * & *\\0 & * & *\endmatrix \right)\mod 
p^2\eqno {(5.4.1.7)}$$

According (5.4.1.3), $k$ must satisfy $k = hsE_{21} + p^2A$ with $A \in M_3(\n Z_p)$.
Condition (3.1.1)

is the following: $$ hTh^t + p^2\{\sqrt{-\Delta}(-hE_{21}A^t + AE_{21}h^t) + \Delta (Ash^t
+ hsA^t)\} + p^4\Delta AE_{21}A = \lambda E_{21}\eqno{(5.4.1.8)}$$

(5.4.1.7) implies $hTh^t \equiv  0 \mod p^4$. Dividing (5.4.1.8) by $p^2$ and considering
it modulo $p$, we get equations in $\n F_p$ (tildes are omitted):$$ hE_{21}A^t =
(hE_{21}A^t)^t;\eqno{(5.4.1.9)}$$

$$ Ash^t + hsA^t = \lambda' E_{21},\eqno{(5.4.1.10)}$$ where $\lambda'=\lambda/p^2 \ne 0$
in $\n F_p$.

We denote by $h'$ the $3\times 2$-matrix  obtained by rejecting the first (zero) column
of $h$. If the rank of $h'$ is $\le 1$ then rank of $ Ash^t + hsA^t\le 2$: a
contradiction. If the rank of $h'$ is 2 then the set of $X$ such that $Xh^t = hX^t$ is
the set $(*|h'\alpha)$, where * means an arbitrary $3\times 1$-matrix, $|$ means the
operation of gluing of matrices and $\alpha$ runs over the set of symmetric $2 \times
2$-matrices. (5.4.1.9) implies an equality  $AE_{21}= (*|h'\alpha)$ which gives us $A =
(*|h'\alpha E_{11})$.

For this $A$ we have $Ash^t = (*|h'\alpha E_{11})sh^t = h'\alpha E_{11}W_2h'{^t}$ and
$hsA^t = h'W_2E_{11}\alpha h'{^t}$ (here $W_2 = \left (\matrix w_1 & w_2\\w_2 & w_1 \endmatrix
\right))$, hence the left hand side of (5.4.1.10) is  $h'(\alpha E_{11}W_2 +
W_2E_{11}\alpha)h'{^t}$

which is of rank at most 2 - a contradiction to (5.4.1.10). $\square$

\medskip

{\bf Proposition 5.4.2.} The special component of $T_{p,0}(V)$ does not coincide with (a)
one of general components of $T_{p,1}(V)$; (b) the general component of $T_{p,2}(V)$; (c)
one of the intermediate components of $T_{p,2}(V)$.

\medskip

{\bf Proof.} As usually, we denote $\goth g=h+k\rd$. Obviously, Theorem 2.16 can be
generalized to the case of different $S(g,\goth i)$ as follows:

\medskip

{\bf (5.4.2.0)} Let $\s_1\in S(g,\goth i_1)$, $\s_2\in S(g,\goth i_2)$. Then $T_{p,\goth
i_1}(V)_{D(\s_1)}=T_{p,\goth i_2}(V)_{D(\s_2)}$ iff (2.17) holds.

\medskip

As earlier we can consider only the $p$-component of adeles, i.e. (2.17) means that there
exists $\goth g\in G_V$ such that $\s_1 \alpha(\goth g) \s_2^{-1}\in G_X(\n Z_p)$. If
$\s_1$, $\s_2$ are fixed then $4g^2$ entries of $\s_1 \alpha(\goth g) \s_2^{-1}$ are
linear forms of $2g^2$ variables $h_{ij}$, $k_{ij}$. We denote by
$L_{\lambda,\mu}=L_{\lambda,\mu}(h_{ij}, k_{ij})$ the linear form that is the
($\lambda,\mu$)-th entry of $\s_1 \alpha(\goth g) \s_2^{-1}$. The set of variables
$h_{ij}$, $k_{ij}$ is $\n Q_p^{2g^2}$, and the condition that all entries of $\s_1
\alpha(\goth g) \s_2^{-1}$ are integer ($\iff$ all $L_{\lambda,\mu}(h_{ij}, k_{ij})$ are
integer) give us a $\n Z_p$-linear subspace in $\n Q_p^{2g^2}$. We denote this subspace
by $R$; we shall see that $R$ is a lattice. $T_{p,\goth i_1}(V)_{D(\s_1)}=T_{p,\goth
i_2}(V)_{D(\s_2)}$ iff there exists ($h_{ij}$, $k_{ij}$) such that $\ord_p \det \goth
g=0$.

Since we are proving the proposition not for any $g$ but only for $g=3$, I recommend to
the reader to get the $4g^2\times2g^2$-matrix of coefficients of $L_{\lambda,\mu}$ by
means of computer calculations. We denote this matrix by $M$. To prove non-equivalence,
we shall consider submatrices of $M$; these submatrices will define us upper bounds of
$R$.

For any component under consideration we consider a fixed value of $\s=\s_1$ or $\s_2$
such that $T_{p,\goth i}(V)_{I(\s)}$ is equal to this component. For the special
component of $T_{p,0}(V)$ we take $\s_1$ from 5.4.1.6.

\medskip

{\bf Case A.} A general component of $T_{p,1}(V)$. We take $\s_2$ from (4.1.1) with
$A=U=C=0$ (for other general component of $T_{p,1}(V)$ the proof is analogous).
Immediately from the condition that all entries of $\s_1 \alpha(\goth g) \s_2^{-1}$ are
integer we get the following table of inferiour bounds of $p$-orders of $h_{ij}$,
$k_{ij}$ ($p$-order is $\ge$ of the number in the table):

\settabs 20 \columns

\medskip

\+ &&&$h$&&&&$k$&&&&&&&&&&(5.4.2.1)\cr
\nopagebreak
\medskip
\+&&1&2&3&&1&2&3\cr
\nopagebreak
\medskip
\+1&&1&0&0&&1&2&2\cr
\nopagebreak
\+2&&$-1$&0&0&&$-1$&0&0\cr
\nopagebreak
\+3&&$-1$&0&0&&$-1$&0&0\cr

\medskip

(2,1), (2,4), (3,1), (3,4)-entries of $\s_1 \alpha(\goth g) \s_2^{-1}$ are linear forms
of only 4 variables $h_{i1}, k_{i1}$, $i=2,3$. The corresponding $4\times 4$-submatrix of
$M$ (i.e. the $4\times 4$-matrix of coefficients of these forms) is $\frac1p M'$ where
entries of $M'$ are in $\n Z_p$ and $\det M'\ne 0$. This means that $\ord_p(h_{i1},
k_{i1})\ge 1$. This implies that $\ord_p(\det \goth g)\ge 1$, i.e. there is no
equivalence in this case.

\medskip

{\bf Case B.} The general component of $T_{p,2}(V)$. We take $\s_2$ from (4.1.1) with
$A=U=C=0$. Analog of (5.4.2.1) for this case is the following:

\medskip

\+ &&&$h$&&&&$k$&&&&&&&&&&(5.4.2.2)\cr

\medskip

\+&&1&2&3&&1&2&3\cr

\medskip

\+1&&0&1&1&&2&1&1\cr

\+2&&0&$-1$&$-1$&&$-2$&$-1$&$-1$\cr

\+3&&0&$-1$&$-1$&&$-2$&$-1$&$-1$\cr

\medskip


(2,2), (2,5), (3,2), (3,5)-entries of $\s_1 \alpha(\goth g) \s_2^{-1}$ are linear forms
of only 4 variables $h_{i2}, k_{i2}$, $i=2,3$, and the $4\times 4$-matrix of coefficients
of these forms has $\det\ne 0\mod p$, so $\ord_p(h_{i2}, k_{i2})\ge 1$ 

The same is true for variables $h_{i3}, k_{i3}$, $i=2,3$ ((2,3), (2,6), (3,3),
(3,6)-entries of $M$). 

Finally, treating the (2,1), (2,4), (3,1), (3,4)-entries of $M$ and variables $h_{i1},
k_{i1}$, $i=2,3$, we get that again the $4\times 4$-matrix of coefficients of these forms
has $\det\ne 0\mod p$. 

Using table 5.4.2.2 we get the non-equivalence for this case.

\medskip

{\bf Case C.} The intermediate component of $T_{p,2}(V)$. We take $\s_2$ from (4.1.1)
with $A=(1 \ \ 0), U=C=0$). Analog of (5.4.2.1) for this case is the following:

\medskip

\+ &&&$h$&&&&$k$&&&&&&&&&&(5.4.2.3)\cr

\medskip

\+&&1&2&3&&1&2&3\cr

\medskip

\+1&&0&0&1&&1&1&1\cr

\+2&&$-1$&$-1$&$-1$&&$-2$&$-2$&$-1$\cr

\+3&&$-1$&$-1$&$-1$&&$-2$&$-2$&$-1$\cr

\medskip


12 entries of $\s_1 \alpha(\goth g) \s_2^{-1}$ (namely, (2,*), (3,*)-entries) are linear
forms of only 12 variables (all except $h_{1*}, k_{1*}$) and the determinant of the
corresponding $12\times12$-submatrix of $M$ is $\ne 0\mod p$. This implies the
non-equivalence for this case. $\square$

\medskip


{\bf Remark.} I think that for other components also there is no coincidence.

\medskip

\medskip

{\bf References} \nopagebreak

\medskip

\medskip

[D71] Deligne P. Travaux de Shimura. Lect. Notes in Math., 1971,
v.244, p. 123 - 165. Seminaire Bourbaki 1970/71, Expos\'e 389.

\medskip

[GZ86] Gross B.H., Zagier D.B., Heegner points and derivatives of $L$-
series. Invent. Math., 1986, v.84, N.2, p. 225 - 320.

\medskip

[K89] Kolyvagin V.A., Finiteness of $E(\n Q )$ and Sh$(E,\n Q )$
for a subclass of  Weil curves. Math. USSR Izvestiya, 1989, v. 32,
No. 3, p. 523 - 541

\medskip

[L01] D. Logachev. Action of Hecke correspondences on Heegner curves on a Siegel
threefold. J. of Algebra, 2001, v. 236, N. 1, p. 307 -- 348.

\medskip

[L04.1] D. Logachev. Relations between conjectural eigenvalues of Hecke operators on
submotives of Siegel varieties. http://arxiv.org/ps/math.AG/0405442

\medskip

[L04.2] D. Logachev. Reduction of a problem of finiteness of Tate-Shafarevich group to a
result of Zagier type. http://arxiv.org/ps/math.AG/0411346

\medskip

[Z85] Zagier, Don,  Modular points, modular curves, modular surfaces and modular forms.
Workshop Bonn 1984 (Bonn, 1984), 225 - 248, Lecture Notes in Math., 1111, Springer,
Berlin, 1985.

\medskip

E-mail: logachev\@ usb.ve

\enddocument